\documentclass[a4paper,twoside,11pt]{article}

\usepackage{a4wide, amsmath, amssymb, mathtools}
\usepackage[utf8]{inputenc}
\usepackage{amsthm}
\usepackage[english]{babel}
\usepackage{hyperref}
\usepackage{authblk}
\numberwithin{equation}{section}

\newcommand{\Z}{\mathbb{Z}}
\newcommand{\Q}{\mathbb{Q}}
\newcommand{\R}{\mathbb{R}}
\newcommand\FF{\mathbb{F}}
\newcommand{\PPP}{\mathbb{P}}
\newcommand{\AAA}{\mathbb{A}}
\newcommand{\Adele}{\mathbf{A}}
\newcommand{\leg}[2]{\left(\frac{#1}{#2}\right)}

\DeclareMathOperator{\Vol}{Vol}
\DeclareMathOperator{\GL}{GL}
\DeclareMathOperator{\rk}{rk}
\DeclareMathOperator{\Msym}{Sym}
\DeclareMathOperator{\Br}{Br}
\DeclareMathOperator{\Pic}{Pic}
\DeclareMathOperator{\HH}{H}
\DeclareMathOperator{\PBM}{PBM}
\DeclareMathOperator{\adm}{adm}
\DeclareMathOperator{\inv}{inv}
\DeclareMathOperator{\Eff}{Eff}
\DeclareMathOperator{\Nef}{Nef}
\DeclareMathOperator{\Mov}{Mov}
\DeclareMathOperator{\Proj}{Proj}
\DeclareMathOperator{\Supp}{Supp}
\DeclareMathOperator{\Tr}{Tr}

\newtheorem{lemma}{Lemma}[section]
\newtheorem{theorem}[lemma]{Theorem}
\newtheorem{prop}[lemma]{Proposition}
\newtheorem{corollary}[lemma]{Corollary}
\newtheorem{conjecture}[lemma]{Conjecture}
\newtheorem{question}[lemma]{Question}
\theoremstyle{definition}
\newtheorem{remark}[lemma]{Remark}
\newtheorem{definition}[lemma]{Definition}

\title{\vspace{-\baselineskip}\sffamily\bfseries On the leading constant in Manin's conjecture}
\author[1]{Stephanie \hspace{-0.15cm} Chan\thanks{Department of Mathematics, University College London, Gower Street, London, WC1E 6BT, UK, stephanie.chan@ucl.ac.uk}}
\author[2]{Peter \hspace{-0.15cm} Koymans\thanks{Mathematisch Instituut, Universiteit Utrecht, Postbus 80.010, 3508 TA Utrecht, The Netherlands, p.h.koymans@uu.nl}}
\author[3]{Nick \hspace{-0.15cm} Rome\thanks{Max Planck Institute for Mathematics, Vivatsgasse 7, 53111 Bonn, Germany, nrome@mpim-bonn.mpg.de}}
\author[4]{Tim \hspace{-0.15cm} Santens\thanks{University of Cambridge, DPMMS, Centre for Mathematical Sciences, Wilberforce
Road, Cambridge, CB3 0WB, UK, ts996@cam.ac.uk}}
\affil[1]{University College London}
\affil[2]{Utrecht University}
\affil[3]{MPIM Bonn}
\affil[4]{University of Cambridge}
\date{\today}

\begin{document}
\maketitle

\begin{abstract}
We disprove a conjecture of Batyrev and Tschinkel \cite{BT98} concerning counting rational points of bounded height. Our results strongly suggest that the transcendental Brauer group plays a role in the leading constant in problems of the same nature such as the Manin--Peyre conjecture.
\end{abstract}

\section{Introduction}
The Manin--Peyre conjecture \cite{Peyre} is a major problem in Diophantine geometry which makes a precise prediction for the asymptotic growth of the number of points of bounded anticanonical height on a smooth Fano variety. Batyrev--Tschinkel \cite{BT98} made a generalised version of the conjecture which extends the problem to singular varieties. In particular, given a Fano variety $V$ which is not too singular, they predicted the following asymptotic for the associated counting problem. Given a variety $X$ over $\Q$, we write $X_{\overline{\Q}}$ (or $(X)_{\overline{\Q}}$) for its base change to an algebraic closure $\overline{\Q}$.

\begin{conjecture}
\label{conj:BT}
Let $V$ be a terminal Fano, $\Q$-Gorenstein variety with $V(\Q)$ not thin such that the smooth locus $U = V^{\mathrm{sm}}$ is geometrically simply connected. Then, for any anticanonical height $H_{-K_V}$, there exists a thin subset $\Omega \subset V(\Q)$ such that the following asymptotic formula holds:
\[
\#\{ x \in V(\Q)\setminus \Omega : H_{-K_V}(x) \leq B\}
\sim
C_{\textup{BT}} B (\log B)^{\rk \Pic U - 1},
\]
where
\[
C_{\textup{BT}} = \frac{\alpha(U)}{(\rk \Pic U - 1)!} \#\HH^1(\Q, \Pic U_{\overline{\Q}})\tau(\overline{U(\Q)}),
\]
and where $\alpha(U)$ and $\tau(\overline{U(\Q)})$ are analogues of Peyre's effective cone constant and Tamagawa number (denoted $\gamma_{-K_V}(V)$ and $\tau_{-K_V}(V)$ in \cite[Definitions 2.3.16, 3.3.10]{BT98}), respectively.
\end{conjecture}

\begin{remark}
This statement differs from that made by Batyrev--Tschinkel in \cite{BT98} in a number of ways. Firstly, they stated this prediction as merely an expectation rather than a full conjecture. Secondly, their conjecture is formulated in the more general setting of heights associated to other line bundles. We have chosen to stick with the anticanonical line bundle for simplicity of presentation. Further, we have updated the language to reflect the way in which the community currently understands such problems (i.e.~to incorporate the notion of an exceptional thin set). 
The assumption that $V$ is terminal Fano is also slightly less general than the generality considered by Batyrev--Tschinkel. The precise relation between the above conjecture and \cite[\S 3.4, Step 4]{BT98} is spelled out in more detail in Section~\ref{ss:BTconj}.

Finally, we note that the Tamagawa measure is being taken of the closure $\overline{U(\Q)} \subseteq \overline{V(\Q)} \cap U(\Adele_\Q)$ in the adelic topology. The latter quantity is what appears originally in Batyrev--Tschinkel and its Tamagawa measure can be strictly larger. However, we demonstrate in Remark \ref{rem:tau} that formulating the conjecture instead with the Tamagawa measure of $\overline{V(\Q)} \cap U(\Adele_\Q)$ leads to a contradiction with Batyrev--Tschinkel's previous work on toric varieties.
\end{remark} 

The purpose of the present paper is to produce an explicit counterexample to Conjecture \ref{conj:BT}. Let $k\geq 2$ be an integer and put $n=4k$. Consider the product of conic bundles defined by the system of equations
\[
\begin{array}{cccccccl}
	t_1 t_2 X_1^2 & -&t_5 \cdots t_n Y_1^2 &+&t_3 t_4 Z_1^2 &= &0, \\
	t_5 t_6 X_2^2 & -&t_1 \cdots t_4 \cdot t_9 \cdots t_n Y_2^2 \quad &+&t_7 t_8 Z_2^2& = &0, \\	
	\vdots &&\vdots &&\vdots&& \\
	t_{n - 3} t_{n - 2} X_k^2 \quad &-&t_1 \cdots t_{n - 4} Y_k^2 &+&t_{n - 1} t_n Z_k^2 &=& 0.
	\end{array}
\]
Let $V_k$ denote the resulting projective algebraic variety (see Section~\ref{sec:construction_variety} for the precise definition). Let $Z$ denote the subset defined by at least two of the $Y_i$ vanishing.

\begin{theorem}
\label{thm:intromain}
There exists $k_0 \in \Z_{>0}$ such that for all $k \geq k_0$ and all thin subsets $Z(\Q) \subseteq \Omega \subset V_k(\Q)$, we have, for $B\gg_{k, \Omega} 1$, the inequality
\[
\#\{ x \in V_k(\Q)\setminus \Omega : H_{-K_{V_k}}(x) \leq B\} > 2C_{\textup{BT}} B (\log B)^{\rk \Pic V_k^{\mathrm{sm}} - 1},
\] 
where $C_{\textup{BT}}$ is defined as in Conjecture \ref{conj:BT}.
\end{theorem}

Unfortunately, the variety $V_k$ is not Fano for $k \geq 3$ and thus Theorem \ref{thm:intromain} does not directly contradict Conjecture \ref{conj:BT}. Moreover, the fact that the anticanonical divisor is big but not ample means that the Northcott property can fail. This is the reason why we must remove the points in $Z$. However, $V_k$ is a Mori dream space. Thus, by appropriating techniques from the Minimal Model Program in birational geometry, we are able to conclude the existence of a Fano variety contradicting Conjecture \ref{conj:BT}.

\begin{corollary}
\label{cor:Fano}
There exists a terminal, $\Q$-Gorenstein Fano variety $V$ with geometrically simply connected smooth locus $V^{\mathrm{sm}}$ such that for all thin subsets $\Omega \subset V(\Q)$ and sufficiently large $B$, we have
\[
\#\{ x \in V(\Q)\setminus \Omega : H_{-K_V}(x) \leq B\} > 2C_{\textup{BT}} B (\log B)^{\rk \Pic V^{\mathrm{sm}} - 1}.
\]
\end{corollary}

To the authors' knowledge, the deduction of this corollary (performed in Sections~\ref{ss:Mori}--\ref{ss:BTconj}) is a completely novel argument in the context of counting rational points. One might object that we have produced a counterexample not to Batyrev--Tschinkel's original conjecture but rather to our reformulation. However, if one replaces $\tau(\overline{U(\Q)})$ by $\tau(\overline{V(\Q)} \cap U(\Adele_\Q))$ in the definition of $C_{\textup{BT}}$ then our varieties still provide a counterexample (cf.\ Theorems~\ref{theorem:BTconstant} and~\ref{thm:maintermcomparison}).
\begin{remark}
Recently, Lehmann--Tanimoto~\cite{almost} have demonstrated that even the exponent of $B$ predicted in the Manin conjecture can be false for varieties which are only ``almost Fano''. In particular, their example has an anticanonical divisor which is big but not nef. However, this example is not a Mori dream space and so is in a sense even further from being Fano than the varieties considered here.
\end{remark}

\subsection{Counting in the presence of transcendental Brauer group elements}
Peyre~\cite{Peyre} gave the first systematic treatment of the leading constant in Manin's conjecture. His original formulation of the predicted constant only contained two factors: the cone constant $\alpha$ and the Tamagawa number $\tau$. However, Peyre worked under the additional assumption that $\Pic X_{\overline{\Q}}$ was a permutation module for the action of the absolute Galois group of $\Q$ (which in particular implies that $\HH^1(\Q, \Pic X_{\overline{\Q}}) = 0$). Batyrev--Tschinkel showed that without this condition when considering toric varieties one needed to include an extra factor $\beta = \# \HH^1(\Q, \Pic X_{\overline{\Q}})=\#\left(\Br_1 X/\Br \Q\right)$, the size of the algebraic part of the Brauer group modulo constants. Moreover, Salberger~\cite{Salberger} demonstrated how such a factor arises naturally from the perspective of universal torsors. 

The transcendental part of the Brauer group is trivial for smooth proper rational varieties. 
When considering singular varieties, as outlined in the discussion in Section~\ref{ss:BTconj}, the relevant quantity is the Brauer group of the smooth locus.
For any of the singular varieties or non-rational Fano varieties on which we have been able to count points, there have been no transcendental Brauer classes in the Brauer groups of the smooth locus.
Therefore it has been reasonably natural to express the cohomological constant $\beta$ as the size of the algebraic part of the Brauer group. However, restricting to this algebraic part seems unjustified in the setting of varieties with transcendental Brauer group elements. Our suggestion is that the underlying mechanism of the counterexample in Theorem \ref{thm:intromain} is the size of the transcendental part of the Brauer group (of the smooth locus) appearing in the count for rational points. Indeed, for our example variety $V_k$ there are no non-constant algebraic Brauer classes. However, the transcendental part has size $ 2^{\frac{k(k-1)}{2}}$ and the point count is large by a factor of size at least $2^{\frac{k^2}{2} - O(k\log k)}$. Inspired by this observation, we pose the following variation on Conjecture \ref{conj:BT} as a question.

\begin{question}
For any terminal, $\Q$-Gorenstein Fano variety $V$ defined over $\Q$ with $V(\Q)$ not thin, let $U=V^{\mathrm{sm}}$ denote the smooth locus. Assume that $U_{\overline{\Q}}$ is simply connected and let $H_{-K_V}$ be an anticanonical height function. Does there exist a thin subset $\Omega \subset V(\Q)$ such that
\[
\#\{ x \in V(\Q)\setminus \Omega : H_{-K_V}(x) \leq B\}
\sim 
 \frac{\alpha(U)}{ (\rk\Pic U - 1)!} \#\left(\frac{\Br U}{\Br \Q}\right)\tau(\overline{U(\Q)})
B (\log B)^{\rk\Pic U - 1}?
\]
\end{question}

\begin{remark}
We note that since $V$ is assumed to be terminal Fano and $U_{\overline{\Q}}$ is simply connected, $\Pic U_{\overline{\Q}}$ is finitely generated and torsion free which, combined with the finiteness of $\Br U_{\overline{\Q}}$ and the Hochschild--Serre spectral sequence, guarantees the finiteness of the quotient $\#\left(\frac{\Br U}{\Br \Q}\right)$.
\end{remark}

Further evidence for the need to introduce the transcendental part of the Brauer group into such predictions is provided by the related problem of counting number fields. Loughran--Santens~\cite{LS} recently produced a conjectural leading constant in this problem inspired by Peyre's work and here one explicitly needs to include the size of the full Brauer group. This also comports with the function field analogue of Manin's conjecture, where the presence of the transcendental Brauer group is related to the triviality of the Griffiths group and the integral Hodge conjecture for curve classes, see \cite[Remark 4.13]{LehmannRiedlTanimoto}.

The major technical achievement of our paper is the ability to produce an explicit lower bound for the number of rational points in the presence of a transcendental Brauer group element. 
As far as the authors are aware, there doesn't exist a single prior example of sharp asymptotic upper or lower bounds for the rational points of bounded height on a Fano variety over a number field with non-trivial geometric Brauer group. In the function field setting, Okamura~\cite{doublesolid} established the geometric Manin conjecture for a general Artin--Mumford double solid providing the first example in that regime where transcendental Brauer group elements arise.
We suggest that producing such an example over $\Q$ is a highly pertinent open problem for the community.

\begin{question}
Can one produce a terminal, $\Q$-Gorenstein Fano variety $V$ defined over $\Q$ such that $\Br V^{\mathrm{sm}}/\Br_1 V^{\mathrm{sm}} \neq 0$, $(V^{\mathrm{sm}})_{\overline{\Q}}$ is simply connected, and such that one can prove an asymptotic formula for the points of bounded anticanonical height on $V$?
\end{question}

\subsection*{Layout of the paper}
We advise the reader that the technical heart of our argument is contained in Sections \ref{sCounting}, \ref{s:geometry}, and \ref{s:thinsets}. Section \ref{sCounting} is devoted to the main analytic counting argument, Section \ref{s:geometry} to interpreting Conjecture \ref{conj:BT} for our particular variety and the establishment of Corollary \ref{cor:Fano}, and Section \ref{s:thinsets} contains the ultimate deduction of our main theorem, including an argument showing that the contribution from accumulating thin sets is negligible. The main results of Section \ref{s:thinsets} are Theorem \ref{thm:maintermcomparison} and Theorem \ref{thm:thinsets}.

In Section \ref{s:construction}, we establish a number of geometric facts about our variety. Section \ref{s:tamagawa} concerns definitions and certain transformations of the Tamagawa number in generality. Section \ref{sTernary} contains one of the major analytic inputs of the paper, a variation on a result of Sofos which allows us to count points (uniformly) on conics under congruence conditions. The local densities arising in the leading constant in this result are analysed in Section \ref{s:densities}, and bounds for the size of the predicted local densities on the total space are developed in Section \ref{s:tamagawabounds}.

\subsection*{Acknowledgements}
The authors would like to thank Carlo Pagano for his contributions and insights towards the beginning of this project.
The second author gratefully acknowledges the support of the Dutch Research Council (NWO) through the Veni grant ``New methods in arithmetic statistics''. The third author was funded by FWF project ESP 441-NBL. The last author was supported by the Herchel Smith Fund.
A part of this work was completed whilst all the authors were attending the programme Universal Statistics in Number Theory at the Centre de Recherches Math{\'e}matiques in Montr{\'e}al. The institute's hospitality, excellent working environment and support are gratefully acknowledged.

\subsection*{AI usage disclosure}
In contrast to the recent trend, our counterexample was arrived at without influence from any Large Language Model. The paper was written without any such intervention until being audited for mistakes and presentational issues by all three major LLMs, one of which was made available to the fourth author by Cambridge University. All mathematical statements, arguments, and all citations were made by the authors. 

\section{Construction of a variety}
\label{s:construction}
\subsection{A toric variety}
Recall that given a vector bundle $\mathcal{V}$ on a scheme $X$ we let $\PPP_X(\mathcal{V}) = \Proj_X \bigoplus_{p = 0}^{\infty} \Msym^p \mathcal{V}^{\vee}$ be the \emph{projective bundle}\footnote{Note that this is dual to the convention used in the Stacks project \cite[\href{https://stacks.math.columbia.edu/tag/01OA}{Tag 01OA}]{Stacks}.} whose points correspond to lines in $\mathcal{V}$. 

Let $n = 4k$. Denote $[k] \coloneqq \{1, \dots, k\}$. For each $i\in[k]$, we let $\mathcal{V}_i \coloneqq \mathcal{O}_{\PPP^{n - 1}}(2k - 3) \oplus \mathcal{O}_{\PPP^{n - 1}} \oplus \mathcal{O}_{\PPP^{n - 1}}(2k - 3)$ and we define the fibre product
\[
	\mathcal{T}_k \coloneqq \PPP_{\PPP^{n-1}}(\mathcal{V}_1) \times_{\PPP^{n-1}} \cdots \times_{\PPP^{n-1}} \PPP_{\PPP^{n-1}}(\mathcal{V}_k).
\] 
The projection $\mathcal{T}_k \to \PPP^{n-1}$ will be denoted $\pi$.

We will use the coordinates $[t_1 : \cdots: t_n]$ on $\PPP^{n - 1}$. To write down equations on the projective bundles we will use the coordinates
\[
\mathcal{V}_i^{\vee} = \mathcal{O}_{\PPP^{n - 1}}(-2k + 3)X_i \oplus \mathcal{O}_{\PPP^{n - 1}} Y_i \oplus \mathcal{O}_{\PPP^{n - 1}}(-2k + 3) Z_i.
\] 
Let $\mathcal{O}_{\PPP(\mathcal{V}_i)}(1)$ (here, and from now on, we will drop the $\PPP^{n - 1}$ subscript in our notation for line bundles) be the tautological bundle on $\PPP_{\PPP^{n - 1}}(\mathcal{V}_i)$ and by abuse of notation we will also use the same symbol for its pullback to $\mathcal{T}_k$. In particular, we have $\pi_*\mathcal{O}_{\PPP(\mathcal{V}_i)}(1) = \mathcal{V}_i^{\vee}$ by definition. Let $L_i$ be its divisor class. Let $H$ be the hyperplane class on $\PPP^{n-1}$.

\begin{lemma}
\label{lem:canonical_toric}
	We have $K_{\mathcal{T}_k} = -4k \pi^*H - \sum_{i = 1}^k (3L_i + (4k - 6) \pi^*H).$
\end{lemma}

\begin{proof}
	This follows from iteratively applying the formula for the canonical line bundle on a projective bundle coming from the relative Euler sequence.
\end{proof}

\begin{remark}
\label{rem:Cox}
	One way to think about $\mathcal{T}_k$ is as a toric variety. Indeed, the vector bundles $\mathcal{V}_i$ are all torus equivariant so the total space is a toric variety.
	 One can understand this toric variety concretely by its Cox coordinates, see \cite[Chapter 5]{CoxLittleSchenck} for an introduction. These are
	 \[
	 t_1, \dots, t_n, X_1, Y_1, Z_1, \dots, X_k, Y_k, Z_k
	 \] 
	 together with an action of $(\lambda_0, \dots, \lambda_k) \in \mathbb{G}_m^{k + 1}$ given by
	 \[
	 	\left((t_i)_{i \in [n]},\ (X_i,\ Y_i,\ Z_i)_{i \in [k]}\right) \to \left((\lambda_0 t_i)_{i \in [n]},\ \left((\lambda_i \lambda_0^{2k - 3})X_i,\ \lambda_i Y_i,\ \lambda_i \lambda_0^{2k - 3}Z_i\right)_{i \in [k]}\right).
	 \]	 
	 The quotient by this action gives a morphism
	 \[
	 \varpi: (\AAA^n \setminus \{0\}) \times (\AAA^3 \setminus \{0\})^k \to \mathcal{T}_k.
	 \]
	 Since $(\lambda_1, \dots, \lambda_k) \in \mathbb{G}_m^k$ acts by standard scalar multiplication on each $\AAA^3 \setminus \{0\}$ factor, $\varpi$ factors through a $\mathbb{G}_m$-torsor (the quotient by $\lambda_0 \in \mathbb{G}_m$), which by abuse of notation we also denote by $\varpi$:
	 \[
	 (\AAA^n \setminus \{0\}) \times (\AAA^3 \setminus \{0\})^k \to (\AAA^n \setminus \{0\}) \times (\PPP^2)^k \xrightarrow{\varpi} \mathcal{T}_k.
	 \]
	 A rational point on $\mathcal{T}_k$ is thus given by a list of $7k$ elements
	 \[
	 t_1, \dots, t_n, X_1, Y_1, Z_1, \dots, X_k, Y_k, Z_k \in \Q
	 \]
	 such that not all $t_i$ are zero and for each $i$ not all of $X_i, Y_i, Z_i$ are zero (modulo the action of $(\Q^{\times})^{k + 1}$).
\end{remark}

\subsection{\texorpdfstring{The variety $V_k$}{The variety Vk}} \label{sec:construction_variety}
We now want to describe a subvariety of $\mathcal{T}_k$. To do this we need equations, i.e.~global sections of line bundles. Note that for $i\in[k]$ and $a \in \Z_{\geq 0},\ b \in \Z_{> 0}$ we have the following inclusion by the push-pull formula:
\begin{equation}
\label{eq:global_sections_Tk}
	\begin{split}
	&\HH^0(\PPP^{n-1}, \mathcal{O}_{\PPP^{n - 1}}(a + b(-2k + 3)))X_i^b \\
	&\oplus \HH^0(\PPP^{n-1}, \mathcal{O}_{\PPP^{n - 1}}(a)) Y_i^b \\
	&\oplus \HH^0(\PPP^{n-1}, \mathcal{O}_{\PPP^{n - 1}}(a + b(-2k + 3)))Z_i^b \\
	&\subseteq \HH^0(\mathcal{T}_k, \pi^{*}\mathcal{O}_{\PPP^{n-1}}(a) \otimes \mathcal{O}_{\PPP(\mathcal{V}_i)}(b)) = \HH^0(\PPP^{n-1}, \mathcal{O}_{\PPP^{n - 1}}(a) \otimes \Msym^b \mathcal{V}_i^{\vee}).
\end{split}
\end{equation}
Applying this with $a = 4k - 4$ and $b = 2$ we see that the following equations define divisors in the respective divisor class $(4k - 4)\pi^*H + 2L_i$:
\begin{equation}
\label{eQuadForms}
\begin{array}{cccccccccl}
	Q_{1}: && t_1 t_2 X_1^2 & -&t_5 \cdots t_n Y_1^2 &+&t_3 t_4 Z_1^2 &= &0, \\
	Q_{2}: &&t_5 t_6 X_2^2 & -&t_1 \cdots t_4 \cdot t_9 \cdots t_n Y_2^2 \quad &+&t_7 t_8 Z_2^2& = &0, \\	
	\vdots &&\vdots &&\vdots &&\vdots&& \\
	Q_{k}: &&t_{n - 3} t_{n - 2} X_k^2 \quad &-&t_1 \cdots t_{n - 4} Y_k^2 &+&t_{n - 1} t_n Z_k^2 &=& 0.
	\end{array}
\end{equation}
We call $V_k$ the intersection in $\mathcal{T}_k$ of these divisors, and set $U_0 \coloneqq \{t_i \neq 0 \textup{ for all } i \in [n]\} \subset \PPP^{n-1}$.

\begin{lemma}\label{lem:complete_intersection}
	The variety $V_k$ is an integral complete intersection of these divisors. In particular, its dimension is $5k - 1$, it is Gorenstein and $K_{V_k} = -2k\pi^*H - \sum_{i = 1}^k L_i$.
\end{lemma}

\begin{proof}
	Consider the projection $\pi: V_k \to \PPP^{n-1}$. Over $U_0$ this map is smooth of relative dimension $k$ by the Jacobian criterion. Moreover, each fibre is geometrically connected so $\pi^{-1}(U_0)$ has dimension $n - 1 + k = 5k - 1$, and is geometrically irreducible and smooth. 
	
	Every irreducible component of $V_k$ has dimension at least $6k - 1 - k = 5k - 1$ as $V_k$ is defined by $k$ equations in the variety $\mathcal{T}_k$ which has dimension $4k - 1 + 2k$. We thus have to show that $V_k \setminus \pi^{-1} U_0$ has dimension lower than $n - 1 + k = 5k - 1$. The variety $V_k$ is then a complete intersection by definition. This also implies that it is integral since it is a generically reduced and irreducible complete intersection.
		
	To do this we split the complement depending on the image under $\pi$. For all $\ell = 1, \dots, n - 1$ let $U_{\ell} \subset \PPP^{n-1}$ be the locus where exactly $\ell$ of the $t_i$ are $0$. We claim that the fibres of $\pi: \pi^{-1} U_{\ell} \to U_{\ell}$ have dimension at most $k + \lfloor \frac{\ell}{2} \rfloor$. Indeed, the fibre is just a product of diagonal conics and such a conic has dimension $1$ unless all coefficients are zero, in which case it has dimension $2$. Assume that $m$ of the conics have all coefficients zero. Without loss of generality we assume that these are the first $m$ conics so $t_1 t_2 = t_3 t_4 = \cdots = t_{4m - 1} t_{4m} = 0$ which implies $2 m \leq \ell$.
	
	We then have
	\[
		\dim \pi^{-1}(U_{\ell}) \leq \dim(U_{\ell}) + k + \left \lfloor \frac{\ell}{2} \right \rfloor = 4k - 1 - \ell + k + \left \lfloor \frac{\ell}{2} \right \rfloor \leq 5k - 2,
	\] 
	as $\ell \geq 1$ which is what we had to prove.
	
	This implies that $V_k$ is Gorenstein since all complete intersections are Gorenstein \cite[\href{https://stacks.math.columbia.edu/tag/0BVA}{Tag 0BVA}]{Stacks}. The adjunction formula and Lemma \ref{lem:canonical_toric} then show that
	\[
	K_{V_k} = K_{\mathcal{T}_k} + \sum_{i = 1}^k ((4k - 4)\pi^*H + 2L_i) = -2k\pi^*H - \sum_{i = 1}^k L_i,
	\]
	as desired.
\end{proof}

Let $V_k^{\mathrm{sm}} \subset V_k$ be the smooth locus. We will consider the following open subvarieties
\begin{align*}
	U^{\mathrm{bot}} &\coloneqq \{\textup{at most one $t_i$ is $0$}\} \subset \PPP^{n-1}, \\
	U^{\mathrm{top}} &\coloneqq \left\{\begin{array}{l}\textup{at most one $t_i$ is $0$ and if such a $t_i$ exists,} \\
		\textup{then $(X_j, Z_j) \neq (0,0)$ for all $1 \leq j \leq k$}\end{array}\right\} \subset V_k.
\end{align*}

\begin{lemma} \label{lem:properties_U_k}\hfill
	\begin{enumerate}
		\item The complements of $U^{\mathrm{bot}} \subset \PPP^{n-1}$ and $U^{\mathrm{top}} \subset V_k$ have codimension $2$.
		\item We have the inclusion $U^{\mathrm{top}} \subseteq V_k^{\mathrm{sm}}$.
		\item The map $\pi: U^{\mathrm{top}} \to U^{\mathrm{bot}}$ is flat.
		\item The variety $V_k$ is normal.
	\end{enumerate}
\end{lemma}

\begin{proof}
	The dimension of the complements can be computed in the same way that the dimension of $V_k$ is computed in Lemma \ref{lem:complete_intersection}, the second claim follows from a simple computation using the Jacobian criterion.
	The third claim follows from the second claim, miracle flatness \cite[\href{https://stacks.math.columbia.edu/tag/00R4}{Tag 00R4}]{Stacks} and the fact that all fibres of this map have dimension $k$, as each fibre is a product of $k$ non-zero conics.
	The last claim follows from Serre's criterion for normality \cite[\href{https://stacks.math.columbia.edu/tag/031S}{Tag 031S}]{Stacks} by combining the first two claims and the fact that $V_k$ is Gorenstein by Lemma \ref{lem:complete_intersection} and thus Cohen--Macaulay by \cite[\href{https://stacks.math.columbia.edu/tag/0C06}{Tag 0C06}]{Stacks}.
\end{proof}
The following fact will be useful.
\begin{lemma}\label{lem:V_k_rational}
	The variety $V_k$ is rational.
\end{lemma}
\begin{proof}
	We work on the affine chart $t_{n} = 1$ and $t_{4i - 3} t_{4 i - 2} X_i = 1$ for all $1 \leq i \leq k$. For all $1 \leq i \leq k$ write $a_i = t_{4i - 3} t_{4 i - 2}$, $b_i = t_{4i - 1} t_{4 i}$ and $c_i = a_i b_i$. The equation $Q_i$ of \eqref{eQuadForms} is then after rescaling by $a_i$ equal to
	\[
	1 + c_i Z_i^2 = a_i \prod_{j \neq i} c_j Y_i^2.
	\]
	Now consider $c_i, Y_i, Z_i$ as free variables for all $1 \leq i \leq k$. We can then solve for $a_i$ and write it as a rational function of the free variables.
	
	In addition consider $t_{4k - 2}$ together with $t_{4i - 2}, t_{4i}$ for $1 \leq i \leq k - 1$ as free variables. Setting $b_i = c_i/a_i$ and $X_i = 1/a_i$ for $1 \leq i \leq k$, we then have $t_{4i - 3} = a_i/t_{4i - 2}$ for $1 \leq i \leq k$, $t_{4i - 1} = b_i/t_{4i}$ for $1 \leq i \leq k - 1$, and $t_{4k - 1} = b_k$ (since $t_{4k} = 1$). All coordinates have thus been solved in terms of rational functions of the $5k - 1$ free variables, which implies the lemma.
\end{proof}

\subsection{Picard group}
The following lemma computes the Picard group of $V_k$. Recall that for a variety $X$ over $\Q$ we write $(X)_{\overline{\Q}}$ for the base change of $X$ from $\Q$ to $\overline{\Q}$.

\begin{lemma}
	\label{lem:Picard_group_computation}
	We have $\Pic \, (V_k)_{\overline{\Q}} = \Pic V_k \cong \Pic V_k^{\mathrm{sm}} = \Pic \,(V_k^{\mathrm{sm}})_{\overline{\Q}} \cong \Z^{k + 1}$. A basis is given by $\pi^* H$ and the $L_i$ for all $i \leq k$.
\end{lemma}

\begin{proof}
	We have $\Pic V_k \subseteq \Pic \, (V_k)_{\overline{\Q}}$ by the Hochschild--Serre spectral sequence as $V_k$ is proper. We also have $\Pic V_k \subseteq \Pic V_k^{\mathrm{sm}} \cong \Pic U^{\mathrm{top}}$ both over $\Q$ and $\overline{\Q}$ as $V_k$ is normal and the complement of $U^{\mathrm{top}}$ has codimension $2$ by Lemma \ref{lem:properties_U_k}. Since $\pi^* H$ and $L_i$ live in $\Pic V_k$, it thus suffices to show that $\pi^\ast H$ and the $L_i$ form a basis for $\Pic \,(U^{\mathrm{top}})_{\overline{\Q}}$. For the sake of ease of notation we will work with the base change of the varieties over $\overline{\Q}$ for the rest of this proof.
	
	Recall that a divisor $E \subset U^{\mathrm{top}}$ is \emph{vertical} with respect to $\pi$ if $\pi(E) \subset U^{\mathrm{bot}}$ is not dense. We claim that every vertical divisor is of the form $\pi^* D$ for a divisor $D \subset U^{\mathrm{bot}}$. As $\pi$ is flat by Lemma \ref{lem:properties_U_k} this is equivalent to the following condition: $\pi^* D$ is irreducible and generically reduced for all irreducible divisors $D \subset U^{\mathrm{bot}}$.
	
	The claim is clearly true for all divisors at which the map $\pi$ is smooth. The only remaining case is when $D = \{t_i = 0\}$ for some $i$. By symmetry we may assume $i = n$. Write $T_i = t_i/t_{n - 1}$ for $i = 1, \dots, n - 2$ and $K \coloneqq \overline{\Q}(T_1, \dots, T_{n - 2})$ for the residue field of the generic point of $D$. Consider the fibre of $\pi$ at this generic point. It is a product of degenerate conics which split into a union of two lines after a quadratic extension. These quadratic extensions are
	\[
	K\left(\sqrt{-T_1 T_2/T_3 T_4} \right), \dots, K \left(\sqrt{- T_{n - 7} T_{n - 6}/T_{n - 5} T_{n - 4}}\right), K\left(\sqrt{T_1 \cdots T_{n - 4}/T_{n - 3} T_{n - 2}}\right). 
	\]
	These quadratic extensions are linearly disjoint so $\pi^* D$ is irreducible and generically reduced.
	
	Let $\eta \in \PPP^{n - 1}_{\overline{\Q}}$ be the generic point, $\bar{\eta} \in \PPP^{n - 1}_{\overline{\Q}}$ a geometric point lying over it and $U_{\eta}^{\mathrm{top}}, U_{\bar{\eta}}^{\mathrm{top}}$ the fibre of $\pi$ at $\eta, \bar{\eta}$ respectively. This is a product of $k$ smooth proper conics so $U_{\bar{\eta}}^{\mathrm{top}} \cong (\PPP_{\Q(\bar{\eta})}^1)^k$ and we have
	\begin{equation}\label{eq:Picard_product_conics}
		\bigoplus_{j = 1}^k L_j \Z \subseteq \Pic U_{\eta}^{\mathrm{top}} \subseteq \Pic U_{\bar{\eta}}^{\mathrm{top}} \cong \bigoplus_{j = 1}^k \frac{L_j}{2} \Z.
	\end{equation} 	
 	Moreover, we have the equalities 
	\[
	\overline{\Q}[U_{\eta}^{\mathrm{top}}]^{\times} = \overline{\Q}(\eta)^{\times} = \overline{\Q}(t_1/t_n, \dots, t_{n-1}/t_n)^{\times}.
	\]
	By the localisation sequence for the Picard group we get an exact sequence
	\[
	\overline{\Q}(t_1/t_n, \dots, t_{n-1}/t_n)^{\times} \to \bigoplus_{\substack{E \subset U^{\mathrm{top}} \\ \textup{vertical}}} \Z[E] \to \Pic U^{\mathrm{top}} \to \Pic U_{\eta}^{\mathrm{top}} \to 0.
	\]
	The fact that all vertical divisors are of the form $\pi^* D$ for $D \subset \PPP^{n - 1}$ implies that the first map is isomorphic to the natural map $\overline{\Q}(\PPP^{n - 1})^{\times} \to \bigoplus_{D \subset U^{\mathrm{bot}}} \Z [D]$. The cokernel of this map is $\Pic U^{\mathrm{bot}} \cong \Pic \PPP^{n - 1} \cong \Z[H]$ by definition and the fact that the complement of $U^{\mathrm{bot}}$ has codimension $2$. We thus get an induced short exact sequence $0 \to \pi^* H\Z \to \Pic U^{\mathrm{top}} \to \Pic U_{\eta}^{\mathrm{top}} \to 0$. 
	
	Consider for each $i \in [k]$ the closed subvariety $\iota_i: C_i \subset (U^{\mathrm{top}})_{\overline{\Q}}$ defined by the following equations
	\begin{align*}
		t_{4i - 3} &= 0 \textup{ and } t_j = 1 \textup{ for all } j \neq 4i - 3, \\
		(X_j, Y_j, Z_j) &= (\sqrt{-1}, 1, 1) \textup{ for all } j \neq i.
	\end{align*}
	This subvariety is isomorphic to the singular conic defined by $- Y_i^2 + Z_i^2 = 0$. We thus have $\Pic C_i \cong \Z^2$, with generators corresponding to a point on either of the two irreducible components. 
	
	The composition $C_i \to U^{\mathrm{top}} \to \PPP^{n - 1}$ factors through a point so the pullback map $\iota_i^*:\Pic U^{\mathrm{top}} \to \Pic C_i$ factors through $\Pic U^{\mathrm{top}}/\pi^* H\Z \cong \Pic U_{\eta}^{\mathrm{top}}$. On the other hand, computing the intersection with $\{X_i = 0\}$ shows that $\iota_i^* L_j$ is $(1,1) \in \Pic C_i \cong \Z^2$ if $i = j$ and $0$ otherwise. Combining this with \eqref{eq:Picard_product_conics} shows that $\Pic U_{\eta}^{\mathrm{top}} \cong \bigoplus_{j = 1}^k L_j \Z$.
\end{proof}

\begin{corollary}
\label{cor:factorial}
	The $\Q$-variety $V_k$ and the $\overline{\Q}$-variety $(V_k)_{\overline{\Q}}$ are factorial.
\end{corollary}

\begin{proof}
Immediate from Lemma \ref{lem:Picard_group_computation} and the fact that a normal variety $X$ is factorial if and only if the restriction map $\Pic X \to \Pic X^{\mathrm{sm}}$ is an isomorphism. 
\end{proof}

\begin{lemma}
\label{lem:algebraic_Brauer}
We have $\HH^1(\Q, \Pic \, (V_k^{\mathrm{sm}})_{\overline{\Q}}) = 0$.
\end{lemma}

\begin{proof}
This follows immediately from Lemma \ref{lem:Picard_group_computation} and the fact that $\HH^1(\Q, \Z) = 0$.
\end{proof}

We will also need some information on certain cones contained in the Picard group which naturally arise in the theory of birational geometry. 
\begin{definition}
	Let $X/K$ be a normal projective variety over a field $K$ of characteristic $0$. Given a divisor class $L \in \Pic X$ we let $|L|$ be the corresponding linear system. For $L \in \Pic X$ we define the \emph{stable base locus} as
	\[
	\mathbf{B}(L) \coloneqq \bigcap_{\substack{D \in |mL| \\ m \in \Z_{> 0}}} D.
	\]
	\begin{enumerate}
		\item $L$ is \emph{semi-ample} if $\mathbf{B}(L) = \varnothing$, i.e.~a multiple of $L$ is globally generated.
		\item $L$ is \emph{movable} if $\mathbf{B}(L) \subseteq X$ does not contain a divisor.
		\item $L$ is $\Q$-\emph{effective} if $\mathbf{B}(L) \neq X$, i.e.~$|m L| \neq \varnothing$ for some $m > 0$.
	\end{enumerate}
\end{definition}

For simplicity we assume that $\Pic X$ is finitely generated. These three types of divisor classes then lead to three cones in $(\Pic X)_{\R} \coloneqq \Pic X \otimes_{\Z} \R$.

\begin{definition} \hfill
	\begin{enumerate}
	\item The \emph{nef cone} $\Nef_X \subseteq (\Pic X)_{\R}$ is the cone generated by all nef $\R$-divisor classes. This is equal to the closure of the cone generated by semi-ample/ample divisor classes by \cite[Theorem 1.4.23]{Lazarsfeld}.
	\item The \emph{movable cone} $\overline{\Mov}_X \subseteq (\Pic X)_{\R}$ is the closure of the cone generated by the movable divisor classes.
	\item The \emph{pseudo-effective cone} $\overline{\Eff}_X \subseteq (\Pic X)_{\R}$ is the closure of the cone of $\Q$-effective divisor classes.
\end{enumerate}
\end{definition}

We will need the following facts about these cones.

\begin{lemma}
\label{lem:cones} 
Assume that $k \geq 2$.
	\begin{enumerate}
		\item The divisor classes $\pi^*H$ and $L_i + (2k - 3) \pi^* H$ for all $i\in[k]$ are globally generated. In particular, the nef cone $\Nef_{V_k}$ contains $\pi^*H$ and $L_i + (2k - 3) \pi^* H$.
		\item The movable cone $\overline{\Mov}_{V_k}$ contains $\Nef_{V_k}$ and the divisor class $\sum_{i = 1}^k L_i + (2k - 3)\pi^* H$. More precisely 
		\[
		\mathbf{B}\Big(\sum_{i = 1}^k L_i + (2k - 3)\pi^* H\Big) \subseteq Z = \{\textup{at least $2$ of the $Y_i$ are $0$}\} \subset V_k.
		\] 
		\item The pseudo-effective cone $\overline{\Eff}_{V_k}$ contains $\pi^*H$ and $L_i$ for all $i\in[k]$.
	\end{enumerate}
\end{lemma}

\begin{remark}
	With more work, one can show that the nef and pseudo-effective cones are generated by the above divisors.
\end{remark}

\begin{proof}
	The divisor $H$ is very ample so its pullback is globally generated. The base locus of $L_i + (2k - 3) \pi^* H$ is contained in the locus $\{X_i = Z_i = t_1^{2k - 3} Y_i = \cdots = t_n^{2k - 3} Y_i = 0\} = \varnothing$ by \eqref{eq:global_sections_Tk} with $a = 2k - 3, b = 1$.
	
	The movable cone contains the nef cone by definition. The linear system $|\sum_{i = 1}^k L_i + (2k - 3)\pi^* H|$ contains the divisors $\{X_j\prod_{i \neq j}^k Y_i = 0\}$, $\{Z_j\prod_{i \neq j}^k Y_i = 0\}$ for all $1 \leq j \leq k$ by repeated application of \eqref{eq:global_sections_Tk} with $a = 2k - 3, b = 1$ if $i = j$ and $a = 0, b = 1$ if $i \neq j$. It also contains $\{t_{\ell}^{2k - 3} \prod_{i} Y_i = 0\}$ by a similar application of \eqref{eq:global_sections_Tk}. Thus the stable base locus is contained in $Z$, which has codimension $2$.
	
	For the last part we note that the divisor classes $\pi^* H$ and $L_i = [\{Y_i = 0\}]$ are effective.
\end{proof}

We can use this lemma to bound Peyre's effective cone constant.

\begin{lemma}
\label{lem:effective_cone_constant}
	Consider the lattice $\Pic V_k$ and the cone $\Lambda = \overline{\Eff}_{V_k} \subset (\Pic V_k)_{\R}$. Define $\mathcal{X}_{\Lambda}(\cdot)$ as in \cite[Definition~2.3.14]{BT98}. Then $\mathcal{X}_{\Lambda}(-K_{V_k}) \leq \frac{1}{2k}$.
\end{lemma}

\begin{proof}
	We identify $\Pic V_k \otimes_{\Z} \R$ with $\R^{k + 1}$ using the basis from Lemma \ref{lem:Picard_group_computation}. Under this identification $\R_{\geq 0}^{k + 1} \subseteq \overline{\Eff}_{V_k}$ since the divisors $\pi^* H$ and $L_i$ for all $i$ are effective by Lemma \ref{lem:cones}. Dualising we find that $\overline{\Eff}_{V_k}^* \subseteq \R_{\geq 0}^{k + 1}$. 
	
	The measure on the dual vector space $\R^{k + 1}$ in \cite[Definition~2.3.14]{BT98} is just the Lebesgue measure and $-K_{V_k} = 2k \pi^* H + \sum_{i = 1}^k L_i$ by Lemma \ref{lem:complete_intersection} so we compute
\[
	\mathcal{X}_{\Lambda}(-K_{V_k}) = \int_{\overline{\Eff}_{V_k}^*} e^{- \langle -K_{V_k}, \mathbf{y} \rangle} \mathrm{d}\mathbf{y} \leq \int_{\R_{\geq 0}^{k + 1}} e^{- \langle -K_{V_k}, \mathbf{y} \rangle} \mathrm{d}\mathbf{y} = \int_{0}^{\infty} e^{-2k x} \mathrm{d}x \prod_{i = 1}^k \int_{0}^{\infty} e^{-y_i} \mathrm{d}y_i 
	= \frac{1}{2k},
\]
	as desired.
\end{proof}

\subsection{Adelic metrics and heights}\label{sec:adelic_metrics}
We will need an anticanonical height and Tamagawa measures on the variety $V_k$. These will both depend on a choice of adelic metric on the canonical\footnote{We define the metric on the canonical divisor instead of the anticanonical since it leads to simpler formulas.} divisor $K_{V_k}$. In this subsection, we make such a choice which will be maintained for the rest of the paper. 

To compute the $p$-adic integrals relevant to the paper we will use the fibration $\pi: V_k \to \PPP^{n-1}$; our goal will be to ultimately work towards Lemma \ref{lem:integrate_fibers_affine}. To do this we need to determine the relative canonical divisor of this fibration.

Note that this fibration is a fibre product of $k$ fibrations $\pi_i: Q_i \to \PPP^{n-1}$, where $Q_i \subset \PPP_{\PPP^{n-1}}(\mathcal{V}_i)$ is defined by the $i$-th equation of \eqref{eQuadForms}. The divisor class $L_i = [\{Y_i = 0\}] \in \Pic V_k$ is obtained by pullback of the divisor class $[\{Y_i = 0\}] \in \Pic Q_i$, which by abuse of notation we will also denote $L_i$.

\begin{lemma}
	For all $i\in[k]$ we have $K_{Q_i/\PPP^{n-1}} = 2\pi_i^* H - L_i$.
\end{lemma}

\begin{proof}
	We have $K_{Q_i/\PPP^{n-1}} = K_{Q_i} - \pi_i^*K_{\PPP^{n-1}} = K_{Q_i} + n \pi_i^* H$. The same (but easier) argument as in Lemma \ref{lem:complete_intersection} shows that $Q_i$ is a complete intersection and $K_{Q_i} = (-4k + 2) \pi_i^* H - L_i$ by the adjunction formula.
\end{proof}

Let $p_i: V_k \to Q_i$ be the natural projection. We can thus write $K_{V_k} = \pi^*K_{\PPP^{n-1}} + \sum_{i = 1}^k p_i^* K_{Q_i/\PPP^{n-1}}$, construct adelic metrics on each term separately and take the tensor product of all of these adelic metrics to define an adelic metric on $K_{V_k}$.

In the following definition we will identify sections of a line bundle $L$ with maps to $\AAA^1$ on the total space of $-L$.

\begin{definition}
	We will consider for each place $v$ of $\Q$ the standard metric 
 	$\lVert \cdot\rVert_{\PPP^{n - 1}, v}^n$ on $K_{\PPP^{n-1}} = - n H$ defined by the formula
	\[
		|\mathbf{t}|^n_v \coloneqq \max(|t_1|_v, \cdots, |t_n|_v)^n = \max(|t_1|^n_v, \cdots, |t_n|^n_v).
	\]
\end{definition}

To define a metric on $p_i^*K_{Q_i/\PPP^{n-1}} = 2\pi^* H - L_i$ we will use that $(2k - 3)\pi^* H + L_i$ is globally generated for $k \geq 2$ by Lemma \ref{lem:cones}. Moreover, the proof shows that $X_i, Z_i, t_j^{2k - 3} Y_i$ for all $1 \leq j \leq n$ are global sections that generate this line bundle. The formula 
\[
	\max(|X_i|_v,\ |\mathbf{t}|^{2k - 3}_v |Y_i|_v,\ |Z_i|_v)
\]	
thus defines a local metric at $v$ on $(3 - 2k)\pi^* H - L_i$. This leads to the following metric.

\begin{definition}
\label{def:local_metric_conic}
	Assume that $k \geq 2$. Consider for each place $v$ of $\Q$ the local metric $\lVert \cdot\rVert_{Q_i, v}$ on $K_{Q_i/\PPP^{n-1}} = 2\pi_i^* H - L_i$ given by the formula
	\[
			\max(|\mathbf{t}|^{1 - 2k}_v|X_i|_v,\ |\mathbf{t}|^{- 2}_v |Y_i|_v,\ |\mathbf{t}|^{1 - 2k}_v|Z_i|_v).
	\]
\end{definition}

Taking the products leads us to the local metric.

\begin{definition}
	Assume that $k \geq 2$. For each place $v$ of $\Q$ we will consider the local metric
	\[
	\lVert \cdot\rVert_v \coloneqq \pi^*\lVert \cdot\rVert_{\PPP^{n - 1}, v}^n \cdot \prod_{i = 1}^k p_i^*\lVert \cdot\rVert_{Q_i, v}
	\]
	 on $K_{V_k} = \pi^*K_{\PPP^{n-1}} + \sum_{i = 1}^k p_i^*K_{Q_i/\PPP^{n-1}}$. It is given by the formula
	\begin{align*}
		&|\mathbf{t}|^{n}_v \prod_{i = 1}^k\max(|\mathbf{t}|^{1 - 2k}_v|X_i|_v,\ |\mathbf{t}|^{- 2}_v |Y_i|_v,\ |\mathbf{t}|^{1 - 2k}_v|Z_i|_v) \\
		&= |\mathbf{t}|^{2k}_v \prod_{i = 1}^k\max(|\mathbf{t}|^{3 - 2k}_v|X_i|_v,\ |Y_i|_v,\ |\mathbf{t}|^{3 - 2k}_v|Z_i|_v).
	\end{align*}
\end{definition}

Let $V_k^{\mathrm{prim}}$ be the set of tuples $(\mathbf{t}, \mathbf{x}_1, \dots, \mathbf{x}_k)$ of primitive vectors $\mathbf{t} = (t_1, \dots, t_n) \in \Z_{\mathrm{prim}}^n$ and $\mathbf{x}_i = (X_i, Y_i, Z_i) \in \Z_{\mathrm{prim}}^3$ satisfying \eqref{eQuadForms}. The Cox quotient morphism from Remark~\ref{rem:Cox} restricts to a surjective $2^{k + 1}$-to-$1$ map
\[
\varpi: V_k^{\mathrm{prim}} \to V_k(\Q),
\]
corresponding to the choices of signs $(\lambda_0, \lambda_1, \dots, \lambda_k) \in \{\pm 1\}^{k + 1}$. For any $(\mathbf{t}, \mathbf{x}_1, \dots, \mathbf{x}_k) \in V_k^{\mathrm{prim}}$, primitivity implies $|\mathbf{t}|_p = 1$ and $\max(|X_i|_p, |Y_i|_p, |Z_i|_p) = 1$ for every prime $p$, so the anticanonical height $H = H_{-K_{V_k}}$ of $P = \varpi(\mathbf{t}, \mathbf{x}_1, \dots, \mathbf{x}_k) \in V_k(\Q)$ is given by
\begin{equation}
\label{eq:height_primitive}
H(\varpi(\mathbf{t}, \mathbf{x}_1, \dots, \mathbf{x}_k)) = |\mathbf{t}|^{2k} \prod_{i = 1}^k \max\left(|\mathbf{t}|^{3 - 2k}|X_i|,\ |Y_i|,\ |\mathbf{t}|^{3 - 2k}|Z_i|\right).
\end{equation}

\section{Tamagawa measures}\label{s:tamagawa}
Let us give a quick recollection on how one defines Tamagawa measures on varieties. We follow the conventions of \cite{ChambertLoirTschinkel2010}.

Let $X/\Q$ be a smooth proper variety such that $\HH^1(X, \mathcal{O}_X) = \HH^2(X, \mathcal{O}_X) = 0$. Let $U \subseteq X$ be an open subvariety and let $(\lVert \cdot\rVert_v)_v$ be an adelic metric on $\omega_U$. A construction due to Peyre for $U = X$ and generalised to all $U$ by Chambert-Loir--Tschinkel \cite[\S2.1.8]{ChambertLoirTschinkel2010} gives for all places $v$ a \emph{local Tamagawa measure} $\tau_{U, v}$ on $U(\Q_v)$. In local coordinates it is given by $|\omega|/\lVert \omega\rVert_v$ where $\omega$ is a nowhere vanishing top degree differential form and $|\omega|$ is the associated measure of \cite[\S2.1.7]{ChambertLoirTschinkel2010}.
One then defines suitable convergence factors $\lambda_v$ \cite[Theorem 2.5]{ChambertLoirTschinkel2010} and shows that the measure $\prod_v \lambda_v \tau_{U, v}$ converges absolutely. This product measure is normalised by $\lambda^* \in \R_{> 0}$ to define the \emph{global Tamagawa measure} $\tau_{U} \coloneqq (\lambda^*)^{-1} \prod_v \lambda_v \tau_{U, v}$ on $U(\Adele_{\Q})$ \cite[Definition 2.8]{ChambertLoirTschinkel2010}.

The following fact shows that there is some flexibility in choosing $U$.

\begin{lemma}
\label{lem:measure_change_codimension_2}
	Let $ V \xhookrightarrow{i}U$ be an open subvariety whose complement has codimension at least $2$ and $\tau_V$ the Tamagawa measure induced by the pullback metric. Then $\tau_U$ is the pushforward of $\tau_V$ along the map $V(\Adele_{\Q}) \to U(\Adele_{\Q})$.
\end{lemma}

\begin{proof}
	For the local measures we have $\tau_{V, v} = \tau_{U, v}\vert_{V(\Q_v)}$ by definition. The complement $U(\Q_v) \setminus V(\Q_v)$ has measure zero so $i_*\tau_{V, v} = \tau_{U, v}$. The global equality follows since the convergence factors $\lambda_v$ and normalisation factor $\lambda^*$ depend on $\Pic U_{\overline{\Q}}$ and $\overline{\Q}[U]^{\times}$ \cite[Definition 2.2]{ChambertLoirTschinkel2010}, both of which do not change when replacing $U$ by $V$.
\end{proof}

\subsection{Singular varieties}
For singular varieties $X$ the appropriate Tamagawa measures were defined in \cite[\S3]{BT98} in terms of a resolution of $X$. For varieties with terminal singularities they can be constructed without reference to a resolution using Tamagawa measures for non-proper varieties \cite{ChambertLoirTschinkel2010}. This is implicit in the discussion on p.~415 of loc.~cit. We now make it explicit.

Throughout this section, $V/\Q$ will be a normal proper variety and we let $V^{\mathrm{sm}}$ be its smooth locus. We assume that $K_V$ is $\Q$-Cartier. This means that there exists $m \in \Z_{>0}$ such that $m K_V$ is Cartier, i.e.~there exists a line bundle $\mathcal{L}$ on $V$ such that $\mathcal{L}\vert_{V^{\mathrm{sm}}} \cong \omega_{{V^{\mathrm{sm}}}}^m$.

Let $f: X \to V$ be a resolution of singularities. Recall \cite[Definition~2.2.9]{BT98} that $V$ is \emph{canonical} if there exists an effective divisor $D \subset X$ such that $f^* \mathcal{L} \cong \omega^m_X(- D)$, which must then be exceptional for $f$. Let $f_{D} \in \HH^0(X, \mathcal{O}_X(D))$ be a global section defining $D$. Batyrev--Tschinkel \cite[\S 3.3.1]{BT98} define a local Tamagawa measure on $X(\Q_v)$ for all places $v$ of $\Q$ as follows. Let $U_v \subseteq X(\Q_v)$ be an open on which there exists a nowhere vanishing top degree differential form $\omega \in \omega_X(U_v)$. The local Tamagawa measure is then defined in terms of a local metric $\lVert\cdot\rVert_{\mathcal{L}^{-1}, v}$ on $\mathcal{L}^{-1}$ as
\[
	\tau_{X, \lVert\cdot\rVert_{\mathcal{L}^{-1}, v}}\vert_{U_v} \coloneqq \lVert f_{D} \cdot \omega^{-\otimes m}\rVert_{{\mathcal{L}^{-1}, v}}^{\frac{1}{m}} \cdot |\omega|.
\]
This is independent of the choice of $\omega$ and thus glues to a measure on $X(\Q_v)$. It depends on the choice of $f_{D}$, which is only well-defined up to multiplication by an element of $\Q^{\times}$. It will become independent of this choice after multiplying over all places by the product formula.

Let $U \subseteq V$ now be any smooth open subvariety over which $f:X \to V$ is an isomorphism, e.g.~by an appropriate choice of resolution of singularities we could choose $U = V^{\mathrm{sm}}$. We want to use the metric $\lVert\cdot\rVert_{\mathcal{L}^{-1}, v}$ to define a local Tamagawa measure on $U(\Q_v)$. To do this we must fix an isomorphism $\omega_U^m \cong \mathcal{L}\vert_{U}$. Such an isomorphism is not unique, we choose the following isomorphism
\[
\omega_U^m \cong \omega_{f^{-1}(U)}^m \cong \omega^m_X\vert_{f^{-1}(U)} \xrightarrow{\cdot f_D^{-1}} \omega^m_X(- D)\vert_{{f^{-1}(U)}} \cong f^*\mathcal{L}\vert_{{{f^{-1}(U)}}} \cong \mathcal{L}\vert_{U}. 
\]
The map $\cdot f_{D}^{-1}$ is an isomorphism since $D \cap f^{-1}(U) = \varnothing$ by the assumption that $f^{-1}(U) \cong U$.

If we use this definition in the construction \cite[\S2.1.8]{ChambertLoirTschinkel2010}, the local Tamagawa measure on $U$ is explicitly defined for $U_v \subseteq U(\Q_v)$ and a nowhere vanishing top degree differential form $\omega \in \omega_U(U_v)$ as 
\[
	\tau_{U, \lVert\cdot\rVert_{\mathcal{L}^{-1}, v}}\vert_{U_v} \coloneqq \lVert f_{D} \cdot\omega^{-\otimes m}\rVert_{{\mathcal{L}^{-1}, v}}^{\frac{1}{m}} \cdot |\omega|.
\]
Let $i: U \cong f^{-1}(U) \subseteq X$. The subspace $X(\Q_v) \setminus U(\Q_v)$ has positive codimension and thus has measure $0$. We thus have the equality of measures
\begin{equation}
\label{eq:local_measures_resolution}
	\tau_{X, \lVert \cdot\rVert_{\mathcal{L}^{-1}, v}} = i_* \tau_{U, \lVert \cdot\rVert_{\mathcal{L}^{-1}, v}}.
\end{equation}
Recall that $V$ is terminal if it is canonical and the support of $D$ contains all exceptional divisors for the resolution $f: X \to V$.

\begin{prop}
\label{prop:Tamagawa_measure_terminal}
	Let $V$ be a normal proper terminal variety such that $K_V$ is $\Q$-Cartier. Fix $m \in \Z_{> 0}$ such that $m K_V$ is Cartier and fix an adelic metric on it. Let $U \subseteq V$ be a smooth open subvariety whose complement has codimension at least $2$. Let $X \to V$ be a resolution of singularities which is an isomorphism over $U$. Assume that $\HH^1(X, \mathcal{O}_X) = \HH^2(X, \mathcal{O}_X) = 0$.
	
	Under these assumptions the Tamagawa measure $\tau_U$ on $U(\Adele_{\Q})$ of \cite[Definition 2.8]{ChambertLoirTschinkel2010} and the Tamagawa measure $\tau_V^{\mathrm{BT}}$ on $X(\Adele_{\Q})$ of \cite[Definitions 3.3.9, 3.3.10]{BT98} are well-defined and $\tau_V^{\mathrm{BT}} = i_* \tau_U$.
\end{prop}
\begin{remark}
	Batyrev--Tschinkel \cite[Definitions 3.3.9, 3.3.10]{BT98} technically work in a different generality, but their definition works in the above setting.
\end{remark}
\begin{proof}
	The assumption that $\HH^1(X, \mathcal{O}_X) = \HH^2(X, \mathcal{O}_X) = 0$ is required for $\tau_U, \tau_V^{\mathrm{BT}}$ to be well-defined, see \cite[Theorem 2.5]{ChambertLoirTschinkel2010}.
	
	 To show this equality of measures it suffices to check the equality of integrals with respect to test functions of the form $f = \prod_v f_v$ with $f_v$ the unit function for all but finitely many $v$. Given \eqref{eq:local_measures_resolution} it remains to show that the convergence factors $\lambda_v$ and normalisation factor $\lambda^*$ are equal for $\tau_U$ and $\tau_V^{\mathrm{BT}}$. 
	 
	 The complement of $U \subseteq V$ has codimension at least $2$ so by the algebraic Hartogs' lemma we have $\overline{\Q}[U]^{\times} = \overline{\Q}[V]^{\times}$. The $\lambda_v, \lambda^*$ for $\tau_U$ in \cite[Definition 2.2, Theorem 2.5]{ChambertLoirTschinkel2010} thus come from the $L$-function of $\Pic U_{\overline{\Q}} \otimes_{\Z} \Q$.
	 
	 On the other hand, the $\lambda_v, \lambda^*$ for $\tau_V^{\mathrm{BT}}$ in \cite[Definitions 3.3.9, 3.3.10]{BT98} come from the $L$-function of $\Pic \, (X \setminus \Supp(D))_{\overline{\Q}} \otimes_{\Z} \Q$, where $\Supp(D)$ is the support of $D$. The codimension assumption ensures that all prime divisors contained in $X \setminus f^{-1}(U)$ are exceptional and thus contained in $\Supp(D)$ by terminality. The complement of $U \cong f^{-1}(U) \subseteq X \setminus \Supp(D)$ thus has codimension at least $2$. This implies that $\Pic \, (X \setminus \Supp(D))_{\overline{\Q}} \cong \Pic U_{\overline{\Q}}$ as $X$ is smooth, which shows that the factors $\lambda_v, \lambda^*$ agree.
\end{proof}

\begin{corollary}
\label{cor:Tamagawa_measure_converges}
	Under the assumptions of Proposition \ref{prop:Tamagawa_measure_terminal} we have $\tau_U(U(\Adele_{\Q})) < \infty$ and $\tau_{U, v}(U(\Q_v)) < \infty$ for all places $v$.
\end{corollary}

\begin{proof}
	Immediate from Proposition \ref{prop:Tamagawa_measure_terminal}, \eqref{eq:local_measures_resolution}, compactness of $X(\Q_v)$, and \cite[Theorem 3.3.8 and subsequent discussion]{BT98}.
\end{proof}

\subsection{Lifting to the affine cone}
Let $X \subseteq \PPP^{n - 1}_K$ be a projective variety over a local field $K$. Instead of computing Tamagawa measures in terms of $X$ it is often more convenient to compute them in terms of the affine cone $\tilde{X} \subseteq \AAA^n_K$. Let $\pi: \tilde{X} \setminus \{0\} \to X$ be the projection.

We will assume that $X = \{f_1 = \cdots = f_m = 0\}$ is a smooth complete intersection of codimension $m$. Let $\deg f_i = d_i$ be the degree. We will also assume $\sum_{i = 1}^m d_i < n$ which ensures that $X$ is Fano. Note that we allow the case $m = 0$.
	
By the adjunction formula we then have a canonical (up to a sign) isomorphism $\omega_X \cong \mathcal{O}(\sum_{i = 1}^m d_i - n)\vert_{X}$. Let $a_1, \dots, a_n \in K^{\times}$ be any set of elements. Let $\lVert \cdot\rVert_X$ be the metric obtained by pulling back a power of the metric $\lVert \cdot\rVert_{\PPP} \coloneqq \max(|a_1x_1|, \dots, |a_nx_n|)$ on $\mathcal{O}_{\PPP^{n-1}}(-1)$. Let $\tau_{X, \lVert \cdot\rVert_X}$ be the corresponding local Tamagawa measure.

The polynomials $f_1, \dots, f_m$ can be used to define a \emph{Leray form} $\omega_{L, X} \in \HH^0(\tilde{X}, \omega_{\tilde{X}})$, see \cite[(5.2.1)]{Peyre} for a coordinate definition. It also comes from repeated applications of the adjunction formula as in \cite[\S2.1.12]{ChambertLoirTschinkel2010}. We can then consider the induced measure $|\omega_{L, X}|$ on $\tilde{X}$. This has the following relation to $\tau_{X, \lVert \cdot\rVert_X}$.

\begin{lemma}
\label{lem:lift_to_Leray_measure}
	Let $0 \leq A \leq B$. For all measurable $g: X(K) \to \R_{\geq 0}$ we have
	\[
	\int_{\{x \in \tilde{X}(K): A \leq \lVert x\rVert_{\PPP} \leq B \}} g(\pi(x)) |\omega_{L, X}| = \int_{X(K)} g \mathrm{d}\tau_{X, \lVert \cdot\rVert_X} \cdot \int_{\{y \in K: A \leq |y| \leq B\}} |y|^{n - \sum_i d_i -1 } \mathrm{d}y.
	\]
\end{lemma}
Here $\mathrm{d}y$ is the Haar measure on $K$, implicit in the definition of Tamagawa measure \cite[\S2.1.1]{ChambertLoirTschinkel2010}.

\begin{proof}
If $f_1, \ldots, f_m \in \Q[X_1, \dots, X_n]$, $a_1 = \cdots = a_n = 1$ and $g = 1$ then the $p$-adic case with $A = 0, B = 1$ is \cite[Lemma 5.4.5]{Peyre}. The real case with $A = 0, B = H$ is \cite[Lemma 5.4.7]{Peyre}. The general case follows from the same argument.
\end{proof}

\subsubsection{\texorpdfstring{$p$-adic case}{p-adic case}}
Assume now that $K = \Q_p$. We will also assume for all $i$ that $a_i = 1$ and that $f_i \in \Z_p[X_1, \dots, X_n]$. We then have the following formula for the Leray measure.

\begin{lemma}
	We have
	\begin{multline*}
	|\omega_{L, X}|(\{x \in \tilde{X}(\Q_p): \lVert x\rVert_{\PPP} = 1\}) \\
	= \lim_{r \to \infty} 
	\frac{\#\{\mathbf{x} \in \AAA^n(\Z/p^r \Z) \textup{ primitive}: f_i(\mathbf{x}) = 0 \textup{ for all } i\in [m]\}}{p^{r \dim \tilde{X}}},
	\end{multline*}
	where the limit on the right-hand side exists.
\end{lemma}

\begin{proof}
	If $f_i \in \Z[X_1, \dots, X_n]$ then this is \cite[Lemma~3.2]{PeyreTschinkel}. Since the proof of \cite[Lemma~3.2]{PeyreTschinkel} is completely local, the argument applies also to this more general case.
\end{proof}

Since the Haar measure on $\Q_p$ is just the standard one for which $\Z_p$ has measure $1$, we may combine this with Lemma \ref{lem:lift_to_Leray_measure} to obtain our next result:

\begin{lemma}
\label{lem:TamaPrim}
	We have
	\[
	\tau_{X, \lVert \cdot\rVert_X}(X(\Q_p)) = \lim_{r \to \infty} 
	\frac{\#\{\mathbf{x} \in \AAA^n(\Z/p^r \Z) \textup{ primitive}: f_i(\mathbf{x}) = 0 \textup{ for all } i\in [m]\}}{\big(1 - \frac{1}{p}\big) p^{r \dim \tilde{X}}},
	\]
	where the limit on the right-hand side exists.
\end{lemma}

\subsection{Fibration}
We now want to apply the above theory to give a concrete formula for the Tamagawa volume of $V_k$. To do this we will use the fibration $\pi: V_k \to \PPP^{n-1}_{\Q}$. 
\begin{definition}
	Let $v$ be a place of $\Q$ and let $\tau_{V_k, v}$ be the local Tamagawa measure on $V_k^{\mathrm{sm}}(\Q_v)$ induced by the local metrics of Section~\ref{sec:adelic_metrics}. We let $\tau_v(V_k) \coloneqq \tau_{V_k, v}(V_k^{\mathrm{sm}}(\Q_v))$ be the local Tamagawa volume.
\end{definition}

Recall $U_0 = \{t_i \neq 0 \textup{ for all } i\in [n]\} \subset \PPP^{n-1}_{\Q}$ and let $\tilde{U}_0 \subset \AAA^n \setminus \{0\}$ be its affine cone. Let $v$ be a place of $\Q$. For each $\mathbf{t} \in U_0(\Q_v)$ we let $Q_{i, \mathbf{t}} \subset \PPP^2_{\Q_v}$ be the $i$-th conic in \eqref{eQuadForms} for these particular coordinates (note that scaling $\mathbf{t}$ changes the defining equation of $Q_{i, \mathbf{t}}$ but not its isomorphism class). This is a smooth conic by the definition of $U_0$.

\begin{definition}\label{def:sigmataup}
	Let $v$ be a place of $\Q$ and $i\in[k]$. Let $\tau_{Q_{i, \mathbf{t}}, v}$ be the local Tamagawa measure on $Q_{i, \mathbf{t}}(\Q_v)$ induced by the metric given by $\max(|\mathbf{t}|_v^{1 - 2k} |X_i|_v, |\mathbf{t}|_v^{-2} |Y_i|_v, |\mathbf{t}|_v^{1 - 2k} |Z_i|_v)$ on $\mathcal{O}_{\PPP^2}(- 1)$, and denote its total volume by $\tau_v(Q_{i, \mathbf{t}}) \coloneqq \tau_{Q_{i, \mathbf{t}}, v}(Q_{i, \mathbf{t}}(\Q_v))$.
	
	For a prime $p$ we write $\sigma_p(Q_{i, \mathbf{t}}) = \tau_p(Q_{i, \mathbf{t}})(1 - \frac{1}{p})$. The factor $(1 - \frac{1}{p})$ is the relevant convergence factor. Note that when $\mathbf{t} \in \Z_p^n$ is primitive (so $|\mathbf{t}|_p = 1$), the metric is $\max(|X_i|_p, |Y_i|_p, |Z_i|_p)$, and by Lemma~\ref{lem:TamaPrim} this $\sigma_p(Q_{i, \mathbf{t}})$ agrees with the $p$-adic density defined via primitive affine cone points in Section~\ref{sTernary}.
\end{definition}

Note that this measure is independent of the choice of lift $\tilde{\mathbf{t}} \in \tilde{U}_0(\Q_v)$ of $\mathbf{t}$. Indeed, if $\tilde{\mathbf{t}}' = \lambda \tilde{\mathbf{t}}$ for $\lambda \in \Q_v^{\times}$ then the map $Q_{i, \tilde{\mathbf{t}}} \to Q_{i, \tilde{\mathbf{t}}'}: (X_i, Y_i, Z_i ) \mapsto (\lambda^{2k - 1} X_i,\ \lambda^2 Y_i,\ \lambda^{2k - 1} Z_i)$ is an isomorphism which preserves the metric. The effect of this scaling on the Haar measure on the affine cone $\AAA^3$ is multiplication by $|\lambda|_v^{4k}$. But it sends the defining equation of $Q_{i, \tilde{\mathbf{t}}'}$ to $\lambda^{4k}$ times the defining equation of $Q_{i, \tilde{\mathbf{t}}}$ so it is compatible with the corresponding Leray measures. The notation $\tau_v(Q_{i, \mathbf{t}})$ is thus well-defined by Lemma \ref{lem:lift_to_Leray_measure}.

We can compute $\tau_v(V_k)$ in terms of these $\tau_v(Q_{i, \mathbf{t}})$. Let $\tau_{\PPP^{n-1}, v}$ be the standard Tamagawa measure on $\PPP^{n-1}(\Q_v)$.

\begin{lemma}\label{lem:integrate_fibers}
	Let $v$ be a place of $\Q$ and $f: \PPP^{n - 1}(\Q_v) \to \R$ non-negative or such that $f \circ \pi$ is $L^1$ with respect to $\tau_{V_k, v}$, then
	\[
	\int_{V_k^{\mathrm{sm}}(\Q_v)} f(\pi(x)) \mathrm{d}\tau_{V_k, v}(x) = \int_{U_0(\Q_v)} f(\mathbf{t})\prod_{i = 1}^k \tau_v(Q_{i, \mathbf{t}}) \mathrm{d}\tau_{\PPP^{n-1}, v}(\mathbf{t}).
	\]
\end{lemma}

\begin{proof}
	The complement of $\pi^{-1}(U_0) \subset V_k^{\mathrm{sm}}$ has positive codimension, so we can replace the integral over $V_k^{\mathrm{sm}}(\Q_v)$ with an integral over $\pi^{-1}(U_0)(\Q_v)$. 
	
	The fibration $\pi^{-1}(U_0) \to U_0$ is smooth so the induced map on $\Q_v$-points is a submersion. We constructed the metrics in Section~\ref{sec:adelic_metrics} exactly such that the assumptions of \cite[Theorem 1.22(c)]{Salberger} hold so we may repeatedly integrate over fibres to get
	\[
	\int_{\pi^{-1}(U_0)(\Q_v)} f(\pi(x)) \mathrm{d}\tau_{V_k, v}(x) = \int_{U_0(\Q_v)} f(\mathbf{t}) \prod_{i = 1}^k \tau_v(Q_{i, \mathbf{t}}) \mathrm{d}\tau_{\PPP^{n-1}, v}(\mathbf{t}).
	\]
	Note that the reference assumes that one is integrating over a compactly supported function, which is not true for this application. The extension to when $f$ is non-negative is a standard application of the monotone convergence theorem. If it is $L^1$ then we write it as a difference of its positive and negative parts.
\end{proof}

For concrete computations it will be convenient to lift this to the affine cone. Let $\mathrm{d}\mathbf{t}$ be the standard Haar measure on $\Q_v^n$. The next result follows from Lemma \ref{lem:integrate_fibers} and Lemma \ref{lem:lift_to_Leray_measure}.

\begin{lemma}
\label{lem:integrate_fibers_affine}
	If $p$ is a prime,
	\[
	\tau_p(V_k) = \left(1 - \frac{1}{p}\right)^{-1}\int_{\substack{\mathbf{t} \in \Z_p^n \\ \mathbf{t} \textup{ prim.}}} \prod_{i = 1}^k \tau_p(Q_{i, \mathbf{t}}) \mathrm{d}\mathbf{t}.
	\]
	At the infinite place, 
	\[
	\tau_\infty(V_k) = \frac{n}{2} \int_{\substack{\mathbf{t} \in \R^n \\ |\mathbf{t}| \leq 1}} \prod_{i = 1}^k \tau_\infty(Q_{i, \mathbf{t}}) \mathrm{d}\mathbf{t}.
	\]
\end{lemma}

\section{Ternary quadratic forms}
\label{sTernary}
In this section we state the relevant counting result on conics that we need.

Attached to a matrix $A \in \GL(3, \R)$, we define the norm on $\R^3$ given by
\[
\lVert\mathbf{x}\rVert_A = \lVert A\mathbf{x}\rVert_\infty
\]
for $\mathbf{x} \in \R^3$. If we are further given a non-singular quadratic form $Q(\mathbf{x}) \in \Z[x_1, x_2, x_3]$, then we define the counting function
\[
N(Q, B) \coloneqq \#\{\mathbf{x} \in \Z^3_{\mathrm{prim}} : Q(\mathbf{x}) = 0, \lVert\mathbf{x}\rVert_A \leq B\}.
\]
We also define $K_0 \coloneqq 1 + \sup_{\mathbf{x} \neq \mathbf{0}} \frac{\lVert\mathbf{x}\rVert_{I_3}}{\lVert\mathbf{x}\rVert_A}$. Finally, we introduce the densities
\[
\sigma_\infty = \lim_{\epsilon \rightarrow 0} \frac{1}{2 \epsilon} \int_{\substack{|Q(\mathbf{x})| \leq \epsilon \\ \lVert\mathbf{x}\rVert_A \leq 1}} 1 \mathrm{d}\mathbf{x}, \qquad \sigma_p = \lim_{n \rightarrow \infty} \frac{\# \{\mathbf{x} \bmod p^n : p \nmid \mathbf{x},\ Q(\mathbf{x}) \equiv 0 \bmod p^n\}}{p^{2n}}
\]
and the constant
\[
c_Q \coloneqq \frac{1}{2}\sigma_\infty \prod_p \sigma_p.
\]
Note that the real density depends both on the norm $\lVert\cdot\rVert_A$ and on $Q$, while the $p$-adic density depends only on $Q$. The following result of Sofos allows one to count the points on a conic uniformly in the size of the coefficients.

\begin{theorem}
\label{tSofos}
Let $Q(\mathbf{x}) \in \Z[x_1, x_2, x_3]$ be a non-singular isotropic quadratic form. Then we have
\[
N(Q, B) = c_Q B + O\left( (B K_0)^{1/2} (\log B K_0) \langle Q \rangle^5 \right),
\]
where $\langle Q \rangle$ denotes the maximum modulus of the coefficients of $Q$. The implied constant is absolute.
\end{theorem}

\begin{proof}
See \cite[Main Theorem 1.1]{Sofos}.
\end{proof}

We also require a generalisation of Theorem \ref{tSofos} involving (projective) congruence conditions. Let $Q\in\Z[x_1,x_2,x_3]$ be a non-singular isotropic ternary quadratic form and let $M \geq 1$. Write $M = \prod_{p \mid M} p^{m_p}$. Let $\boldsymbol{\upsilon} \in \PPP^2(\Z/M\Z)$. Define 
\[
\rho_M: \{(x, y, z) \in \Z^3 : \gcd(x, y, z, M) = 1\}\rightarrow \PPP^2(\Z/M\Z)
\] 
to be the projection map. This also defines a natural map on $\{(x, y, z) \in (\Z/M\Z)^3 : \gcd(x, y, z, M) = 1\}$, which we continue to denote by $\rho_M$ by abuse of notation. Define
\[
N_{\boldsymbol{\upsilon}}(Q, B) = \#\{\mathbf{x} \in \Z^3_{\mathrm{prim}}: Q(\mathbf{x})=0,\ \lVert\mathbf{x}\rVert_A \leq B, \ \rho_M(\mathbf{x})= \boldsymbol{\upsilon}\}.
\]
Define
\begin{equation}
\label{eSofosEuler}
c_{Q,\boldsymbol{\upsilon}} \coloneqq \frac12\sigma_\infty(Q) \prod_{p \nmid M} \sigma_p(Q) \prod_{p \mid M} \sigma_p(Q;\boldsymbol{\upsilon}),
\end{equation}
where, writing $\boldsymbol{\upsilon}_p \in \PPP^2(\Z/p^{m_p}\Z)$ for the reduction of $\boldsymbol{\upsilon}$ modulo $p^{m_p}$,
\[
\sigma_p(Q;\boldsymbol{\upsilon}) \coloneqq \lim_{n\to\infty }p^{-2n} \#\left\{\mathbf{x}\in(\Z/p^n\Z)^3: Q(\mathbf{x})\equiv0\bmod p^n,\ p\nmid \mathbf{x},\ \rho_{p^{m_p}}(\mathbf{x} \bmod p^{m_p}) = \boldsymbol{\upsilon}_p \right\}.
\]
For related estimates concerning integral points on ternary quadratic forms under local congruence conditions, see \cite{Huang}, which refines the circle method to handle congruence restrictions.

\begin{prop}
\label{prop:sofos+}
With the notation above, one has 
\[
N_{\boldsymbol{\upsilon}}(Q, B) = c_{Q,\boldsymbol{\upsilon}}B + O\left(M^3(BK_0)^{1/2}\log(BK_0)\langle Q\rangle^5\right),
\]
where the implied constant is independent of $Q$, $B$, $K_0$, $M$ and $\boldsymbol{\upsilon}$.
\end{prop}

\begin{proof}
The proof only notes the changes needed in the residue class decomposition of \cite[\S 2--6]{Sofos}. We first assume that $Q$ has the special form treated in Sections 2--5 of \cite{Sofos}, i.e.~$Q(0, 1, 0) = 0$. We use the notation $\mathbf{q} = (q_1, q_2, q_3)$ from \cite[(2.1)]{Sofos}. By \cite[Lemma 3.1]{Sofos}, up to an $O(1)$ contribution, primitive points on the conic are represented by
\[
\mathbf x = \eta \mathbf{q}(s, t)/\lambda, \qquad \lambda = \gcd(q_1(s, t), q_2(s, t), q_3(s, t)), \qquad \eta \in \{\pm1\},
\]
with $(s, t)$ primitive and $t>0$. Fix a primitive triple $\mathbf a= (a_1, a_2, a_3)\in\Z^3$ which reduces to $\boldsymbol{\upsilon}$ in $\PPP^2(\Z/M\Z)$. The condition $\rho_M(\mathbf x )=\boldsymbol{\upsilon}$ is equivalent to the requirement that 
\[
u q_i(s, t) \equiv \eta\lambda a_i \bmod \lambda M, \qquad i=1,2,3
\]
for some $u \in (\Z/M\Z)^\times$. We now follow the proof of \cite[Lemmas 3.2 and 4.2]{Sofos} with this condition included in the set being counted. The M{\"o}bius inversion steps remain unaffected by the imposition of a congruence condition. Recalling that $\lambda \mid \Delta_Q/\gcd(b, e)$ by \cite[Lemma 2.2(i)]{Sofos}, the analogue of \cite[(5.2)]{Sofos} is therefore
\[
N_{\boldsymbol{\upsilon}}(Q, B) =
\sum_{\substack{\eta\in\{\pm1\} \\ u \in (\Z/M\Z)^\times}}
\sum_{\lambda\mid \frac{\Delta_Q}{(b,e)}}
\#\left\{(s, t) \in \Z^2_{\mathrm{prim}}:
\begin{array}{l}
\lambda = \gcd(\mathbf{q}(s, t)), \\
u \mathbf{q}(s, t) \equiv \eta \lambda \mathbf{a}\bmod \lambda M, \\
\lVert\mathbf q(s, t)\rVert_A \leq B \lambda, t > 0
\end{array}
\right\}
+O(1).
\]
The first term equals
\begin{align*}
&
\sum_{\substack{\eta\in\{\pm1\} \\ u \in (\Z/M\Z)^\times}}
\sum_{k\lambda\mid \frac{\Delta_Q}{(b,e)}}\mu(k)
\sum_{\substack{(\sigma,\tau)\bmod k\lambda M\\ \gcd(\sigma,\tau,k\lambda M)=1\\
k\lambda \mid \gcd(\mathbf{q}(\sigma, \tau))\\
uq_i(\sigma,\tau)\equiv
 \eta\lambda a_i\bmod{\lambda M}}}
\#\left\{(s, t)\in\Z^2_{\mathrm{prim}}:
\begin{array}{l}
t>0,\\
(s, t)\equiv(\sigma,\tau)\bmod{k\lambda M},\\
\lVert\mathbf q(s, t)\rVert_A \leq B\lambda
\end{array}
\right\} \\
&=
\sum_{\substack{\eta\in\{\pm1\} \\ u \in (\Z/M\Z)^\times}}
\sum_{k\lambda\mid \frac{\Delta_Q}{(b,e)}}\mu(k)
\sum_{\substack{(\sigma,\tau)\bmod k\lambda M \\
\gcd(\sigma,\tau,k\lambda M)=1 \\
k\lambda \mid \gcd(\mathbf{q}(\sigma, \tau)) \\
uq_i(\sigma,\tau) \equiv \eta\lambda a_i \bmod \lambda M}}
\sum_{\substack{m\leq (2BK_0/k)^{1/2} \\ \gcd(m, k \lambda M) = 1}} \mu(m) M_{\bar m\sigma,\bar m\tau}
\left(\frac{B\lambda}{m^2}, k\lambda M\right),
\end{align*}
where $\bar m$ denotes the inverse of $m$ modulo $k\lambda M$, and
\[
M_{a, b}(T, n) \coloneqq \#\left\{(s, t)\in\Z^2 : t>0,\ (s, t)\equiv(a,b)\bmod n,\ \lVert\mathbf q(s, t)\rVert_A \leq T\right\}.
\]
For each fixed residue class $(\sigma_0, \tau_0) \bmod k\lambda$, there are exactly $M^2$ residue classes $(\sigma, \tau) \bmod k\lambda M$ reducing to $(\sigma_0, \tau_0) \bmod k\lambda$. The congruence condition
\[
uq_i(\sigma, \tau) \equiv \eta\lambda a_i \bmod \lambda M
\]
selects a subset of these $M^2$ lifts. Hence, compared to \cite{Sofos}, we are summing over a factor of at most $M^2$ more residue classes. We then apply \cite[Lemma 4.3]{Sofos} to each such residue class modulo $k \lambda M$ with
\[
T = \frac{B \lambda}{m^2}, \qquad n = k\lambda M.
\]
The proof of that lemma is a lattice-point count in a translate of $n\Z^2$, so the estimate is uniform in the chosen residue class. Therefore the error terms in the estimates following \cite[(5.2)]{Sofos} are enlarged by at most a factor $M^3$ after also summing over $u$. Define
\begin{align*}
V &\coloneqq \{(s, t) \in \R^2 : t > 0, \ \lVert\mathbf q(s, t)\rVert_A \leq 1\}, \\
c'_{Q, \boldsymbol{\upsilon}} &\coloneqq \sum_{\substack{\eta\in\{\pm1\} \\ u \in (\Z/M\Z)^\times}}
\sum_{k\lambda\mid \frac{\Delta_Q}{(b,e)}}\mu(k)
\sum_{\substack{(\sigma,\tau)\bmod k\lambda M \\
\gcd(\sigma,\tau,k\lambda M)=1 \\
k\lambda \mid \gcd(\mathbf{q}(\sigma, \tau)) \\
uq_i(\sigma,\tau) \equiv \eta\lambda a_i \bmod \lambda M}}
\sum_{\substack{m \geq 1 \\ \gcd(m, k \lambda M) = 1}} \frac{\mu(m) \Vol(V)}{k^2 \lambda M^2 m^2},
\end{align*}
which converges. This gives, for the forms $Q$ with $Q(0, 1, 0) = 0$, the asymptotic formula
\[
N_{\boldsymbol{\upsilon}}(Q, B) = c'_{Q, \boldsymbol{\upsilon}}B +
O_\epsilon\left(
M^3 (BK_0)^{1/2}\log(BK_0)
\min\{|\Delta_Q|^{1/4},\delta_Q^{1/2}\}
(|\Delta_Q|+\langle Q\rangle)\langle Q\rangle^\epsilon\right).
\]
We now explain why $c'_{Q, \boldsymbol{\upsilon}} = c_{Q, \boldsymbol{\upsilon}}$. In the case when we do not impose additional congruences as in \cite[\S 5]{Sofos}, $c'_{Q} = c_{Q}$ follows from equidistribution of rational points on $\PPP^1$ for the height $H(s,t)=\lVert\mathbf q(s, t)\rVert_A/\gcd(\mathbf{q}(s,t))$; see \cite[Proposition~6.2.15, Corollary 6.2.18]{Peyre} and \cite[Remark 10]{FMT}. The condition $\rho_M(\mathbf x) = \boldsymbol{\upsilon}$ cuts out a good open subset of the finite adelic points in the sense of \cite[Definition 3.1]{Peyre}. Therefore the same equidistribution input may be applied to this restricted adelic set. It follows that the leading constant $c'_{Q,\boldsymbol{\upsilon}}$ is Peyre's constant with an additional congruence condition $\rho_M(\mathbf x) = \boldsymbol{\upsilon}$ imposed on the adelic factor. This agrees with the densities in the formula \eqref{eSofosEuler} by \cite[Corollary 3.5]{PeyreTschinkel} for the $p$-adic densities and by \cite[Theorem 3]{Heath-Brown} and \cite[Remark 6.3]{HSS2026} for the infinite density. The factor $\frac{1}{2}$ is Peyre's effective cone constant.

Finally, the reduction to a general non-singular isotropic ternary quadratic form is exactly the reduction in \cite[\S 6]{Sofos}. The unimodular change of variables in \cite[Lemma 6.1]{Sofos} is a bijection on primitive integer vectors and on residue classes modulo $M$. The coefficient and norm bounds are the same as those used in the proof of \cite[Theorem 1.1]{Sofos}, giving the claimed error term.
\end{proof}

We will need a lower bound for the real density appearing in Theorem \ref{tSofos} or rather Proposition \ref{prop:sofos+}. We shall sometimes write $c_{Q, \boldsymbol{\upsilon}, \lVert\cdot\rVert_A}$ and $N_{\boldsymbol{\upsilon}}(Q, t, \lVert\cdot\rVert_A)$ etc.~in order to keep track of the implicit norm present in these functions. After we apply Proposition \ref{prop:sofos+}, we shall encounter similar sums that come with a complicated weight. This weight can be removed by partial summation, which we codify in our next result. 

\begin{lemma}
\label{lSofos}
Let $M \geq 1$. Let $\boldsymbol{\upsilon} \in \PPP^2(\Z/M\Z)$. Let $Q(\mathbf{x}) \in \Z[x_1, x_2, x_3]$ be a non-singular isotropic quadratic form. Let $\mathbf{s} = (s_1, s_2, s_3)$ consist of non-zero real numbers and let $S$ be the corresponding diagonal matrix. Then we have for all $\epsilon > 0$
\begin{multline*}
\sum_{\substack{B^{1/2} \leq \max(|as_1|, |bs_2|, |cs_3|) \leq B \\ \gcd(a, b, c) = 1 \\ \rho_M(a, b, c)= \boldsymbol{\upsilon} }} \frac{\mathbf{1}_{Q(a, b, c) = 0}}{\max(|as_1|, |bs_2|, |cs_3|)} \\
= \frac{c_{Q, \boldsymbol{\upsilon}, \lVert\cdot\rVert_S} \log B}{2} + O_\epsilon\left(\frac{M^3(1 + \max_i |s_i|^{-1})^{1/2 + \epsilon} \langle Q \rangle^5}{B^{\frac{1}{4} - \epsilon}}\right).
\end{multline*}
\end{lemma}

\begin{proof}
We define
\begin{align*}
a_n &= \# \left\{(a, b, c) : \begin{array}{l} \max(|as_1|, |bs_2|, |cs_3|) = n, \ \gcd(a, b, c) = 1, \\
Q(a, b, c) = 0, \ \rho_M(a, b, c)= \boldsymbol{\upsilon} \end{array}\right\}, \\
A(t) &= \sum_{\substack{n\in\R\\ 0< n \leq t\\a_n\neq 0}} a_n = N_{\boldsymbol{\upsilon}}(Q, t, \lVert\cdot\rVert_S).
\end{align*}
We observe that
\[
\sum_{\substack{B^{1/2} \leq \max(|as_1|, |bs_2|, |cs_3|) \leq B \\ \gcd(a, b, c) = 1 \\ \rho_M(a, b, c)= \boldsymbol{\upsilon} }} \frac{\mathbf{1}_{Q(a, b, c) = 0}}{\max(|as_1|, |bs_2|, |cs_3|)} = \int_{B^{1/2}}^B \frac{\mathrm{d}A(t)}{t}.
\]
Hence the lemma follows from an application of partial summation.
\end{proof}

\begin{lemma}
\label{lReal}
Let $s_1, s_2, s_3$ be non-zero real numbers and set
\[
S = 
\begin{pmatrix}
s_1 & 0 & 0 \\
0 & s_2 & 0 \\
0 & 0 & s_3
\end{pmatrix}.
\]
Let $Q(x, y, z) = ax^2 + by^2 + cz^2 \in \R[x, y, z]$ be a quadratic form with $abc\neq 0$ and $a,b,c$ not all of the same sign. Assume that 
\begin{equation}
\label{eTrickyAssu}
\frac{\sqrt{|a|/|b|}}{2} \leq |s_1| \leq 2\sqrt{|a|/|b|}, \quad \frac{1}{2} \leq |s_2| \leq 2, \quad \frac{\sqrt{|c|/|b|}}{2} \leq |s_3| \leq 2\sqrt{|c|/|b|}.
\end{equation}
Then
\[
\frac{\pi}{\sqrt{|ac|}}
\leq\sigma_\infty(Q,\lVert\cdot\rVert_S)
\leq\frac{4\pi}{\sqrt{|ac|}}.
\]
\end{lemma}

\begin{proof}
Make the change of variables
\[
X=\sqrt{|a|/|b|}x,\qquad Y=y,\qquad
Z=\sqrt{|c|/|b|}z.
\]
The quadratic form becomes $|b|q(X,Y,Z)$, where $q$ has coefficients $\pm1$, and the Jacobian is $|b|/\sqrt{|ac|}$. Thus 
\begin{gather*}
\sigma_\infty(Q, \lVert\cdot\rVert_S) = \frac{|b|}{\sqrt{|ac|}} \lim_{\epsilon \rightarrow 0} \frac{1}{2\epsilon} \int_{\substack{|bq(\mathbf{x})| \leq \epsilon \\ \mathbf{x} \in \mathcal{B}}} 1 \mathrm{d}\mathbf{x} = \frac{1}{\sqrt{|ac|}} \lim_{\epsilon \rightarrow 0} \frac{1}{2\epsilon} \int_{\substack{|q(\mathbf{x})| \leq \epsilon \\ \mathbf{x} \in \mathcal{B}}} 1 \mathrm{d}\mathbf{x}, \\
\mathcal{B} = \left\{
|X|\leq\frac{\sqrt{|a|/|b|}}{|s_1|}, \quad
|Y|\leq\frac1{|s_2|}, \quad
|Z|\leq\frac{\sqrt{|c|/|b|}}{|s_3|}
\right\}.
\end{gather*}
By \eqref{eTrickyAssu}, this box satisfies
\begin{equation}
\label{eq:boxinclusion}
[-1/2, 1/2]^3 \subseteq \mathcal{B} \subseteq [-2, 2]^3.
\end{equation}

Up to a permutation of coordinates and multiplication by $-1$, we have $q = X^2 + Y^2 - Z^2$, and \eqref{eq:boxinclusion} is still preserved. 
Fix $R \in \{1/2, 2\}$. We first fix $Z$ and integrate over $(X,Y)$. For $\sqrt{\epsilon} \leq |Z| \leq \sqrt{R^2 - \epsilon}$, we automatically have $\{(X, Y)\in\R^2:|X^2 + Y^2 - Z^2| \leq \epsilon\}\subseteq [-R, R]^2$. The integral over $(X,Y)$ satisfying $Z^2 - \epsilon \leq X^2 + Y^2 \leq Z^2 + \epsilon$ is $2\pi\epsilon$. Integrating over $Z$ gives
\[
\int_{\substack{|q(\mathbf{x})| \leq \epsilon \\ \mathbf{x} \in [-R, R]^3}} 1 \mathrm{d}\mathbf{x}
= 2 \int_{\sqrt{\epsilon}}^{\sqrt{R^2 - \epsilon}} 2\pi\epsilon \, \mathrm{d}Z + O(\epsilon^{3/2})
= 4\pi R \epsilon + O(\epsilon^{3/2}).
\]
Dividing by $2\epsilon$ and letting $\epsilon \to 0$ yields $2\pi R$, as required.
\end{proof}

\begin{lemma}
\label{lSmoothConic}
Let $p$ be an odd prime. Let $Q \in \Z[x_1, x_2, x_3]$ be a diagonal quadratic form and suppose that none of the coefficients are divisible by $p$. Then we have
\[
\sigma_p(Q) = 1 - \frac{1}{p^2}.
\]
\end{lemma}

\begin{proof}
A smooth conic over $\FF_p$ has $p + 1$ projective points, and hence $p^2 - 1$ points $\mathbf{x}$ with $p \nmid \mathbf{x}$. Every such point lifts to $p^{2n - 2}$ points in $(\Z/p^n\Z)^3$ by Hensel's lemma, and thus we obtain the result.
\end{proof}

\begin{lemma}
\label{lRamifiedConic}
Let $p$ be an odd prime. Let $Q(x_1, x_2, x_3) = ax_1^2 + bx_2^2 + cx_3^2 \in \Z[x_1, x_2, x_3]$. Assume that $p$ divides $a$ exactly once and that $p$ does not divide $b$ or $c$. Moreover, assume that $-b/c$ is a square modulo $p$. Then we have
\[
\sigma_p(Q) = 2 \left(1 - \frac{1}{p}\right).
\]
\end{lemma}

\begin{proof}
Let $r$ be such that $r^2 \equiv -b/c \bmod p$. The non-zero solutions modulo $p$ of the equation $Q(\mathbf{x}) = 0$ are then exactly given by $(v, u, ru)$ with $v \in \FF_p$, $u \in \FF_p^\ast$ and $(v, u, -ru)$ with $v \in \FF_p$, $u \in \FF_p^\ast$ and $(v, 0, 0)$ with $v \in \FF_p^\ast$. Looking at valuations, we see that the points of the shape $(v, 0, 0)$ do not lift to a solution modulo $p^2$. Moreover, every other point lifts to $p^{2n - 2}$ points in $(\Z/p^n\Z)^3$ by Hensel's lemma.
\end{proof}
\section{Local densities}
\label{s:densities}
We will need further control on the relevant local densities from
Section \ref{sTernary}.

\begin{lemma}
\label{lFormulaIntegralConic}
Let $p$ be an odd prime and let $a, b, c \in \Z_p$ be integral
and non-zero. Define $f(x, y, z) = ax^2 + by^2 + cz^2$, and let
$V$ be the variety cut out by $f = 0$.
Viewing $V \subset \PPP^2$ as a projective variety with the metric
$\lVert\cdot\rVert = \max(|x|_p, |y|_p, |z|_p)$ on $\mathcal{O}_{\PPP^2}(-1)$,
we have
	\[
\tau_{V,\lVert\cdot\rVert}(V(\Q_p))
\leq (2+v_p(b))p^{\lfloor \frac{v_p(a)+v_p(c)}{2}\rfloor}.
\]
If $v_p(a)+v_p(b)+v_p(c)=2$, then
$\tau_{V,\lVert\cdot\rVert}(V(\Q_p)) \leq 3$.
\end{lemma}

\begin{proof}
Write $a=p^{2r}a_0$ and $c=p^{2s}c_0$ with $v_p(a_0),v_p(c_0)\in\{0,1\}$, put $\beta=v_p(b)$, and let $V_0$ be the conic $f_0=a_0x^2+by^2+c_0z^2=0$ (also endowed with the metric $\lVert\cdot\rVert$). By Lemma~\ref{lem:lift_to_Leray_measure}, $\tau_{V,\lVert\cdot\rVert}(V(\Q_p))$ is the Leray measure $\omega_f$ of $C_f\coloneqq\{\mathbf{x}\in\Z_p^3: f(\mathbf{x})=0\}$. The map $\phi(x,y,z)=(p^rx,y,p^sz)$ satisfies $f=f_0\circ\phi$ and $\phi(C_f)\subseteq C_{f_0}$. Since $|\det\phi|_p=p^{-r-s}$, we have $\omega_f(E)=p^{r+s}\omega_{f_0}(\phi(E))$ for $E\subseteq C_f$. Hence
\[
\tau_{V,\lVert\cdot\rVert}(V(\Q_p))\leq p^{r+s}\,\tau_{V_0,\lVert\cdot\rVert}(V_0(\Q_p)).
\]
For a binary form $q(x,z)$ and $n \in \Z_p \setminus \{0\}$, let $r_q(n)$ denote the Leray measure of $\{(x,z)\in\Z_p^2: q(x,z)=n\}$. For units $u_1,u_2\in\Z_p^\times$ we claim the following:
\begin{enumerate}
\item[(R1)] If $-u_1u_2$ is a square, then $r_{u_1x^2+u_2z^2}(n)=(1-\frac1p)(v_p(n)+1)$.
\item[(R2)] If $-u_1u_2$ is not a square, then $r_{u_1x^2+u_2z^2}(n)\leq 1+\frac1p$.
\item[(R3)] $r_{u_1x^2+pu_2z^2}(n)\leq 2$.
\end{enumerate}
For (R1), write $-u_1u_2=w^2$ with $w\in\Z_p^\times$. The change of variables $X=u_1x-wz$, $Z=x+\frac{w}{u_1}z$ lies in $\GL_2(\Z_p)$, since its determinant is $2w$, and it turns the form into $XZ$. Solving for $Z$ then gives
\[
r_{u_1x^2+u_2z^2}(n)=\int_{|n|_p\leq|X|_p\leq1}|X|_p^{-1}\mathrm{d}X=\Big(1-\frac1p\Big)(v_p(n)+1).
\]
For (R2), the reduction of the form modulo $p$ is anisotropic, so $|q(x,z)|_p=\max(|x|_p,|z|_p)^2$. Solutions exist only if $v_p(n)=2j$ is even, and they are $p^j$ times the solutions of $q=p^{-2j}n$. This scaling does not change the Leray measure, and the latter solutions reduce to the $p+1$ smooth solutions modulo $p$, so $r_q(n)=1+\frac1p$. For (R3), if $p\nmid n$ then $x$ is a unit, and $u_1x^2\equiv n \bmod p$ has at most two solutions; Hensel's lemma gives $r_q(n)\leq2$. If $p\mid n$ then $p\mid x$, and the substitution $x=px'$ shows
\[
r_{u_1x^2+pu_2z^2}(n)=r_{u_2z^2+pu_1x'^2}(n/p),
\]
so we conclude by induction on $v_p(n)$.

Put $q=a_0x^2+c_0z^2$. Since $\mathrm{d}f_0=\mathrm{d}q+2by\mathrm{d}y$, the Leray form of $f_0$ equals $\mathrm{d}y\wedge\omega_q$ away from the line $x=z=0$, which meets $C_{f_0}$ only in $0$. Hence
\[
\tau_{V_0,\lVert\cdot\rVert}(V_0(\Q_p))=\int_{\Z_p} r_{a_0x^2+c_0z^2}(-by^2)\mathrm{d}y,
\]
and similarly with the roles of the variables permuted. We distinguish four cases.
\begin{enumerate}
\item[(i)] Suppose $a_0,c_0$ are units. We have $v_p(-by^2)=\beta+2v_p(y)$ and $\int_{\Z_p}v_p(y)\mathrm{d}y=\frac{1}{p-1}$, so (R1) and (R2) give
\[
\tau_{V_0,\lVert\cdot\rVert}(V_0(\Q_p))\leq\max\Big\{1+\frac1p,\ \Big(1-\frac1p\Big)(\beta+1)+\frac2p\Big\}\leq \beta+1+\frac1p.
\]
\item[(ii)] Suppose exactly one of $v_p(a_0),v_p(c_0)$ equals $1$. Then (R3) applies to $a_0x^2+c_0z^2$, so $\tau_{V_0,\lVert\cdot\rVert}(V_0(\Q_p))\leq 2$.
\item[(iii)] Suppose $v_p(a_0)=v_p(c_0)=1$ and $\beta=0$. We integrate over $x$ instead. Then (R3) applies to $by^2+c_0z^2$, so $\tau_{V_0,\lVert\cdot\rVert}(V_0(\Q_p))\leq 2$.
\item[(iv)] Suppose $v_p(a_0)=v_p(c_0)=1$ and $\beta\geq1$. Then $f_0=pf_1$ with $f_1=\frac{a_0}{p}x^2+\frac{b}{p}y^2+\frac{c_0}{p}z^2$, and $\omega_{f_0}=p\,\omega_{f_1}$. Applying (i) to $f_1$, with $\beta-1$ in place of $\beta$, gives
\[
\tau_{V_0,\lVert\cdot\rVert}(V_0(\Q_p))\leq p\Big(\beta+\frac1p\Big)\leq(1+\beta)p.
\]
\end{enumerate}
Note that $\lfloor\frac{v_p(a)+v_p(c)}{2}\rfloor=r+s+1$ if $v_p(a_0)=v_p(c_0)=1$, and $r+s$ otherwise. Combined with the inequality $\tau_{V,\lVert\cdot\rVert}(V(\Q_p))\leq p^{r+s}\tau_{V_0,\lVert\cdot\rVert}(V_0(\Q_p))$, cases (i)--(iv) prove the first bound.

For the final statement, the valuations of $a,b,c$ are, up to permutation, $(0,1,1)$ or $(0,0,2)$. In the first case we have $r=s=0$ and we are in case (ii) or (iii), so $\tau_{V,\lVert\cdot\rVert}(V(\Q_p))\leq2$. In the second case we integrate over the variable whose coefficient has valuation $2$; the integration formula above holds, with the same proof, for every diagonal form with non-zero coefficients in $\Z_p$, in particular for $f$. The other two coefficients are then units, so (R1) and (R2) give
\[
\tau_{V,\lVert\cdot\rVert}(V(\Q_p)) \leq \max\Big\{1 + \frac{1}{p}, \ 3\Big(1 - \frac{1}{p}\Big) + \frac{2}{p}\Big\} = 3 - \frac{1}{p} \leq 3,
\]
as desired.
\end{proof}

For later purposes, we shall also require a weighted version of Lemma \ref{lFormulaIntegralConic}.

\begin{lemma}
\label{lWeightedIntegralConic}
	Let $p$ be an odd prime and let $a, b, c \in \Z_p$ be non-zero. Define $f(x, y, z) = ax^2 + by^2 + cz^2$, and let $W \subset \PPP^2$ be the projective variety cut out by $f = 0$ endowed with the metric $\lVert\cdot\rVert = \max(|x|_p, |y|_p, |z|_p)$ on $\mathcal{O}_{\PPP^2}(-1)$.	
	Then for all $\delta > 0$, we have the upper bound
	\[
	\int_{W(\Q_p)} \left(\frac{|y|_p}{\max(|x|_p, |y|_p, |z|_p)}\right)^{\delta - 1} \mathrm{d}\tau_{W, \lVert\cdot\rVert}(x, y, z) \leq C (2+v_p(b))
p^{\lfloor\frac{v_p(a)+v_p(c)}2\rfloor}
	\]
	for some constant $C > 0$ depending only on $\delta$.
\end{lemma}

\begin{proof}
By Lemma~\ref{lem:lift_to_Leray_measure} with $A = B = 1$, up to the factor $\int_{|y|_p = 1} \mathrm{d}y = 1 - \frac{1}{p}$, we may lift the integral to primitive integral representatives $(x,y,z) \in \Z_p^3$ on the affine cone, where $\max(|x|_p, |y|_p, |z|_p) = 1$. Decompose according to $r=v_p(y)$,
ignoring the finite set with $y=0$, which has measure zero.
Let $W_r$ be the projective variety cut out by
\[
ax^2+bp^{2r}y^2+cz^2=0,
\]
endowed with the same metric. We split the integral over $W(\Q_p)$ depending on $v_p(y)$. On the part with $v_p(y)=r$,
the substitution $y=p^ry'$ contributes $p^{-r}$ to the Leray
measure, while the weight contributes $p^{(1-\delta)r}$.
Since $y'$ is a unit, the resulting representatives for $W_r$
are primitive. Enlarging to all primitive representatives
for $W_r$ and applying Lemma~\ref{lem:lift_to_Leray_measure} with $A = B = 1$ (so that the factors $1 - \frac{1}{p}$ on both sides cancel) and Lemma~\ref{lFormulaIntegralConic}, we get the bound
\begin{align*}
\int_{W(\Q_p)}
\left(
\frac{|y|_p}{\max(|x|_p,|y|_p,|z|_p)}
\right)^{\delta-1}
\mathrm{d}\tau_{W, \lVert\cdot\rVert}(x,y,z)
&\leq
\sum_{r\geq0}p^{-\delta r}\tau_{W_r, \lVert\cdot\rVert}(W_r(\Q_p))\\
&\leq
p^{\lfloor \frac{v_p(a)+v_p(c)}{2}\rfloor}
\sum_{r\geq0}(2+v_p(b)+2r)p^{-\delta r}\\
&\leq
C (2+v_p(b))
p^{\lfloor\frac{v_p(a)+v_p(c)}2\rfloor}
\end{align*}
for some constant $C>0$ depending only on $\delta$.
\end{proof}

From the weighted integral bound in Lemma \ref{lWeightedIntegralConic}, we will now derive a finiteness result for a weighted integral over the total space.

\begin{lemma}
\label{Weighted_Integral_Total}
Let $p$ be an odd prime. For all $k \geq 1$ and all $\delta > 0$,
\[
\int_{\substack{\mathbf{t} \in \Z_p^n \\
 \mathbf{t} \textup{ prim.}}}
\prod_{i=1}^k
\left(
\int_{Q_{i,\mathbf{t}}(\Q_p)}
\left(
\frac{|y_i|_p}{\max(|x_i|_p,|y_i|_p,|z_i|_p)}
\right)^{\delta-1}
\mathrm{d}\tau_{Q_{i,\mathbf{t}}, p}(x_i,y_i,z_i)
\right)
\mathrm{d}\mathbf{t}<\infty.
\]
\end{lemma}

\begin{proof}
We are free to ignore cases where coordinates of $\mathbf{t}$ are zero as these are null sets. For primitive $\mathbf{t} \in \Z_p^n$, we have $|\mathbf{t}|_p = 1$, so the metric on $Q_{i, \mathbf{t}}$ from Definition~\ref{def:sigmataup} is equal to $\max(|x_i|_p, |y_i|_p, |z_i|_p)$ on $\mathcal{O}_{\PPP^2}(-1)$.

Write $\mathbf{t}=(t_1,\dots,t_n)$. We split the integral by
fixing $v_p(t_i)=\alpha_i$, apply Lemma~\ref{lWeightedIntegralConic} to each inner integral, and then freely drop the primitivity condition $\min_i \alpha_i = 0$. If we write the equation for
each conic $Q_{i,\mathbf{t}}$ as
$a_iX^2+b_iY^2+c_iZ^2=0$, then we note that by inspecting the
explicit defining equations of $Q_{i,\mathbf{t}}$ in
\eqref{eQuadForms},
\[
v_p(a_i)+v_p(b_i)+v_p(c_i)
=v_p(t_1\cdots t_n)
=\alpha_1+\cdots+\alpha_n,
\]
and
\[
\sum_{i=1}^k\frac{v_p(a_i)+v_p(c_i)}{2}
=\frac{v_p(t_1\cdots t_n)}{2}
=\frac{\alpha_1+\cdots+\alpha_n}{2}.
\]
Since the Haar measure of $\mathbf{t}\in\Z_p^n$ with
$v_p(t_i)=\alpha_i$ for each $i$ is upper bounded by
$1/p^{\alpha_1+\cdots+\alpha_n}$, combining the above formulas
with $k$ applications of Lemma~\ref{lWeightedIntegralConic} yields
the upper bound
\[
\sum_{\alpha_1,\dots,\alpha_n\geq0}
\frac{
C^k(2+\alpha_1+\cdots+\alpha_n)^k
p^{\frac{\alpha_1+\cdots+\alpha_n}{2}}
}{
p^{\alpha_1+\cdots+\alpha_n}
}.
\]
We set $\beta\coloneqq\alpha_1+\cdots+\alpha_n$, and we note
that a given $\beta$ can occur in at most $(1+\beta)^n$
different ways. This gives the upper bound
\[
C^k\sum_{\beta\geq0}
\frac{(2 + \beta)^{k+n}}{p^{\beta/2}}<\infty,
\]
as desired.
\end{proof}

\section{Upper bounds for densities}
\label{s:tamagawabounds}
\subsection{Finite places}
Our aim in this subsection is to produce an estimate for the local Tamagawa number $\tau_p(V_k)$ at an odd finite place $p$. By Peyre's definition of $\tau_p(V_k)$, and by integration over fibres as in Lemma \ref{lem:integrate_fibers_affine}, we may write this as
\[
\tau_p(V_k) = \left(1 - \frac{1}{p}\right)^{-1} \int_{\substack{\mathbf{t} \in \Z_p^n \\ \mathbf{t} \textup{ prim.}}} \prod_{i = 1}^k \tau_p(Q_{i, \mathbf{t}}) \mathrm{d} \mathbf{t}.
\]
We split the integral in pieces by fixing the $p$-adic valuation of each $t_i$. Therefore we get
\[
\tau_p(V_k) = \left(1 - \frac{1}{p}\right)^{-1} \sum_{\substack{(k_1, \dots, k_n) \in \Z_{\geq 0}^n \\ \min(k_1, \dots, k_n) = 0}} \int_{\substack{\mathbf{t} \in \Z_p^n \\ v_p(t_i) = k_i}} \prod_{i = 1}^k \tau_p(Q_{i, \mathbf{t}}) \mathrm{d} \mathbf{t}.
\]
We evaluate the $p$-adic integral according to the value of $k_1 + \cdots + k_n$.

\begin{lemma}
\label{lem:taupestimate}
Let $p$ be an odd prime.
\begin{enumerate}
\item We have
\[
\int_{\mathbf{t} \in (\Z_p^{\times})^n } \prod_{i = 1}^k \tau_p(Q_{i, \mathbf{t}}) \mathrm{d} \mathbf{t}
=
1 - \frac{3k}{p} + O\left( \frac{k^2}{p^2}\right).
\]
\item Let $\mathbf{k} \in \Z_{\geq 0}^n$. If there exists exactly one $j$ such that $k_j=1$ and $k_{j'} = 0$ for all $j' \neq j$, then
\[
\int_{\substack{\mathbf{t} \in \Z_p^n \\ v_p(t_i) = k_i}} \prod_{i = 1}^k \tau_p(Q_{i, \mathbf{t}}) \mathrm{d} \mathbf{t}
=
\frac{1}{p} + O\left( \frac{k}{p^2} \right).
\]
\item We have
\[
\sum_{\substack{\mathbf{k} = (k_1, \dots, k_n) \in \Z_{\geq 0}^n \\ k_1 + \dots + k_n = 2}} \int_{\substack{\mathbf{t} \in \Z_p^n \\ v_p(t_i) = k_i}} \prod_{i = 1}^k \tau_p(Q_{i, \mathbf{t}}) \mathrm{d} \mathbf{t} =O\left(\frac{3^kk^2}{p^2} \right).
\]
\item For any integer $m\geq 3$, we have
\[
\sum_{\substack{\mathbf{k} = (k_1, \dots, k_n) \in \Z_{\geq 0}^n \\ k_1 + \dots + k_n \geq m}} \int_{\substack{\mathbf{t} \in \Z_p^n \\ v_p(t_i) = k_i}} \prod_{i = 1}^k \tau_p(Q_{i, \mathbf{t}}) \mathrm{d} \mathbf{t} =O\left(\frac{n^{n-1}m^{k+n-1}}{p^{\lceil m/2\rceil}} \right).
\]
\end{enumerate}
\end{lemma}

Combining the various cases provides the following estimate.

\begin{theorem}
\label{thm:local_density_total_variety}
Let $p$ be an odd prime. Then we have for each $k \geq 2$
\[
\tau_p(V_k) = 1 + \frac{k + 1}{p} + O_k\left(\frac{1}{p^{2}}\right).
\]
\end{theorem}

\begin{proof}[Proof of Lemma \ref{lem:taupestimate}]
Firstly, suppose that $\mathbf{t} \in (\Z_p^{\times})^n $. Then each conic has smooth good reduction and is geometrically integral, thus each conic contributes a local factor of $\left(1 + \frac{1}{p}\right)$. In total when integrating over all unit vectors $\mathbf{t}$, we get a contribution of size
\[
\frac{(p - 1)^n}{p^n} \left( 1 + \frac{1}{p} \right)^k = \left( 1 - \frac{1}{p^2} \right)^k \left( 1 - \frac{1}{p} \right)^{3k} = 1 - \frac{3k}{p} + O\left(\frac{k^2}{p^2}\right),
\]
where we recall that $n = 4k$. 

Secondly, suppose that there exists exactly one index $j$ such that $k_j = 1$, and $k_{j'} = 0$ for all $j' \neq j$. Inspecting the defining equations of $V_k$ in \eqref{eQuadForms}, this imposes a quadratic residue symbol on exactly $k$ linearly independent combinations of the $t_i$ and $-1$, which has density $1/2^k$. But in this case we get $\tau_p(Q_{i, \mathbf{t}}) = 2$ for each $i\in[k]$. Therefore we get, for each fixed index $j$, a contribution of size
\[
\frac{1}{p} \frac{(p - 1)^n}{p^n} \frac{2^k}{2^k} = \frac{1}{p} + O\left(\frac{k}{p^2}\right),
\]
as desired.

Thirdly, suppose that
\[
\sum_{i = 1}^n k_i = 2.
\]
This leads to two cases. The first case is $k_m = 2$ and the second case is $k_m = k_{m'} = 1$ for two distinct $m$ and $m'$. In either case, Lemma~\ref{lFormulaIntegralConic} gives
$\tau_p(Q_{i,\mathbf t})\leq 3$, and the total
contribution is at most
\[
\frac{3^k}{p^2}\binom{n+1}{2}
\ll \frac{3^kk^2}{p^2}.
\]
Fourthly, suppose that
\[
J \coloneqq \sum_{i = 1}^n k_i \geq m.
\]
Passing to the primitive vector $\mathbf{t}' = p^{-\min_j k_j}\mathbf{t}$ via scale invariance $\tau_p(Q_{i, \mathbf{t}}) = \tau_p(Q_{i, \mathbf{t}'})$ and applying Lemma \ref{lFormulaIntegralConic} to $Q_{i, \mathbf{t}'}$ to obtain the crude upper bound
\[
\tau_p(Q_{i, \mathbf{t}}) \leq (J + 2) p^{\lfloor\frac{1}{2}(v_p(t_{4i - 3}) + v_p(t_{4i - 2}) + v_p(t_{4i - 1}) + v_p(t_{4i}))\rfloor},
\]
we bound the contribution of this case by
\begin{align*}
\left( 1 - \frac{1}{p} \right)^n \sum_{J=m}^{\infty}\sum_{\substack{(k_1, \dots, k_n) \in \Z_{\geq 0}^n \\ k_1 + \dots + k_n =J}} \frac{(J + 2)^k}{p^{\lceil J/2\rceil}} &= \left( 1 - \frac{1}{p} \right)^n \sum_{J = m}^\infty \frac{ (J + 2)^k \binom{n + J - 1}{n - 1}}{p^{\lceil J/2\rceil}} \\
&\ll \frac{m^{k+n-1}}{p^{ \lceil m/2\rceil}}\sum_{r = 0}^\infty \frac{ (r+2)^k}{3^{r/2}} \binom{n + r}{n - 1}
\ll \frac{n^{n-1}m^{k+n-1}}{p^{\lceil m/2\rceil}}.
\end{align*}
This ends the proof of the lemma.
\end{proof}

We record here a simple variation when one fixes the residue class of the points on the conics. Here we shall view $Q_{i, \mathbf{t}}$ as a projective variety.

\begin{definition} 
Let $p$ be a prime, let $\mathbf{t} \in (\Z_p\setminus\{0\})^n$ be primitive, and let $\boldsymbol{\upsilon}\in\PPP^2(\Z/p\Z)$. Analogous to Definition \ref{def:sigmataup}, let $\tau_p(Q_{i, \mathbf{t}}; \boldsymbol{\upsilon})$ denote the volume of those points in $Q_{i,\mathbf{t}}(\Q_p)$ for which a primitive, integral representative reduces to $\boldsymbol{\upsilon}$ modulo $p$, with respect to the Tamagawa measure induced by the metric given by $\max(|\mathbf{t}|_p^{1 - 2k} |X_i|_p, |\mathbf{t}|_p^{-2} |Y_i|_p, |\mathbf{t}|_p^{1 - 2k} |Z_i|_p) = \max(|X_i|_p, |Y_i|_p, |Z_i|_p)$ on $\mathcal{O}_{\PPP^2}(- 1)$.

We also write
\[
\sigma_p(Q_{i, \mathbf{t}}; \boldsymbol{\upsilon}) = \left(1 - \frac{1}{p}\right) \tau_p(Q_{i, \mathbf{t}}; \boldsymbol{\upsilon}),
\]
which agrees with the $p$-adic density defined in \eqref{eSofosEuler} by \cite[Corollary 3.5]{PeyreTschinkel}.
\end{definition}

\begin{lemma}
\label{lem:densitywithcongs}
Let $p$ be an odd prime. Let $\mathbf{t} \in (\Z_p^\times)^n$, let $i \in [k]$ and let $\boldsymbol{\upsilon} \in Q_{i, \mathbf{t}}(\Z/p\Z)$. Then we have
\[
\tau_p(Q_{i,\mathbf{t}}; \boldsymbol{\upsilon}) = \frac{1}{p}.
\] 
\end{lemma}

\begin{proof}
Since $\mathbf{t} \in (\Z_p^\times)^n$, the conic $Q_{i, \mathbf{t}}$ is geometrically integral. This conic therefore has $p + 1$ projective points, while each $\boldsymbol{\upsilon}$ determines one such point. Because the various $\boldsymbol{\upsilon}$ all have the same Tamagawa measure, we deduce $\tau_p(Q_{i,\mathbf{t}}; \boldsymbol{\upsilon}) = \tau_p(Q_{i,\mathbf{t}})/(p + 1) = 1/p$.
\end{proof}

Theorem \ref{thm:local_density_total_variety} provides excellent control on the $p$-adic densities for large primes $p$. We will now deal with the small primes $p$.

\begin{lemma}
\label{lDensityLift}
Let $p$ be an odd prime. Let $a, b, c \in \Z_p$ and let $Q_{a, b, c}$ be the conic given by
\[
ax^2 + by^2 + cz^2 = 0.
\]
Let $n \geq 1$, and assume that $\max(v_p(a), v_p(b), v_p(c)) < n$. Then $\sigma_p(Q_{a, b, c})$ depends only on $a, b, c \bmod p^n$. Moreover, for all $\boldsymbol{\upsilon} \in \PPP^2(\Z/p\Z)$, the quantity $\sigma_p(Q_{a, b, c}; \boldsymbol{\upsilon})$ depends only on $a, b, c \bmod p^n$ and $\boldsymbol{\upsilon}$.

Finally, let $n \geq 1$, let $a, b, c \in \Z_2$ and assume that $\max(v_2(a), v_2(b), v_2(c)) < n$. Then, for all $\boldsymbol{\upsilon} \in \PPP^2(\Z/2\Z)$, the quantity $\sigma_2(Q_{a, b, c}; \boldsymbol{\upsilon})$ depends only on $a, b, c \bmod 2^{n + 2}$ and $\boldsymbol{\upsilon}$.
\end{lemma}

\begin{proof}
Let $a', b', c' \in \Z_p$ be such that $a \equiv a' \bmod p^n$, $b \equiv b' \bmod p^n$, and $c \equiv c' \bmod p^n$. Since $\max(v_p(a), v_p(b), v_p(c)) < n$, we deduce $v_p(a) = v_p(a')$, $v_p(b) = v_p(b')$, and $v_p(c) = v_p(c')$. Using that $p$ is odd, it thus follows from Hensel's lemma that each ratio $a/a'$, $b/b'$ and $c/c'$ is the square of some unit congruent to $1$ modulo $p$. The resulting bijection between $Q_{a, b, c}(\Q_p)$ and $Q_{a', b', c'}(\Q_p)$ is a diagonal change of coordinates, where each diagonal entry is $1$ modulo $p$. This gives the first part of the lemma.

For the second part of the lemma, the ratios $a/a'$, $b/b'$ and $c/c'$ are now $1$ modulo $8$. Hence they are squares of units in $\Z_2$. The resulting diagonal change of coordinates has odd entries, so it preserves every class in $\PPP^2(\Z/2\Z)$, and the lemma follows as before.
\end{proof}

\begin{definition}
\label{def:admissible}
For a prime $p$ and an integer $m \geq 1$, let $\adm_p(m)$ be the set of all $\mathbf{b} = (b_1, \dots, b_n) \in (\Z/p^m\Z)^n$ satisfying the following conditions:
\begin{itemize}
\item $\sum_{i = 1}^n v_p(b_i) < m$ if $p$ is odd, and $2 + \sum_{i = 1}^n v_2(b_i) < m$ if $p = 2$.
\item The ideal $(b_1, \dots, b_n)$ equals $\Z/p^m\Z$.
\item For all lifts $\mathbf{t} \in \Z_p^n$ of $\mathbf{b}$ and all $i \in [k]$, the conic $Q_{i, \mathbf{t}}(\mathbf{x}) = 0$ is soluble over $\Q_p$.
\end{itemize}
For $\mathcal{E} \geq 2$ and an integer $m \geq 1$, put $C(\mathcal{E}, m) \coloneqq \prod_{p \leq \mathcal{E}} p^m$ and define
\[
\adm(\mathcal{E}, m) \coloneqq \{\mathbf{b} \in (\Z/C(\mathcal{E}, m)\Z)^n : \mathbf{b} \bmod p^m \in \adm_p(m) \textup{ for all primes } p \leq \mathcal{E}\}.
\]
\end{definition}

Let $p$ be a prime and $\mathbf{b} \in \adm_p(m)$. Since $\sum_{i = 1}^n v_p(b_i) < m$, every lift $\mathbf{t} \in \Z_p^n$ of $\mathbf{b}$ satisfies $v_p(t_i) = v_p(b_i)$ for all $i$. Moreover, $\mathbf{t}$ is primitive since the ideal $(b_1, \dots, b_n)$ is the unit ideal, $\sum_{i = 1}^n v_p(t_i) < m$, and even $\sum_{i = 1}^n v_2(t_i) < m - 2$ if $p = 2$. Hence the coefficients of the conics $Q_{i, \mathbf{t}}$, which are products of the $t_j$, have valuation less than $m$ (resp.\ $m - 2$), and their residues modulo $p^m$ depend only on $\mathbf{b}$. For $\boldsymbol{\upsilon} = (\boldsymbol{\upsilon}_1, \dots, \boldsymbol{\upsilon}_k)$ with each $\boldsymbol{\upsilon}_i \in Q_{i, \mathbf{b}}(\Z/p\Z)$, define
\[
\mathfrak{S}_p(\mathbf{b},\boldsymbol{\upsilon}) = \prod_{i = 1}^k \sigma_p(Q_{i, \mathbf{t}}; \boldsymbol{\upsilon}_i),
\]
where $\mathbf{t}$ is any lift of $\mathbf{b}$. By Lemma~\ref{lDensityLift}, applied with $n = m$ (resp.\ with $n = m - 2$ if $p = 2$), this does not depend on the choice of lift. Moreover, we define
\[
\mathfrak{S}_p(\mathbf{b}) = \sum_{\boldsymbol{\upsilon}} \mathfrak{S}_p(\mathbf{b}, \boldsymbol{\upsilon}) = \prod_{i \in [k]} \sigma_p(Q_{i, \mathbf{t}}).
\]

\subsection{Real place}
Finally, we will give an upper bound for the real density. For $i = 1, \ldots, k$, let $Q_{i, \mathbf{t}}$ denote the $i^{\textup{th}}$ conic in the system \eqref{eQuadForms}. Recall from Definition~\ref{def:local_metric_conic} the real metric $\lVert (X_i, Y_i, Z_i) \rVert_{Q_i, \infty} = \max(\vert \mathbf{t} \vert^{-2k+1} \vert X_i \vert, \vert \mathbf{t}\vert^{-2} \vert Y_i \vert, \vert \mathbf{t} \vert^{-2k+1} \vert Z_i \vert)$ on $Q_{i, \mathbf{t}}$, so the local real norm on $V_k$ is $\vert \mathbf{t} \vert^{4k} \prod_{i = 1}^k \lVert \mathbf{x}_i \rVert_{Q_i, \infty}$.

\begin{lemma}
\label{lRealBound}
For $k \geq 2$, we have $\tau_\infty(V_k) \leq (2^{16}k)^k$.
\end{lemma}

\begin{proof}
By Lemma~\ref{lem:integrate_fibers_affine}, we have
\[
\tau_\infty(V_k) = \frac{n}{2} \int_{\vert \mathbf{t} \vert \leq 1}\prod_{i = 1}^k \tau_\infty(Q_{i, \mathbf{t}})\mathrm{d} \mathbf{t}.
\]
For $|\mathbf{t}| \leq 1$, the domain $\lVert \mathbf{x}_i \rVert_{Q_i, \infty} \leq 1$ is contained in $[-1, 1]^3$, so Lemma~\ref{lem:lift_to_Leray_measure} gives
\[
\tau_\infty(Q_{i, \mathbf{t}}) = \frac{1}{2} \int_{\{\mathbf{x}_i \in \tilde{Q}_{i, \mathbf{t}}(\R) : \lVert \mathbf{x}_i \rVert_{Q_i, \infty} \leq 1\}} |\omega_{L, Q_{i, \mathbf{t}}}|
\leq \frac{1}{2} \int_{\tilde{Q}_{i, \mathbf{t}}(\R) \cap [-1, 1]^3} |\omega_{L, Q_{i, \mathbf{t}}}|.
\]
Assume that $Q_{i, \mathbf{t}}(\R) \neq \varnothing$, otherwise $\tau_\infty(Q_{i, \mathbf{t}}) = 0$. 
Write $Q_{i,\mathbf t}(X_i,Y_i,Z_i) = a_iX_i^2+b_iY_i^2+c_iZ_i^2$ and make the change of variables
\[
u = \sqrt{|a_i|} X_i, \qquad v = \sqrt{|b_i|} Y_i, \qquad w = \sqrt{|c_i|} Z_i,
\]
which has Jacobian $\sqrt{|a_i b_i c_i|}$ and maps $[-1, 1]^3$ into $\{|u| \leq 1,\ |v| \leq \sqrt{|b_i|},\ |w| \leq 1\}$. Since $u^2 + w^2 = v^2 \leq |b_i|$ implies $|w| \leq \sqrt{|b_i|}$, up to multiplying by $-1$ and permuting $(u, v, w)$ we may assume that the transformed quadratic form is $q(u, v, w) = u^2 + v^2 - w^2$ on $\{|u| \leq 1,\ |v| \leq \sqrt{|b_i|}\}$. On each of the two sheets $w=\pm\sqrt{u^2+v^2}$, we have $|\omega_{L, q}| = \frac{\mathrm{d}u\mathrm{d}v}{2|w|}$, so the two sheets together contribute $\frac{\mathrm{d}u\mathrm{d}v}{|w|}$. We have 
\[
|w| \geq \max(|u|, |v|) \geq |u|^{1 - 2\epsilon} |v|^{2 \epsilon}
\]
for all $0 < \epsilon \leq 1/2$. Using this with $\epsilon \coloneqq \frac{1}{4(k - 1)} \in (0, 1/4]$ (for $k \geq 2$), we conclude that
\[
\tau_\infty(Q_{i, \mathbf{t}})
\leq \frac{1}{2\sqrt{|a_i b_i c_i|}} \int_{|u| \leq 1} \int_{|v| \leq \sqrt{|b_i|}} \frac{\mathrm{d}u\mathrm{d}v}{|u|^{1 - 2\epsilon}|v|^{2\epsilon}}
\leq \frac{8k}{|a_i c_i|^{1/2} |b_i|^{\epsilon}}
\]
thanks to the inequality $\frac{1}{\epsilon (1 - 2 \epsilon)} = \frac{8 (k - 1)^2}{2k - 3} \leq 8k$. Taking the product over $i \in [k]$, we have
\[
\prod_{i=1}^k \tau_\infty(Q_{i, \mathbf{t}})
\leq
\frac{(8k)^k}{\prod_{i=1}^k|a_ic_i|^{1/2}|b_i|^{\epsilon}}
=\frac{(8k)^k}{|t_1\cdots t_n|^{3/4}}.
\]
Integrating over $\mathbf{t}$ (and using $\frac{n}{2} = 2k \leq 2^k$ for $k \geq 2$), we conclude that
\[
\tau_\infty(V_k)
\leq
(16k)^k
\int_{[-1,1]^n}|t_1\cdots t_n|^{-3/4} \mathrm{d}\mathbf t=(16k)^k 8^n
=(2^{16}k)^k.
\]
This proves the lemma.
\end{proof}

\section{The character sum}
\label{sCounting}
We will control the rational points on $V_k(\Q)$ of height bounded by $B$ by using the $2^{k + 1}$-to-$1$ map $\varpi: V_k^{\mathrm{prim}} \to V_k(\Q)$ and the induced height on $V_k^{\mathrm{prim}}$ given by \eqref{eq:height_primitive}. Thus, it suffices to count primitive representatives in
\begin{equation}
\label{eNHV}
\left\{\begin{array}{c} \mathbf{t} = (t_1, \dots, t_n), \mathbf{x}_1 = (X_1, Y_1, Z_1), \dots, \mathbf{x}_k = (X_k, Y_k, Z_k) \in V_k^{\mathrm{prim}} : \\
|\mathbf{t}|^{2k} \prod_{i = 1}^k \max(|\mathbf{t}|^{-2k + 3} |X_i|, |Y_i|, |\mathbf{t}|^{-2k + 3} |Z_i|) \leq B \end{array}\right\}.
\end{equation}
We shall proceed to give a lower bound for this set.

\subsection{Introducing the character sum}
Recall that $n = 4k$, and $Q_{i, \mathbf{t}}(\mathbf{x})$ is the $i$-th quadratic form appearing in \eqref{eQuadForms}. On each quadratic form, we consider the height function
\[
H_{i, \mathbf{t}}(x, y, z) \coloneqq \max\left(|\mathbf{t}|^{-2k + 3} |x|, |y|, |\mathbf{t}|^{-2k + 3} |z|\right).
\]
For $B$ larger than an absolute constant, we introduce integers $\mathcal{E} \geq 2$ and $m \geq 3$ satisfying
\[
\mathcal{E},m \leq (\log \log B)^{1/3},
\]
and recall $C(\mathcal{E}, m) = \prod_{p \leq \mathcal{E}} p^m$ from Definition~\ref{def:admissible}.
To lower bound \eqref{eNHV}, we restrict the count to vectors $\mathbf t$ for which
$p^2\mid t_1\cdots t_n$ implies $p\leq\mathcal E$, and then split
into residue classes $\mathbf{b} \in \adm(\mathcal{E}, m)$ from
Definition~\ref{def:admissible}. 
Restricting to ``admissible'' congruence classes modulo $C(\mathcal{E},m)$ is not too lossy as most congruence classes are admissible. We also impose size restrictions on $\mathbf{t}$ and $\mathbf{x}$.

Let $\mathbf{b} \in \adm(\mathcal{E}, m)$. We are going to estimate each sum
\begin{equation}
\label{ePreSofos}
\sum_{\substack{\mathbf{t} \in \Z_{\mathrm{prim}}^n \\ p^2 \mid t_1\cdots t_n \Rightarrow p \leq \mathcal{E} \\ t_i \equiv b_i \bmod C(\mathcal{E},m) \\ |\mathbf{t}| \leq B^{1/100n^2}}} \sum_{\substack{\mathbf{x}_1, \dots, \mathbf{x}_{k - 1} \in \Z_{\mathrm{prim}}^3 \\ B^{\frac{1}{4n}}\leq H_{i, \mathbf{t}}(\mathbf{x}_i) \leq B^{\frac{1}{2n}}}} \prod_{i = 1}^{k - 1} \mathbf{1}_{Q_{i, \mathbf{t}}(\mathbf{x}_i) = 0} \sum_{\substack{\mathbf{x}_k \in \Z_{\mathrm{prim}}^3 \\ H_{k, \mathbf{t}}(\mathbf{x}_k) \leq \frac{B}{|\mathbf{t}|^{2k} \prod_{i = 1}^{k - 1} H_{i, \mathbf{t}}(\mathbf{x}_i)}}} \mathbf{1}_{Q_{k, \mathbf{t}}(\mathbf{x}_k) = 0}.
\end{equation}
Write $f_{\mathcal{E}}$ for the multiplicative function that is given on prime powers by
\[
f_{\mathcal{E}}(p^v) = 
\begin{cases}
2 \left(1 + \frac{1}{p}\right)^{-1} &\textup{if } p > \mathcal{E}, v = 1, \\
0 &\textup{if } p > \mathcal{E}, v> 1, \\
1 &\textup{if } p \leq \mathcal{E}.
\end{cases}
\]
We would now like to apply Proposition \ref{prop:sofos+} followed by an application of Lemma \ref{lReal} to control the infinite density. However, note that the condition \eqref{eTrickyAssu} does not hold automatically. In order to guarantee this assumption, we define
\[
\gamma_n\coloneqq 2^{\frac{1}{n}}, \qquad h \coloneqq \left \lfloor \frac{\log B}{100n^2\log\gamma_n} \right \rfloor = \left\lfloor\frac{\log B}{100n\log 2}\right\rfloor.
\]
We then proceed to restrict all $t_i$ simultaneously to intervals of the shape $\gamma_n^{j - 1} \leq |t_i| < \gamma_n^{j}$ for the same $j$ (where $j$ is in the range $h/2 \leq j \leq h$). This ensures that all $t_i$ are sufficiently close to each other, and hence \eqref{eTrickyAssu} is met.

Moreover, for later purposes (to control the contribution from thin sets), we also impose congruence conditions on the conics. More precisely, let $\boldsymbol{\upsilon} = (\boldsymbol{\upsilon}_1, \dots, \boldsymbol{\upsilon}_k)$ with $\boldsymbol{\upsilon}_i \in Q_{i, \mathbf{b}}(\Z/C(\mathcal{E}, 1)\Z)$ for each $i \in [k]$. We fix $\varepsilon_i \in \{\pm 1\}$ and we restrict to intervals of the shape $\gamma_n^{j - 1} \leq \varepsilon_i t_i < \gamma_n^{j}$. Define 
\[\mathcal{U}_{\mathcal{E},m}(B; V_k,\boldsymbol{\varepsilon}, \mathbf{b},\boldsymbol{\upsilon})
\coloneqq
\left\{\begin{array}{c} \mathbf{t} = (t_1, \dots, t_n), \mathbf{x}_1 = (X_1, Y_1, Z_1), \dots, \mathbf{x}_k = (X_k, Y_k, Z_k) \in V_k^{\mathrm{prim}} : \\
\gamma_n^{j - 1} \leq \varepsilon_i t_i < \gamma_n^{j}\textup{ for all }i\in[n]\textup{ for some }h/2 \leq j \leq h,\\
\mathbf{t} \equiv \mathbf{b}\bmod C(\mathcal{E},m)\textup{ and }
\rho_{C(\mathcal{E}, 1)}(\mathbf{x}_i) = \boldsymbol{\upsilon}_i \textup{ for all } i \in [k],\\
p^2 \nmid t_1 \cdots t_n\textup{ for all }p> \mathcal{E},\\
B^{\frac{1}{4n}} \leq H_{i, \mathbf{t}}(\mathbf{x}_i) \leq B^{\frac{1}{2n}}\textup{ for }i < k,\ 
H_{k, \mathbf{t}}(\mathbf{x}_k) \leq \frac{B}{|\mathbf{t}|^{2k} \prod_{i = 1}^{k - 1} H_{i, \mathbf{t}}(\mathbf{x}_i)}
\end{array}
\right\}.
\] 

We apply Proposition~\ref{prop:sofos+} to the inner sum of \eqref{ePreSofos}, and estimate the resulting densities via Lemma \ref{lReal}, Lemma \ref{lSmoothConic} and Lemma \ref{lRamifiedConic}. Then there exists an absolute constant $c_1> 0$ such that the main term of $n^{c_1 n} \cdot \# \mathcal{U}_{\mathcal{E},m}(B; V_k,\boldsymbol{\varepsilon}, \mathbf{b},\boldsymbol{\upsilon})$ is lower bounded by 
\begin{multline}
\label{eSqueezeU}
B \prod_{p \leq \mathcal{E}} \sigma_p(Q_{k, \mathbf{b}}; \boldsymbol{\upsilon}_k)\times \sum_{\frac{h}{2} \leq j \leq h} \gamma_n^{-(2j + 2kj)} \\
\sum_{\substack{\mathbf{t} \in \Z_{\mathrm{prim}}^n \\ t_i \equiv b_i \bmod C(\mathcal{E},m) \\ \gamma_n^{j - 1} \leq \varepsilon_i t_i< \gamma_n^{j}}} f_{\mathcal{E}}(|t_1 \cdots t_n|) \cdot \mathbf{1}_{Q_{k, \mathbf{t}}(\mathbf{z}) = 0 \textup{ E.L.S.}} 
\sum_{\substack{\mathbf{x}_1, \dots, \mathbf{x}_{k - 1} \in \Z_{\mathrm{prim}}^3 \\ B^{\frac{1}{4n}}\leq H_{i, \mathbf{t}}(\mathbf{x}_i) \leq B^{\frac{1}{2n}} \\ \rho_{C(\mathcal{E}, 1)}(\mathbf{x}_i) = \boldsymbol{\upsilon}_i}} \prod_{i = 1}^{k - 1} \frac{\mathbf{1}_{Q_{i, \mathbf{t}}(\mathbf{x}_i) = 0}}{H_{i, \mathbf{t}}(\mathbf{x}_i)},
\end{multline}
and moreover, the main term of $n^{-c_1 n} \cdot \# \mathcal{U}_{\mathcal{E},m}(B; V_k,\boldsymbol{\varepsilon}, \mathbf{b},\boldsymbol{\upsilon})$ is upper bounded by \eqref{eSqueezeU} as well. 

We now apply Lemma \ref{lSofos} to the $k - 1$ remaining inner sums, and we estimate the resulting densities in the same way as above. This yields that, up to an error of $O_k(B^{1-\delta_k})$ for some $\delta_k>0$ and up to factors of the size $n^{c_2n}$ for some absolute constant $c_2 > 0$, the quantity $\# \mathcal{U}_{\mathcal{E},m}(B; V_k, \boldsymbol{\varepsilon}, \mathbf{b},\boldsymbol{\upsilon})$ is both lower bounded and upper bounded by 
\begin{equation}
\label{eq:asymplb}
B (\log B)^{k - 1}\prod_{p \leq \mathcal{E}} \mathfrak{S}_p(\mathbf{b}, \boldsymbol{\upsilon})\times \sum_{\frac{h}{2} \leq j \leq h}2^{-j} \sum_{\substack{\mathbf{t} \in \Z_{\mathrm{prim}}^n \\ t_i \equiv b_i \bmod C(\mathcal{E},m) \\ \gamma_n^{j - 1} \leq\varepsilon_i t_i< \gamma_n^{j}}} f_{\mathcal{E}}(|t_1 \cdots t_n|)^k \prod_{i = 1}^k \mathbf{1}_{Q_{i, \mathbf{t}}(\mathbf{z}) = 0\textup{ E.L.S.}}.
\end{equation}

At this point we are reduced to a (weighted) count of how frequently the system of conics is soluble, for which we adapt the approach of the first three authors' paper \cite{CKR}.
This begins by detecting the indicator function, which we do in our next lemma. Note that the everywhere locally soluble condition is automatically satisfied at all primes $p \leq \mathcal{E}$ since $\mathbf{b} \in \adm(\mathcal{E}, m)$.

For each $r\in[n]$, write
\[
t_r = u_r v_r,
\]
where $u_r=\varepsilon_r\cdot \gcd(C(\mathcal{E},m),b_r)$ is divisible only by primes $p \leq \mathcal{E}$ and $v_r$ is divisible only by primes $p > \mathcal{E}$. Also set $t_{\{-\}}=u_{\{-\}}=-1$. We specify a subspace $W_r$ of the power set of $[n] \cup \{-\}$ by taking the space generated by $\{g_{i, r} : i\in[k]\}$, where 
\[g_{i, r}\coloneqq
\begin{cases}
[n]-\{4i-3, 4i-2\}&\textup{if }r\in \{4i-3,4i-2\},\\
[n]-\{4i-1,4i\}&\textup{if }r\in \{4i-1,4i\},\\
\{\{-\}, 4i-3,4i-2,4i-1, 4i\}&\textup{if }r\notin \{4i-3,4i-2, 4i-1,4i\}.
\end{cases}
\]
We denote the symmetric difference of sets by $+$.

\begin{lemma}
\label{lem:local_solubility_characters}
Let $\mathcal{E} \geq 2$ and $m \geq 3$. Let $\mathbf{t} \in \Z^n$ be such that $\mathbf{t} \equiv \mathbf{b} \bmod C(\mathcal{E},m)$ with $\mathbf{b} \in \adm(\mathcal{E},m)$, and assume that $p^2 \nmid t_1 \cdots t_n$ for all primes $p > \mathcal{E}$. Then we have
\[
\prod_{i = 1}^k \mathbf{1}_{Q_{i, \mathbf{t}}(\mathbf{z}) = 0\textup{ E.L.S.}} = 
\prod_{r=1}^n \prod_{p \mid v_r}\prod_{i = 1}^k\frac{1}{2}\left(1 + \leg{\prod_{s \in g_{i,r}} t_s}{p}\right).
\]
\end{lemma}

\begin{proof}
Since $t_i\equiv b_i\bmod C(\mathcal{E},m)$ ensures that the system is soluble at all primes $p \leq \mathcal{E}$, this follows from a direct inspection of the equations defining the $Q_{i, \mathbf{t}}$. The solubility at $\infty$ follows from Hilbert reciprocity.
\end{proof}

We will evaluate the sum
\begin{equation}
\label{eq:charsumsetup}
N_1(j, \boldsymbol{\varepsilon}, \mathbf{b})\coloneqq \sum_{\substack{\mathbf{t} \in \Z_{\mathrm{prim}}^n \\ t_i \equiv b_i \bmod C(\mathcal{E},m) \\ \gamma_n^{j - 1} \leq \varepsilon_i t_i < \gamma_n^{j}}} f_{\mathcal{E}}(|t_1 \cdots t_n|)^k 
\prod_{r=1}^n \prod_{p \mid v_r}\prod_{i = 1}^k\frac{1}{2}\left(1 + \leg{\prod_{s \in g_{i,r}} t_s}{p}\right).
\end{equation}
Noting that $f_{\mathcal{E}}(|t_1\cdots t_n|)$ is non-zero only when $v_1\cdots v_n$ is squarefree, we expand the inner product as
\begin{equation}
\label{eq:expand1}
\begin{split}
\prod_{p \mid v_r}\prod_{i = 1}^k\frac{1}{2}\left(1 + \leg{\prod_{s \in g_{i,r}} t_s}{p}\right)
&= \prod_{p \mid v_r}\frac{1}{|W_r|} \sum_{T \in W_r} \left(\frac{\prod_{s \in T} t_s}{p}\right) \\
&=\frac{1}{|W_r|^{\omega(v_r)}} \prod_{p \mid v_r}\sum_{T \in W_r} \left(\frac{\prod_{s \in T} t_s}{p}\right) \\
&=\frac{1}{|W_r|^{\omega(v_r)}}\sum_{(w_{r,T})_T} \prod_{T \in W_r} \left(\frac{\prod_{s \in T} t_s}{w_{r,T}}\right),
\end{split}
\end{equation}
where the sum over $(w_{r,T})_T$ ranges over factorisations of the form 
\[v_r= \prod_{T \in W_r} w_{r,T}.\]
We define $\mathcal{I}$ to be the set of pairs $(r, T)$, where $r\in [n]$ and $T \in W_r$.
We now make the substitution
\begin{equation}
\label{eq:tsub}
t_s = u_sv_s =u_s \prod_{T \in W_s} w_{s, T}.
\end{equation}
Inserting \eqref{eq:expand1} and \eqref{eq:tsub} into \eqref{eq:charsumsetup}
\begin{align*}
N_1(j, \boldsymbol{\varepsilon}, \mathbf{b})&=
\sum_{\substack{\mathbf{t} \in \Z_{\mathrm{prim}}^n \\ t_i \equiv b_i \bmod C(\mathcal{E},m) \\ \gamma_n^{j - 1} \leq \varepsilon_i t_i < \gamma_n^{j}}} 
f_{\mathcal{E}}(|t_1 \cdots t_n|)^k\sum_{(w_{r,T})_{(r,T)\in\mathcal{I}}} \prod_{(r,T) \in \mathcal{I}} \frac{1 }{|W_r|^{\omega(w_{r,T})}}\left(\frac{\prod_{s \in T} t_s}{w_{r,T}}\right)\\
&=\sum_{(w_{r,T})_{(r,T)\in\mathcal{I}}}^{\flat} F((w_{r,T})_{(r,T)\in\mathcal{I}}),
\end{align*}
where 
\[F((w_{r,T})_{(r,T)\in\mathcal{I}})=\frac{f_{\mathcal{E}}\left(\prod_{(r,T)\in\mathcal{I} } w_{r, T}\right)^k}{\prod_{(r,T)\in\mathcal{I}}|W_r|^{\omega(w_{r,T})}}\times \prod_{(r,T) \in \mathcal{I}} \prod_{s\in T}\left(\frac{u_s}{w_{r,T}}\right)
\times \prod_{(r,T_1) \in \mathcal{I}} \prod_{\substack{(s,T_2) \in \mathcal{I}\\ s\in T_1}}\left(\frac{w_{s, T_2}}{w_{r,T_1}}\right),\]
and $\flat$ denotes the conditions
\begin{align*}
&(w_{r,T})_{(r,T)\in\mathcal{I}}\in \Z_{\geq1}^{|\mathcal{I}|}\textup{ are pairwise coprime},&&
 u_r\prod_{T \in W_r} w_{r,T}\equiv b_r \bmod C(\mathcal{E},m),\\
 &\gamma_n^{j - 1} \leq \varepsilon_r u_r\prod_{T \in W_r} w_{r,T} < \gamma_n^{j},&
& p\mid w_{r,T} \Rightarrow p > \mathcal{E}.
\end{align*}

\subsection{Linked indices and oscillation}
We start with a definition of which several variations have appeared in previous literature.

\begin{definition}
\label{def:linked}
We say that $(r, T_1) \in \mathcal{I}$ and $(s, T_2) \in \mathcal{I}$ are linked if
\[
|\{s\} \cap T_1| + |\{r\} \cap T_2| =1.
\]
Otherwise, we say that they are unlinked. 
\end{definition}

We first treat the case when $(r, T_1)$ and $(s, T_2)$ are linked, and $w_{r, T_1}$ and $w_{s, T_2}$ are not too small. Below is a result taken from~\cite[Lemma~15]{FK4rank}, which we will use.

\begin{lemma}[Large sieve]
\label{lemma:largesieve}
Let $\epsilon > 0$. Let $a_m$ and $b_n$ be complex numbers such that $|a_m|, |b_n| \leq 1$. Then for every $M, N \geq 1$, we have
\[
\sum_{m \leq M} \sum_{n \leq N} a_m b_n \mu^2(2m) \mu^2(2n) \leg{n}{m} \ll_{\epsilon} MN (M^{-\frac{1}{2} + \epsilon} + N^{-\frac{1}{2} + \epsilon}).
\]
\end{lemma}

\begin{lemma}
\label{lemma:LS}
Let $n$ be given. For all real numbers $B\geq 3$, all $\frac{h}{2} \leq j \leq h$, all $\boldsymbol{\varepsilon} \in \{\pm 1\}^n$, all integers $\mathcal{E} \geq 2$ and $m \geq 1$, all $\mathbf{b} \in \adm(\mathcal{E}, m)$ and all linked indices $(r, T_1)$ and $(s, T_2)$, there holds
\[
\sum_{(w_{i, T})_{(i, T) \neq (r, T_1), (s, T_2)}}
\left|\sum_{w_{r, T_1},w_{s, T_2}>\exp((\log B)^{1/2})}^\flat F((w_{i, T})_{i, T})\right|
\ll_n 2^j \exp\left(-\tfrac{1}{3}(\log B)^{1/2}\right).
\]
\end{lemma}

\begin{proof}
Freezing the values of $(w_{i, T})_{(i, T) \neq (r, T_1), (s, T_2)}$, we apply Lemma \ref{lemma:largesieve} with $\epsilon = \frac{1}{12}$ to
\[
\sum_{w_{r, T_1}} \sum_{w_{s, T_2}} F((w_{i, T})_{i, T}) \ll \frac{\gamma_n^j}{\prod_{T_1' \neq T_1} w_{r, T_1'}} \cdot \frac{\gamma_n^j}{\prod_{T_2' \neq T_2} w_{s, T_2'}} \exp\left(-\tfrac{5}{12}(\log B)^{1/2}\right).
\]
Summing first over $w_{r, T'}$ and $w_{s, T'}$ contributes at most $(\log B)^{O_n(1)}$. Then summing over the variables $(w_{i, T})_T$ with $i \neq r, s$ contributes at most $\gamma_n^j (\log B)^{O_n(1)}$ for each such $i$. Since $\gamma_n^{nj} = 2^j$, this gives the desired bound.
\end{proof}

\subsection{Blocking sets}
Our current goal is to study the expression 
\[
\prod_{\substack{(r,T_1), (s,T_2) \in \mathcal{I}\\ s\in T_1}}\left(\frac{w_{s, T_2}}{w_{r,T_1}}\right)
\] 
from a combinatorial perspective. With this combinatorial description, we will be able to find the main term in our character sum.

We write $\mathcal{P}(A)$ for the power set of $A$, which we will view as a vector space over $\FF_2$ using symmetric difference of sets. Throughout this paper, we shall maintain the convention that the notation $\{r, s\}$ denotes a $2$-set, i.e.~$r$ and $s$ are implicitly assumed to be distinct. 

\begin{definition}[$\PBM$] \label{def:Ralt_PBM}
Denote by $\PBM$ the set of maps $g: [n] \rightarrow \mathcal{P}([n]\cup\{-\})$ such that
\begin{itemize}
\item $g(r) \in W_r$,
\item and the alternating properties
\[
|\{r\} \cap g(s)| = |\{s\} \cap g(r)|, \quad |\{r\} \cap g(r)| = 0
\]
\end{itemize}
hold for all $r,s\in[n]$.
\end{definition}

For $i \in [k]$, define $A^{(1)}_i=\{4i-3,4i-2\}$, $A^{(2)}_i=\{4i-1,4i\}$, and $A_i=A^{(1)}_i\cup A^{(2)}_i$. Then
\begin{align*}
W_{4i-3}=W_{4i-2}&=\langle[n]-A_i^{(1)},\{-\}+A_{i'}: i'\neq i\rangle=\langle (k-1)\{-\}+A_i^{(2)},\{-\}+A_{i'}: i'\neq i\rangle,\\
W_{4i-1}=W_{4i}&=\langle[n]-A_i^{(2)},\{-\}+A_{i'}: i'\neq i\rangle=\langle (k-1)\{-\}+A_i^{(1)},\{-\}+A_{i'}: i'\neq i\rangle.
\end{align*}
An explicit inspection of these generators has the following consequence.

\begin{remark}
\label{rem:Wr}
We have $|W_r| = 2^k$. Moreover, the map $T \mapsto T \cap [n]$ is injective when restricted to $W_r$.
\end{remark}

An $\FF_2$-basis of $\PBM$ is given by
\[
g_{\{i_1,i_2\}}(r)=\begin{cases}
\{-\} + A_{i_2}&\textup{if }r\in A_{i_1},\\
\{-\} + A_{i_1}&\textup{if }r\in A_{i_2},\\
\varnothing&\textup{otherwise},\end{cases}\qquad \{i_1,i_2\}\subseteq [k],
\]
and
\[
g'_{i}(r)=\begin{cases}
(k-1)\{-\} + A^{(2)}_i&\textup{if }r\in A^{(1)}_i,\\
(k-1)\{-\}+ A^{(1)}_i&\textup{if }r\in A^{(2)}_i,\\
\varnothing&\textup{otherwise},\end{cases}\qquad i\in [k],
\]
so the dimension of $\PBM$ is $k(k-1)/2+k=k(k+1)/2$. Every $g \in \PBM$ satisfies
\[
\sum_{r\in[n]}|\{-\}\cap g(r)|+\sum_{\{r,s\}\subseteq [n]}|\{r\}\cap g(s)|\equiv 0\bmod 2.
\]
This is easily verified on the generators of $\PBM$.

Recall the definition of unlinked indices from Definition \ref{def:linked}. 

\begin{definition}
A blocking set is a subset $\mathcal{B} \subseteq \mathcal{I}$ satisfying the following conditions:
\begin{itemize}
\item For every $r \in [n]$, there exists some $T \in W_r$ such that $(r, T) \in \mathcal{B}$.
\item Every $(r, T) \in \mathcal{B}$ (with $r \in [n]$ and $T \in W_r$) and every $(s, T') \in \mathcal{B}$ are unlinked.
\end{itemize}
\end{definition}

\begin{lemma}
\label{lBlock}
Let $\mathcal{B}$ be a blocking set. Then there exists a unique element $g \in \PBM$ with
\[
\mathcal{B} =\left\{(r, g(r)) : r \in [n]\right\}.
\]
Conversely, any set of this shape is a blocking set.
\end{lemma}

\begin{proof}
The converse is clear, so let us prove the first part. For the first part, let $\mathcal{B}$ be a blocking set. By the definition of blocking set, we may choose a subset $J_r \in W_{r}$ satisfying $(r, J_r) \in \mathcal{B}$. The elements $(r, J_r), (s, J_s)\in\mathcal{B}$ must be unlinked for all $r, s \in[n]$, so $|\{r\} \cap J_s| = |\{s\} \cap J_r|$. Since $W_{r}$ is generated by sets that do not contain $r$, we also have $|\{r\} \cap J_r| = 0$. This determines $J_r \cap [n]$. By Remark~\ref{rem:Wr}, this therefore determines $J_r$ uniquely. Hence there exists a unique $g \in \PBM$ such that $g(r) = J_r$ for every $r \in [n]$.
\end{proof}

We split the sum according to whether $w_{r, T} > \exp((\log B)^{1/2})$ for every $(r,T)\in\mathcal I$. For an arbitrary subset $\mathcal{J} \subseteq \mathcal{I}$, we define the restricted sum
\[
N_1(j, \boldsymbol{\varepsilon}, \mathbf{b},\mathcal{J}) \coloneqq \sum_{\substack{(w_{r, T})_{r, T} \\ w_{r, T} > \exp((\log B)^{1/2}) \Leftrightarrow (r, T) \in \mathcal{J} }}^{\flat} F((w_{r, T})_{r, T}).
\]
It is clear that
\[
N_1(j, \boldsymbol{\varepsilon}, \mathbf{b}) = \sum_{\mathcal{J} \subseteq \mathcal{I}} N_1(j, \boldsymbol{\varepsilon}, \mathbf{b},\mathcal{J}).
\]
The main contribution in this decomposition will come from those $\mathcal{J}$ that are blocking sets. 

\begin{lemma}
\label{lemma:notblocking}
For all $\frac{h}{2} \leq j \leq h$, 
\[
\sum_{\mathcal J\textup{ is not a blocking set}}
N_1(j,\boldsymbol\varepsilon,\mathbf b,\mathcal J)
\ll_{n} 2^j \exp\left(-\tfrac{1}{3}(\log B)^{1/2}\right).
\]
\end{lemma}

\begin{proof}
The total number of possible $\mathcal{J}$ depends only on $n$, so it suffices to bound each $N_1(j, \boldsymbol{\varepsilon}, \mathbf{b}, \mathcal{J})$ individually. The set $\mathcal{J}$ imposes the summation conditions
\[
\begin{aligned}
w_{r, T} &> \exp((\log B)^{1/2}) && \textup{for } (r, T) \in \mathcal{J}, \\
w_{r, T} &\leq \exp((\log B)^{1/2}) && \textup{for } (r, T) \notin \mathcal{J}.
\end{aligned}
\]
For $N_1(j, \boldsymbol{\varepsilon}, \mathbf{b},\mathcal{J})$ to be non-zero, there must be some $(r,T)\in\mathcal J$ for each $r\in[n]$. Otherwise this contradicts the conditions $j\geq h/2$ and
$|u_r|\prod_{T \in W_r} w_{r, T}\geq \gamma_n^{j-1}$.
Suppose that $(r, T) \in \mathcal{J}$ is linked to some $(s, T') \in \mathcal{J}$. Then we bound the contribution from all such $\mathcal{J}$ by Lemma~\ref{lemma:LS}, which gives the desired bound since there are only finitely many subsets $\mathcal J$ in terms of $n$.
For the remaining $\mathcal{J}$, all $(r, T) \in \mathcal{J}$ and $(s, T') \in \mathcal{J}$ are unlinked. Hence $\mathcal{J}$ is a blocking set by definition.
\end{proof}
 
Now, in the light of Lemma~\ref{lemma:notblocking}, assume that $\mathcal{J} = \mathcal{B}$ is a blocking set. Lemma~\ref{lBlock} allows us to associate a unique $g \in \PBM$ to each $\mathcal{B}$. 
In the remaining cases we have $w_{s, T'} \leq \exp((\log B)^{1/2})$ for all $(s, T') \notin \mathcal{B}$. For $(r,T)\in\mathcal{B}$, there is a stronger lower bound
\begin{equation}\label{eq:autolb}
w_{r,T}\geq
\frac{\gamma_n^{j-1}}{C(\mathcal E,m)\exp((|W_r|-1)(\log B)^{1/2})}
\gg_{n,\epsilon}B^{\frac{1}{200n^2}-\epsilon}
\end{equation}
from $j\geq h/2$ and $|u_r|\prod_{T' \in W_r} w_{r, T'}\geq \gamma_n^{j-1}$. Thus, for $B$ sufficiently large in terms of $n$, the condition $w_{r, T} > \exp((\log B)^{1/2})$ for $(r, T) \in \mathcal{B}$ is automatic from \eqref{eq:autolb} and can be dropped.

\subsection{A product of local factors}
For all $g\in\PBM$, all places $v$ over $\Q$ and all $\mathbf{a} \in (\Q_v^\times)^n$, analogously to \cite[(4.3)]{CKR}, define 
\begin{equation}
\label{eq:Tvdef}
T_v(g, \mathbf{a})\coloneqq\prod_{\substack{r \in [n]\\ \{-\}\in g(r)}} (-1,a_r)_v\prod_{\substack{\{r,s\} \subseteq [n]\\ s\in g(r)}}(a_r,a_s)_v
= \prod_{r \in [n]} (-1, a_r)_v^{|\{-\} \cap g(r)|} \prod_{\{r,s\} \subseteq [n]} (a_r, a_s)_v^{|\{s\} \cap g(r)|}.
\end{equation}

For a prime $p$, if $\mathbf{b} \in \adm_p(m)$, then $T_p(g, \mathbf{a})$ is the same for all choices $\mathbf{a} \in \Z_p^n$ with $\mathbf{a} \equiv \mathbf{b} \bmod p^m$. We will then write $T_p(g, \mathbf{b}) \coloneqq T_p(g, \mathbf{a})$.

\begin{lemma}
\label{lemma:tp1}
Let $p$ be an odd prime and $g\in\PBM$. Suppose
$\mathbf t\in\Z_p^n$ and $\sum_i v_p(t_i) = 1$.
Let $r\in[n]$ be the unique index with $v_p(t_r)=1$. Then
\begin{equation}
\label{eq:ufkappa}
T_p(g, \mathbf{t})= \prod_{s \in g(r)} \left(\frac{t_s}{p} \right).
\end{equation}
In particular, if $Q_{i,\mathbf t}$ is locally soluble for every
$i\in[k]$, then $T_p(g,\mathbf t)=1$.
\end{lemma}

\begin{proof}
Since $v_p(t_r)=1$ and $t_s\in\Z_p^{\times}$ for all $s\neq r$, we have
\[
(t_r,t_s)_p=\left(\frac{t_s}{p}\right)\textup{ for all }s\neq r,
\qquad
(-1,t_r)_p=\left(\frac{-1}{p}\right),
\]
and all other terms are trivially $1$.
This gives \eqref{eq:ufkappa}.

By Lemma~\ref{lem:local_solubility_characters}, the local solubility of $Q_{i, \mathbf{t}}$ is equivalent to 
\begin{equation}
\label{eq:originalg}
\prod_{s \in g_{i, r}} \left(\frac{t_s}{p}\right) =1.
\end{equation}
Since $W_r$ is generated by $g_{i,r}$, combining \eqref{eq:originalg} and \eqref{eq:ufkappa} shows that $T_p(g, \mathbf{t} ) = 1$.
\end{proof}

Given $g \in \PBM$ and its associated blocking set $\mathcal{B} = \{(r, g(r)) : r \in [n]\}$, we write 
\[w_r \coloneqq w_{r, g(r)},\qquad w'_r \coloneqq \prod_{T \in W_r \setminus \{g(r)\}} w_{r, T},\qquad W' \coloneqq \prod_{r=1}^n w'_r.\] We split
\[
N_1(j, \boldsymbol{\varepsilon}, \mathbf{b}, \mathcal{B}) = N_2(j, \boldsymbol{\varepsilon}, \mathbf{b}, \mathcal{B}) + N_3(j, \boldsymbol{\varepsilon}, \mathbf{b}, \mathcal{B})
\]
according to whether $W' = 1$ or $W' > 1$. Fix a tuple $(w_{r, T})_{(r, T) \notin \mathcal{B}}$ with $\mu^2(W') = 1$ and $p > \mathcal{E}$ for all $p \mid W'$, and consider $(w_r)_{r \in [n]}$ satisfying the remaining conditions $\flat$
\begin{align*}
&w_r\in\Z_{\geq 1} \textup{ odd}, &&\gamma_n^{j - 1}/(\varepsilon_r u_r w'_r) \leq w_{r} < \gamma_n^{j}/(\varepsilon_r u_r w'_r), \\ 
&w_r\equiv b_r/(u_r w'_r) \bmod \frac{C(\mathcal{E},m)}{|u_r|},&& p \mid w_r \Rightarrow p > \mathcal{E}.
\end{align*}
When $W'=1$, we have $t_r=u_rw_r$ and
$w_{r,T}=1$ for $(r,T)\notin \mathcal{B}$.
We have
\[
\prod_{r=1}^n
\left(\frac{\prod_{s\in g(r)}t_s}{w_r}\right)
=
\prod_{p>\mathcal E}T_p(g,\mathbf t)
=
T_\infty(g,\boldsymbol\varepsilon)
\prod_{p\leq\mathcal E}T_p(g,\mathbf b).
\]
The first equality follows from Lemma~\ref{lemma:tp1}, since each $p>\mathcal E$ divides $t_1\dots t_n$ at most once. The second follows from Hilbert reciprocity.
This gives
\begin{equation}
\label{eq:afterTp}
N_2(j, \boldsymbol{\varepsilon}, \mathbf{b}, \mathcal{B}) =
T_\infty(g, \boldsymbol{\varepsilon})
\prod_{p\leq \mathcal{E}} T_p(g, \mathbf{b})
\sum_{\substack{(w_r)_{r\in[n]}\\ W'=1}}^{\flat} \mu^2\left(\prod_{r=1}^n w_r\right) \prod_{r=1}^n \frac{f_{\mathcal{E}}(w_r)^k}{|W_r|^{\omega(w_r)}}.
\end{equation}
For $N_3(j, \boldsymbol{\varepsilon}, \mathbf{b}, \mathcal{B})$, after fixing $(w_{r, T})_{(r, T) \notin \mathcal{B}}$, quadratic reciprocity shows that up to a fixed sign, the character product in $F((w_{r, T})_{(r, T) \in \mathcal{I}})$ equals $\prod_{r=1}^n \psi_r(w_r)$, where 
\[
\psi_r(w) \coloneqq \mathbf{1}_{(w, W') = 1} \prod_{\substack{(s, T) \notin \mathcal{B} \\ r \in T + g(s)}} \leg{w}{w_{s, T}}.
\]
We obtain
\begin{equation}
\label{eq:N3bound}
|N_3(j, \boldsymbol{\varepsilon}, \mathbf{b}, \mathcal{B})|
\leq \sum_{\substack{(w_{r, T})_{(r, T) \notin \mathcal{B}}\\ W'>1,\ \mu^2(W')=1\\ 1 \leq w_{r, T} \leq \exp((\log B)^{1/2})\\ p \mid W' \Rightarrow p > \mathcal{E}}}
 \left|\sum_{(w_r)_{r\in[n]}}^{\flat} \mu^2\left(\prod_{r=1}^n w_r\right) \prod_{r=1}^n \frac{f_{\mathcal{E}}(w_r)^k}{|W_r|^{\omega(w_r)}} \psi_r(w_r)\right|.
\end{equation}
We now evaluate the inner sum over $(w_r)_{r\in[n]}$ in \eqref{eq:afterTp} and \eqref{eq:N3bound} simultaneously.

\begin{lemma}\label{lemma:subwi}
There exist constants $\eta_k>0$ such that for every sufficiently large $B$ relative to $k$ and every fixed tuple $(w_{r, T})_{(r, T) \notin \mathcal{B}}$ satisfying $1 \leq w_{r, T} \leq \exp((\log B)^{1/2})$, $p > \mathcal{E}$ for all $ p \mid W'$, $\mu^2(W') = 1$,
\[
\sum_{(w_r)_{r\in[n]}}^\flat
\mu^2\left(\prod_{r=1}^n w_r\right)
\prod_{r=1}^n\frac{f_{\mathcal E}(w_r)^k}{|W_r|^{\omega(w_r)}}\psi_r(w_r)
=\mathbf{1}_{W'=1}\frac{2^j(1-\gamma_n^{-1})^n}{C(\mathcal E,m)^n}
\left(1+O_k(\mathcal E^{-1})\right)
+O_k\left(2^j B^{-\eta_k}\right)
\]
holds uniformly for all $2 \leq \mathcal{E} \leq (\log\log B)^{1/3}$ and $3 \leq m \leq (\log\log B)^{1/3}$.
\end{lemma}

\begin{proof}
Set $C=C(\mathcal E,m)$, $C_r=C/|u_r|$ and $X_r=\gamma_n^{j-1}/(|u_r|w'_r)$. The assumptions ensure that each $C_r$ has precisely the primes $p\leq\mathcal E$ as prime divisors and that $(b_r/(u_r w'_r),C_r)=(W', C_r)=1$. Thus the congruence $w_r\equiv b_r/(u_r w'_r)\bmod C_r$ automatically implies $p > \mathcal{E}$ for every $p \mid w_r$. 

Let $a$ be the multiplicative function supported on squarefree integers coprime to $C$ with $a(p) = (1+p^{-1})^{-k} = f_{\mathcal{E}}(p)^k/|W_r|$ for $p>\mathcal{E}$ thanks to Remark~\ref{rem:Wr}. Define $r(\mathbf d)$ by
\[
R(\mathbf s)
=\sum_{\mathbf{d}\in\Z_{\geq 1}^n}
\frac{r(\mathbf d)}
 {d_1^{s_1}\cdots d_n^{s_n}}
=\prod_{p>\mathcal E}
\left(1+a(p)
\sum_{r=1}^n\frac{\psi_r(p)}{p^{s_r}}\right)
\prod_{r=1}^n\left(1-\frac{\psi_r(p)}{p^{s_r}}\right).
\] 
Since $a(p) = 1 + O_k(p^{-1})$, the Euler product converges absolutely and uniformly when $\Re(s_r)\geq 3/4$, so 
\begin{equation}\label{eq:EP34}
 \sum_{\mathbf d}\frac{|r(\mathbf d)|}{(d_1\cdots d_n)^{3/4}}\ll_k 1,
\end{equation} and $r(\mathbf d)$ is supported on $(d_r, C)=1$. When $W'=1$, we have
\begin{equation}\label{eq:R1}
R(1,\ldots,1)
=\prod_{p>\mathcal E}
(1-p^{-1})^n
\left(1+\frac np(1+p^{-1})^{-k}\right)
=1+O_k(\mathcal E^{-1}).
\end{equation}
Expanding $\prod_{p>\mathcal E}\prod_{r=1}^n(1-\psi_r(p)p^{-s_r})^{-1} R(\mathbf s)$ gives, for all $w_1,\dots,w_n$ coprime to $C$,
 \[
\mu^2(w_1\cdots w_n)\prod_{r=1}^n a(w_r)\psi_r(w_r)
=
\sum_{\substack{d_rv_r=w_r\\r\in[n]}}
r(\mathbf d)\prod_{r=1}^n\psi_r(v_r).
\]
Substituting this into the sum over $(w_r)_{r\in[n]}$ yields
\[
\sum_{(w_r)_{r\in[n]}}^\flat
\mu^2\left(\prod_{r=1}^n w_r\right)
\prod_{r=1}^n\frac{f_{\mathcal E}(w_r)^k}{|W_r|^{\omega(w_r)}}\psi_r(w_r)
=\sum_{\mathbf d\in\Z_{\geq 1}^n}
r(\mathbf d)
\prod_{r=1}^n \sum_{\substack{X_r/d_r\leq v_r<\gamma_nX_r/d_r\\ v_r\equiv b_r/(u_r w'_rd_r)\bmod C_r}}\psi_r(v_r).
\]
Since $(C_r, W') = 1$ and $X_r/C_r = \gamma_n^{j-1}/(Cw'_r)$, periodicity of $\psi_r$ modulo $W'$ along the progression modulo $C_r$ gives, for every $d_r\geq 1$,
\[
\sum_{\substack{X_r/d_r\leq v_r<\gamma_nX_r/d_r\\ v_r\equiv b_r/(u_r w'_rd_r)\bmod C_r}}\psi_r(v_r)
=
\mathbf{1}_{\psi_r\textup{ principal}}\frac{\varphi(W')}{W'}\frac{\gamma_n^j(1-\gamma_n^{-1})}{Cw'_r d_r}
+O\left(W'\min\left\{1,\frac{X_r}{d_r}\right\}\right).
\]
Observe that, thanks to Remark~\ref{rem:Wr}, $(T + g(s)) \cap [n] \neq \varnothing$ for every $(s, T) \notin \mathcal{B}$ (as then $T \in W_s$ and $T \neq g(s)$). Thus, since the $w_{s, T}$ are pairwise coprime and squarefree, $\psi_r$ is principal for all $r \in [n]$ if and only if $W' = 1$. Define $X_*=\min_rX_r$. Multiplying over $r\in[n]$ and summing over $\mathbf d\in\Z_{\geq 1}^n$ using \eqref{eq:EP34} and \eqref{eq:R1} gives
\[
\mathbf{1}_{W'=1}\frac{2^j(1-\gamma_n^{-1})^n}{C^n}\left(1+O_k(\mathcal E^{-1})\right)
+O_k\left(2^j (W')^n X_*^{-1/4}\right).
\]
Since $\mathcal E,m\leq(\log\log B)^{1/3}$, $w_{r, T} \leq \exp((\log B)^{1/2})$ and $j\geq h/2$, we have $W'=B^{o(1)}$ and $X_*\geq \gamma_n^{j-1}/(CW')\geq B^{\frac{1}{300n^2}}$ for sufficiently large $B$. Inserting $(W')^n X_*^{-1/4}\ll_k B^{-\eta_k}$ for some $\eta_k>0$ completes the proof.
\end{proof}

\subsection{The final estimate}
We now assemble the various bounds in order to estimate $\# \mathcal{U}_{\mathcal{E},m}(B; V_k,\boldsymbol{\varepsilon}, \mathbf{b},\boldsymbol{\upsilon})$. We start with a positivity argument.

\begin{lemma}
\label{lemma:charorth}
Given $\boldsymbol{\varepsilon}\in\{\pm 1\}^n$ and $\mathbf{b}\in\adm(\mathcal{E},m)$, we have
\[
\sum_{g\in\PBM}T_\infty(g,\boldsymbol\varepsilon)
\prod_{p\le\mathcal E}T_p(g,\mathbf b)
\in \{0,\lvert\PBM\rvert\}.
\]
\end{lemma}

\begin{proof}
By \eqref{eq:Tvdef}, each $T_v(g,\mathbf a)$ is a product of Hilbert symbols with exponents linear in $g$ over $\FF_2$. Therefore
\[g\mapsto T_\infty(g,\boldsymbol\varepsilon)
\prod_{p\le\mathcal E}T_p(g,\mathbf b)
\]
defines a character in $\PBM$. The claim follows from orthogonality of characters.
\end{proof}

\begin{lemma}
\label{lemma:stepU}
There exists a constant $c_3 > 0$ such that for all $k \geq 2$, for all $\mathcal{E}$ and $m$ sufficiently large relative to $k$, and for all $B \geq 100$ sufficiently large in terms of $k$, $\mathcal{E}$ and $m$, for all $\boldsymbol{\varepsilon}\in\{\pm 1\}^n$, $\mathbf{b} \in \adm(\mathcal{E},m)$ and $\boldsymbol{\upsilon}_i \in Q_{i, \mathbf{b}}(\Z/C(\mathcal{E}, 1)\Z)$ for each $i \in [k]$, we have
\[
k^{-c_3k} B (\log B)^kS_{\mathcal{E},m}(V_k,\boldsymbol{\varepsilon}, \mathbf{b},\boldsymbol{\upsilon}) \leq \# \mathcal{U}_{\mathcal{E},m}(B; V_k,\boldsymbol{\varepsilon}, \mathbf{b},\boldsymbol{\upsilon}) 
\leq k^{c_3k} B (\log B)^k S_{\mathcal{E},m}(V_k,\boldsymbol{\varepsilon}, \mathbf{b},\boldsymbol{\upsilon}),
\]
where 
\[
S_{\mathcal{E},m}(V_k,\boldsymbol{\varepsilon}, \mathbf{b},\boldsymbol{\upsilon}) \coloneqq 
\sum_{g \in \PBM} T_\infty(g, \boldsymbol{\varepsilon})\prod_{p \leq \mathcal{E}} \frac{\mathfrak{S}_p(\mathbf{b},\boldsymbol{\upsilon}) T_p(g, \mathbf{b})}{p^{mn}}\geq 0.
\]
\end{lemma}

\begin{proof}
Applying Lemma~\ref{lemma:subwi} to \eqref{eq:afterTp} and \eqref{eq:N3bound} and summing over the tuples $(w_{r, T})_{(r, T)\notin \mathcal{B}}$ (of which there are at most $\exp(|\mathcal{I}|(\log B)^{1/2}) = B^{o(1)}$) and summing over all blocking sets $\mathcal{B}$, we obtain after reducing $\eta_k > 0$
\begin{multline}
\label{eq:lovelysum}
\sum_{\mathcal{B}\textup{ is a blocking set}} N_1(j, \boldsymbol{\varepsilon}, \mathbf{b}, \mathcal{B}) \\
= \frac{2^{j }(1-\gamma_n^{-1})^n}{C(\mathcal{E},m)^n}
\left(\sum_{g \in \PBM} 
T_\infty(g, \boldsymbol{\varepsilon})
 \prod_{p\leq \mathcal{E}} T_p(g, \mathbf{b}) 
\left(1+O_k(\mathcal E^{-1})\right) + O_k(B^{-\eta_k})\right).
\end{multline}
By Lemma~\ref{lemma:charorth}, the main term in \eqref{eq:lovelysum} is non-negative.
If $ \mathcal{U}_{\mathcal{E},m}(B; V_k,\boldsymbol{\varepsilon}, \mathbf{b},\boldsymbol{\upsilon}) $ is non-empty, each $\mathfrak{S}_p(\mathbf{b},\boldsymbol{\upsilon}) $ is positive. Given $(\mathbf{t}, \mathbf{x})\in \mathcal{U}_{\mathcal{E},m}(B; V_k,\boldsymbol{\varepsilon}, \mathbf{b},\boldsymbol{\upsilon})$, for $p>\mathcal{E}$, the assumption $p^2\nmid t_1\dots t_n$ and Lemma~\ref{lemma:tp1} implies that $T_p(g,\mathbf{t})=1$. Hilbert reciprocity therefore gives
\[
T_\infty(g,\boldsymbol\varepsilon)
\prod_{p\leq\mathcal E}T_p(g,\mathbf b)=1
\]
for every $g\in\PBM$.
Thus $S_{\mathcal E,m}(V_k,\boldsymbol\varepsilon,\mathbf b,
\boldsymbol\upsilon)>0$ whenever $ \mathcal{U}_{\mathcal{E},m}(B; V_k,\boldsymbol{\varepsilon}, \mathbf{b},\boldsymbol{\upsilon}) $ is non-empty.

We combine \eqref{eq:lovelysum} with \eqref{eq:asymplb}, and then sum over $h/2\leq j\leq h$. We have
\[
\#\{j\in\Z:h/2\leq j\leq h\}\asymp\frac{\log B}{n},
\qquad
(1-\gamma_n^{-1})^n=n^{-n}\exp(O(n)).
\]
These factors, together with the factors $n^{O(n)}$ in
\eqref{eq:asymplb}, are absorbed into $k^{\pm c_3k}$ after enlarging
the absolute constant $c_3$. The error terms from Lemma~\ref{lemma:notblocking} and \eqref{eq:asymplb} can also be absorbed after further enlarging $c_3$. Taking $\mathcal E$ sufficiently large
in terms of $k$, and then $B$ sufficiently large, proves the claim.
\end{proof}

\begin{lemma}
\label{lemma:sigmahat}
Let $k \geq 2$, let $p>k^{12k}$ be a prime, and let $m\geq 2$. For all $g\in\PBM$, we have
\[
\frac{1}{p^{mn}}\sum_{\mathbf b\in\adm_p(m)}
\mathfrak S_p(\mathbf b)T_p(g,\mathbf b)
=1+O\left(\frac{k^{5k}}{p^{2}}\right).
\]
\end{lemma}

\begin{proof}
By Lemma~\ref{lDensityLift} and admissibility, the sum equals
\[
\left(1-\frac{1}{p}\right)^k\sum_{0\leq J<m}
\int_{\substack{\mathbf t\in\Z_p^n\textup{ prim.}\\
\sum_iv_p(t_i)=J}}
\prod_{i=1}^k\tau_p(Q_{i,\mathbf t})
T_p(g,\mathbf t)\mathrm{d}\mathbf{t}.
\]
If $J\in\{0,1\}$, then $T_p(g, \mathbf{t}) = 1$ by Lemma~\ref{lemma:tp1}. By Lemma~\ref{lem:taupestimate}, the integral over the case $J\in\{0,1\}$ equals $1-\frac{3k}{p}+\frac{n}{p}+O(\frac{k^2}{p^2}) = 1+\frac{k}{p}+O(\frac{k^2}{p^2})$; multiplying by $(1-\frac{1}{p})^k = 1-\frac{k}{p}+O(\frac{k^2}{p^2})$, the contribution of this case is $1+O(\frac{k^2}{p^2})$. For $J\geq2$, use $|T_p(g,\mathbf t)|\leq1$ and Lemma~\ref{lem:taupestimate}. This gives an error $O(k^{5k}p^{-2})$. 
\end{proof}

Inspired by the results obtained thus far, we introduce the set $\mathcal{U}_{\mathcal{E},m}(B; V_k)$ (which are exactly those rational points we intend to use for our lower bound) given by
\[ 
\mathcal{U}_{\mathcal{E},m}(B; V_k)
\coloneqq
\left\{\begin{array}{c} \mathbf{t} = (t_1, \dots, t_n), \mathbf{x}_1 = (X_1, Y_1, Z_1), \dots, \mathbf{x}_k = (X_k, Y_k, Z_k) \in V_k^{\mathrm{prim}} : \\
\gamma_n^{j - 1} \leq |t_i| <\gamma_n^{j}\textup{ for all }i\in[n]\textup{ for some }h/2 \leq j \leq h,\\
\mathbf{t}\bmod C(\mathcal{E},m)\in\adm(\mathcal{E},m),\ 
p^2 \nmid t_1 \cdots t_n\textup{ for all }p> \mathcal{E},\\
B^{\frac{1}{4n}}\leq H_{i, \mathbf{t}}(\mathbf{x}_i) \leq B^{\frac{1}{2n}}\textup{ for }i < k,\ 
H_{k, \mathbf{t}}(\mathbf{x}_k) \leq \frac{B}{|\mathbf{t}|^{2k} \prod_{i = 1}^{k - 1} H_{i, \mathbf{t}}(\mathbf{x}_i)}
\end{array}
\right\}.
\] 
Also introduce
\[ 
\mathcal{U}'_{\mathcal{E},m}(B; V_k)
\coloneqq
\left\{(\mathbf{t},\mathbf{x})\in\mathcal{U}_{\mathcal{E},m}(B; V_k):
\begin{array}{c} 
t_i>0\textup{ for all }i\in[n],\ \mathbf{t}\equiv \mathbf{1}\bmod 8,\\ \mathbf{t}\bmod p^m\in((\Z/p^m\Z)^{\times})^n\textup{ for all }2 < p \leq k^{12k}
\end{array}
\right\}.
\] 
For $k \geq 2$, $\mathcal{E} > k^{12k}$ and $m \geq 3$ sufficiently large relative to $k$, and $B$ sufficiently large in terms of $k, \mathcal{E}, m$, applying Lemma~\ref{lemma:stepU} and summing over the relevant $(\boldsymbol{\varepsilon}, \mathbf{b}, \boldsymbol{\upsilon})$ gives
\begin{equation} \label{eq:preUu}
\begin{split}
k^{-c_3k} B (\log B)^k \times S_{\mathcal{E},m}(V_k) &\leq \# \mathcal{U}_{\mathcal{E},m}(B; V_k) \leq k^{c_3k} B (\log B)^k \times S_{\mathcal{E},m}(V_k), \\
k^{-c_3k} B (\log B)^k \times S'_{\mathcal{E},m}(V_k) &\leq \# \mathcal{U}'_{\mathcal{E},m}(B; V_k) \leq k^{c_3k} B (\log B)^k \times S'_{\mathcal{E},m}(V_k),
\end{split}
\end{equation}
where 
\begin{equation}
\label{eq:defSem}
S_{\mathcal{E},m}(V_k) \coloneqq
\sum_{g \in \PBM} \sum_{\boldsymbol{\varepsilon} \in \{\pm 1\}^n} T_\infty(g, \boldsymbol{\varepsilon}) \prod_{p \leq \mathcal{E}}\sum_{\mathbf{b} \in \adm_p(m)} \frac{ \mathfrak{S}_p(\mathbf{b}) T_p(g, \mathbf{b})}{p^{mn}},
\end{equation}
and
\begin{multline}
\label{eq:defSemprime}
S'_{\mathcal{E},m}(V_k) \coloneqq
\frac{1}{2^{mn}} \sum_{\substack{\mathbf{b}\in\adm_2(m)\\ \mathbf{b}\equiv\mathbf{1}\bmod 8}}\mathfrak{S}_2(\mathbf{b}) \times \prod_{2<p \leq k^{12k}} \frac{1}{p^{mn}} \sum_{\substack{\mathbf{b}\in\adm_p(m)\\ \mathbf{b}\in((\Z/p^m\Z)^{\times})^n}}\mathfrak{S}_p(\mathbf{b}) \\
{}\times \sum_{g \in \PBM}\prod_{ k^{12k}<p \leq \mathcal{E}} \frac{1}{p^{mn}} \sum_{\mathbf{b}\in\adm_p(m)}\mathfrak{S}_p(\mathbf{b}) T_p(g, \mathbf{b}).
\end{multline}

\begin{lemma}
\label{lemma:strictlypositive}
For sufficiently large $k$, every
$\mathcal E>k^{12k}$ and $m\geq3$, we have
\[
S_{\mathcal{E},m}(V_k)\geq S'_{\mathcal{E},m}(V_k)
\geq 2^{\frac{1}{2}k^2 - c_4 k \log k},
\]
where $c_4>0$ is an absolute constant. Moreover, we have
\[
S'_{\mathcal{E},m}(V_k)
\geq 2^{\frac{1}{2}k(k + 1) - c_4} 
\frac{\mathfrak{S}_2(\mathbf{1})}{2^{3n}} \prod_{2 < p \leq k^{12 k}} \sum_{\substack{\mathbf{b} \in \adm_p(m) \\ \mathbf{b} \in ((\Z/p^m\Z)^{\times})^n}} \frac{\mathfrak{S}_p(\mathbf{b})}{p^{mn}}
\prod_{p > k^{12k}} \sum_{\mathbf{b} \in \adm_p(m)} \frac{\mathfrak{S}_p(\mathbf{b})}{p^{mn}}.
\]
\end{lemma}

\begin{proof}
By Lemma~\ref{lemma:charorth}, we may freely drop congruence classes in order to lower bound $S_{\mathcal{E},m}(V_k)$. We restrict the sum $S_{\mathcal{E},m}(V_k)$ to $\boldsymbol{\varepsilon}=\mathbf{1}$, $\mathbf{b} \in \adm(\mathcal{E},m)$, $\mathbf{b}\equiv \mathbf{1}\bmod 8$, and $\mathbf{b}\bmod p^m\in((\Z/p^m\Z)^{\times})^n $ for all $2<p\leq k^{12k}$. Then
\[S_{\mathcal{E},m}(V_k)\geq S'_{\mathcal{E},m}(V_k).\]

We first show the second claim.
For the terms with $ p> k^{12k}$, Lemma~\ref{lemma:sigmahat} implies that
\begin{equation}\label{eq:largeprimel}\frac{1}{p^{mn}} \sum_{\mathbf{b} \in \adm_p(m)} \mathfrak{S}_p(\mathbf{b})T_p(g, \mathbf{b})
=1+O\left(\frac{k^{5k}}{p^{2}}\right).\end{equation}
The error term is $O(k^{-19k})$ for $p>k^{12k}$, so for sufficiently large $k$, we may assume that the right hand side of \eqref{eq:largeprimel} is at least $1/2$.
For the terms with $k^{12k} < p \leq \mathcal{E}$ in \eqref{eq:defSemprime}, we apply \eqref{eq:largeprimel} twice (once with an arbitrary $g$ and once with $g = 0$, in which case $T_p(0, \mathbf{b}) = 1$). This gives
\[
\frac{1}{p^{mn}} \sum_{\mathbf{b} \in \adm_p(m)} \mathfrak{S}_p(\mathbf{b})T_p(g, \mathbf{b})
= \left(\frac{1}{p^{mn}} \sum_{\mathbf{b} \in \adm_p(m)} \mathfrak{S}_p(\mathbf{b})\right) \left(1+O\left(\frac{k^{5k}}{p^{2}}\right)\right).
\]
For the terms with $p> \mathcal{E}$, apply \eqref{eq:largeprimel} with $g = 0$ to obtain
\[
\frac{1}{p^{mn}} \sum_{\mathbf{b} \in \adm_p(m)} \mathfrak{S}_p(\mathbf{b})
=1+O\left(\frac{k^{5k}}{p^{2}}\right).
\]
The product of the multiplicative error terms is
\begin{equation}\label{eq:multer}
\prod_{ p>k^{12k}}
\left(1+O\left(\frac{k^{5k}}{p^{2}}\right)\right)
=\exp\left(O\left((k^{7k+1}\log k)^{-1}\right)\right).
\end{equation}
This yields
\begin{multline*}
S'_{\mathcal{E},m}(V_k) \gg
 \frac{1}{2^{mn}} \sum_{\substack{\mathbf{b}\in\adm_2(m)\\ \mathbf{b}\equiv\mathbf{1}\bmod 8}}\mathfrak{S}_2(\mathbf{b}) \times \prod_{2<p \leq k^{12k}} \frac{1}{p^{mn}} \sum_{\substack{\mathbf{b}\in\adm_p(m)\\ \mathbf{b}\in((\Z/p^m\Z)^{\times})^n}}\mathfrak{S}_p(\mathbf{b}) \\
{}\times \sum_{g \in \PBM}\prod_{p > k^{12k}} \frac{1}{p^{mn}} \sum_{\mathbf{b}\in\adm_p(m)}\mathfrak{S}_p(\mathbf{b}).
\end{multline*}
Since $\lvert\PBM\rvert = 2^{\frac{1}{2}k(k+1)}$, the second part follows.

We now deduce the first part of the lemma. 
For the primes up to $k^{12k}$, we evaluate the sums over odd primes using Lemma~\ref{lSmoothConic} and the $2$-adic sum directly. For $p>k^{12k}$, we replace the sum by \eqref{eq:largeprimel}, which collectively introduces an error bounded by \eqref{eq:multer}. Since $\lvert\PBM\rvert=2^{\frac12k(k+1)}$, we have
\[
S'_{\mathcal{E},m}(V_k)\gg 2^{\frac{1}{2}k(k+1)}\times 2^{-12k} \prod_{2<p \leq k^{12k}} \left( 1 - \frac{1}{p^2} \right)^k \left( 1 - \frac{1}{p} \right)^{4k}
\gg \frac{2^{\frac{1}{2}k(k+1)-c_4k}}{(\log k^{12k})^{4k}},
\]
where the implied constant and $c_4>0$ are absolute. The result follows from enlarging $c_4$.
\end{proof}

Combining \eqref{eq:preUu} and Lemma~\ref{lemma:strictlypositive}, we deduce the following lower bound.

\begin{theorem}
\label{tLowerBound}
There exist constants $c_3,c_4>0$ such that for all sufficiently large $k$, for all $\mathcal{E}$ and $m$ sufficiently large in terms of $k$, and for all $B \geq 100$ sufficiently large in terms of $k$, $\mathcal{E}$ and $m$, we have
\[
\# \mathcal{U}_{\mathcal{E},m}(B; V_k)\geq k^{-c_3k} B (\log B)^k\times S_{\mathcal{E},m}(V_k),
\]
and
\[
\# \mathcal{U}'_{\mathcal{E},m}(B; V_k)\geq k^{-c_3k} B (\log B)^k\times S'_{\mathcal{E},m}(V_k),
\]
where $S_{\mathcal{E},m}(V_k)$ and $S'_{\mathcal{E},m}(V_k)$, defined in~\eqref{eq:defSem} and~\eqref{eq:defSemprime}, satisfy
\[
S_{\mathcal{E},m}(V_k) \geq S'_{\mathcal{E},m}(V_k) \geq 2^{\frac{1}{2}k^2-c_4k\log k},
\]
and
\[
S'_{\mathcal{E},m}(V_k)
\geq 2^{\frac{1}{2}k(k + 1) - c_4} 
\frac{\mathfrak{S}_2(\mathbf{1})}{2^{3n}} \prod_{2 < p \leq k^{12 k}} \sum_{\substack{\mathbf{b} \in \adm_p(m) \\ \mathbf{b} \in ((\Z/p^m\Z)^{\times})^n}} \frac{\mathfrak{S}_p(\mathbf{b})}{p^{mn}}
\prod_{p > k^{12k}} \sum_{\mathbf{b} \in \adm_p(m)} \frac{\mathfrak{S}_p(\mathbf{b})}{p^{mn}}.
\]
\end{theorem}

\section{Geometry}\label{s:geometry}
The goal of this section is to check that the geometry of $V_k$ is compatible with the assumptions of \cite{BT98}. We will assume that $k \geq 2$. 

\begin{lemma}
\label{lem:Hp(X,OX)}
	For a resolution of singularities $X \to V_k$ we have $\HH^p(X, \mathcal{O}_X) = 0$ for all $p > 0$.
\end{lemma}
\begin{proof}
	This follows immediately from the rationality of $V_k$ in Lemma \ref{lem:V_k_rational}.
\end{proof}

It follows from \cite[Theorem 2.5]{ChambertLoirTschinkel2010} that global Tamagawa measures are well-defined on any smooth open $U \subseteq V_k$. We will use this freely in this section.

\subsection{Singularities}
The variety $V_k$ is singular. In this case the analogue of Peyre's constant was defined in \cite{BT98}.
We start by recalling definitions of certain classes of singularities. Let $f: X \to V_k$ be a resolution of singularities, which we may assume to be an isomorphism on the smooth locus $V_k^{\mathrm{sm}}$ and such that $D \coloneqq X \setminus V_k^{\mathrm{sm}}$ has simple normal crossings. The complement $(X \setminus V_k^{\mathrm{sm}})_{\overline{\Q}}$ is a union of irreducible divisors $D_i$ for a finite index set $I$. We have $K_X = f^* K_{V_k} + \sum_{i \in I} a_i D_i$ with $a_i \in \Q$. The rational numbers $a_i$ are called the \emph{discrepancies} of the resolution. Recall that $V_k$ is \emph{canonical} if $a_i \geq 0$ for all $i$ and \emph{terminal} if $a_i > 0$ for all $i$. If $V_k$ is canonical and $a_i = 0$ then $D_i$ is a \emph{crepant divisor}. 

Recall that we fixed an adelic metric on $K_{V_k}$. This induces an adelic metric on $K_X - \sum_{i \in I} a_i D_i$ by pullback. Choose an adelic metric on $\sum_{i \in I} a_i D_i$ and for each prime $p$ let $\tau_{X, p}$ be the induced Tamagawa measure on $X(\Q_p)$. 

We will now put the following additional conditions on $p$. Note that a positive density of primes satisfy these conditions by spreading out, the Chebotarev density theorem and Lang--Weil.
\begin{definition}
	The prime $p$ is \emph{split and good} if the following conditions hold:
	\begin{itemize}
		\item All divisors $D_i$ are defined over $\Q_p$.
		\item There exists a smooth proper $\Z_p$-integral model $\mathcal{X}$ of $X$ such that the $p$-adic metric on $K_X$ is equal to the one defined by this model.
		\item The closure of $\mathcal{D} \subset \mathcal{X}$ of $D$ has simple normal crossings over $\Z_p$ and every non-empty open stratum has a $\Q_p$-point.
		\item If $\mathcal{D}_i$ is the closure of $D_i$ then the metric on $\sum_{i \in I} a_i D_i$ is the one induced by these integral models $\mathcal{D}_i$.
	\end{itemize}
\end{definition}

If $p$ is a split and good prime then for all $i$ let $\lVert f_{\mathcal{D}_i}\rVert_p: X(\Q_p) \to \R_{\geq 0}$ be the map obtained by evaluating the model metric of $\mathcal{O}(D_i)$ at a section $f_{\mathcal{D}_i} \in \HH^0(\mathcal{X}, \mathcal{O}(D_i))$ defining $\mathcal{D}_i$. This function is well-defined as $f_{\mathcal{D}_i}$ is unique up to multiplication by an element of $\Z_p^{\times}$.

\begin{lemma}
\label{lem:p-adic_density_on_X}
	Assume that the prime $p$ is split and good, then
	\[
	\tau_p(V_k) = \int_{V_k^{\mathrm{sm}}(\Q_p)} \mathrm{d}\tau_{V_k, p}(x) = \int_{X(\Q_p)} \prod_{i \in I} \lVert f_{\mathcal{D}_i}\rVert_p(x)^{a_i} \mathrm{d}\tau_{X, p}(x).
	\]
\end{lemma}

\begin{proof}
	The function $x \mapsto \prod_{i \in I} \lVert f_{\mathcal{D}_i}\rVert_p(x)^{a_i}$ is by definition the value at $x$ of the $p$-adic metric of the section $\prod_{i \in I}f_{\mathcal{D}_i}^{a_i} \in \mathcal{O}(\sum_{i \in I} a_i D_i)$ which defines the $\Q$-divisor $\sum_{i \in I} a_i D_i$. This equality is then by definition of the local Tamagawa measure, see also \cite[\S2.4.3]{ChambertLoirTschinkel2010}.
\end{proof}
\begin{lemma}
\label{lem:Vk_terminal}
	The variety $V_k$ has terminal singularities. 
\end{lemma}

\begin{proof}
	By spreading out there exists an $N \in \Z_{> 0}$ and an $\Z[N^{-1}]$-integral model $\mathcal{X} \to \mathcal{V}_k$ of $X \to V_k$ with the following properties: $\mathcal{X}$ is smooth over $\Z[N^{-1}]$, $\mathcal{X}$ and $\mathcal{V}_k$ are proper over $\Z[N^{-1}]$ and the closure $\mathcal{D} \subset \mathcal{X}$ of $D$ is a simple normal crossing divisor.
	Let $p$ be a split and good prime for the given integral model $\mathcal{X}_{\Z_p}$. A positive density of primes have these properties. For each subset $J \subseteq I$ we let $\mathcal{D}_J \coloneqq \bigcap_{j \in J} \mathcal{D}_j$ and $\mathcal{D}_J^{\circ} \subseteq \mathcal{D}_J$ the complement of all $\mathcal{D}_{J'}$ with $J \subset J'$.
	
	We can then apply Denef's formula \cite[Proposition~4.5]{ChambertLoirTschinkel2010} to compute the right hand side of Lemma \ref{lem:p-adic_density_on_X}. It implies that $\tau_p(V_k) < \infty$, which holds by Theorem \ref{thm:local_density_total_variety}, only if $a_i > -1$ for all $i \in I$ and if this condition holds then 
	\[
	\tau_p(V_{k}) = \frac{\#\mathcal{V}^{\mathrm{sm}}(\FF_p)}{p^{\dim V_k}} + \frac{1}{p^{\dim V_k}} \sum_{\varnothing \neq J \subseteq I} \# \mathcal{D}_J^{\circ}(\FF_p) \prod_{i \in J}\frac{1 - \frac{1}{p}}{p^{a_i} - \frac{1}{p}}.
	\]
	A standard application of the Lefschetz trace formula, Lemma \ref{lem:Hp(X,OX)} and the Riemann hypothesis for varieties (e.g.~as in the proof of \cite[Theorem 3.3.8]{BT98}) shows that for a prime $\ell \neq p$ we have
	\[
	\frac{\#\mathcal{V}^{\mathrm{sm}}(\FF_p)}{p^{\dim V_k}} = 1 + \frac{\Tr(F_p : \Pic \,(V_k^{\mathrm{sm}})_{\overline{\Q}} \otimes \Q_{\ell} )}{p} + O\left(\frac{1}{p^{\frac{3}{2}}} \right),
	\]
	where $F_p$ is the Frobenius element.
	The Picard group $\Pic \, (V_k^{\mathrm{sm}})_{\overline{\Q}}$ has rank $k + 1$ and a trivial Galois action by Lemma \ref{lem:Picard_group_computation}, so the trace is equal to $k + 1$.
	
	For $\varnothing \neq J \subseteq I$ we can use the Lang--Weil bounds to see that $\# \mathcal{D}_J(\FF_p) \ll_J p^{\dim V_k - \# J}$. Moreover, if $\# J = 1$ then $\mathcal{D}_J$ is absolutely irreducible by assumption so we get $\# \mathcal{D}_J(\FF_p) = p^{\dim V_k - 1}(1 + O(p^{-\frac{1}{2}}))$.
	
	Putting everything together we find that 
	\[
	\tau_p(V_k) = 1 + \frac{k + 1}{p} + O\left(\frac{1}{p^{\frac{3}{2}}}\right) + \frac{B}{p^{1 + A}}+ O\left(\frac{1}{p^{1 + A + \delta}}\right).
	\]
	The constants are given by $A = \min_{i \in I} a_i, B = \#\{i \in I: a_i = A\}$ and $\delta > 0$.
	By comparing this with our computation of $\tau_p(V_k)$ in Theorem \ref{thm:local_density_total_variety} we deduce that $A \geq \frac{1}{2}$, which proves the lemma.
\end{proof}

\begin{remark}
	We wrote that the discrepancies $a_i$ are rational numbers since this is the general situation considered in the study of singularities where $K_{V_k}$ might only be $\Q$-Cartier. But note that since $K_{V_k}$ is Cartier the discrepancies $a_i$ are integers. The fact that $V_k$ is terminal thus implies that we have $a_i \geq 1$. This implies that the singularities of $V_k$ are quite mild.
\end{remark}

\subsection{Brauer group}
The goal of this section is to compute the Brauer group of $V_k^{\mathrm{sm}}$ and relate it to the main terms which appear in the character sums. We will do this by reducing this to the computation of the relevant \emph{subordinate Brauer group}, which was determined in \cite[\S5.2]{CKR}. Let us recall the relevant definitions of \cite[\S2]{LRS}. Subordinate Brauer groups are only defined for proper varieties in \cite{LRS}, it will be convenient to allow non-proper varieties. We begin by recalling the basic theory of Brauer groups.

The Brauer group of a smooth variety $X$ over a field $K$ of characteristic $0$ is the group $\Br X \coloneqq \HH^2_{\textup{\'et}}(X, \mathbb{G}_m)$, where the cohomology theory is \'etale cohomology. There is an inclusion $\Br X \subseteq \Br K(X)$ and for every divisor $D \subset X$ there is a \emph{residue map} $\partial_D: \Br K(X) \to \HH^1(K(D), \Q/\Z)$. We have the following fundamental exact sequence, known as \emph{Grothendieck purity} \cite[Theorem 3.7.3]{ColliotTheleneSkorobogatov2021}:
\[
	0 \to \Br X \to \Br K(X) \to \bigoplus_{\substack{D \subset X \\ \textup{divisor}}} \HH^1(K(D), \Q/\Z).
\]

\begin{definition}
	Let $f:Y \to X$ be a dominant map of smooth varieties over $K$. We define the \emph{subordinate Brauer group}
	\[
	\Br_{\mathrm{sub}}(Y/X) \coloneqq \left\{b \in \Br K(X) \hspace{-0.05cm} : \hspace{-0.05cm}
	\begin{matrix}
		\partial_D(b) \in \ker\left(\HH^1(K(D), \Q/\Z) \to \HH^1(K(E), \Q/\Z)\right) \textup{ for all } \\
		 D \subset X \textup{ and } E \subseteq f^{-1}(D) \textup{ of multiplicity } 1 \textup{ dominant over } D
	\end{matrix}
	\right\}.
	\]
\end{definition}

Note that this definition is invariant when modifying $Y$ or $X$ in codimension at least $2$ by \cite[Theorem 3.7.2]{ColliotTheleneSkorobogatov2021}.

We also need to recall the following definition \cite[Definition 2.1]{Loughran2018Varieties}.
\begin{definition}
	Let $R$ be a DVR with perfect residue ring and $Y$ be a finite type separated scheme over $R$ with smooth locus $Y^{\mathrm{sm}}$. We say that $Y$ is almost smooth if its generic fibre is smooth and for all \'etale $R$-algebras $R'$ we have $Y(R') = Y^{\mathrm{sm}}(R')$.
	
	If $X$ is a normal scheme then a map $f: Y \to X$ is almost smooth if it is almost smooth over the local ring of every codimension $1$ point.
\end{definition}
This definition is useful because of the following lemmas.

\begin{lemma}
	The projection map $\pi: V_k \to \PPP^{n - 1}$ is almost smooth.
\end{lemma}

\begin{proof}
	The fibre of a codimension $1$ point is either a product of smooth conics or a product of degenerate conics whose discriminant has valuation $1$ with respect to the valuation induced by the point.
	
	 A product of almost smooth varieties over a DVR is almost smooth, directly by the definition. It thus suffices to show that conics of the above type are almost smooth, which is true by \cite[Proposition 3.1.2]{Bosch} since such conics have a regular total space.
\end{proof}

The following lemma generalises \cite[Lemma 2.5]{LRS}.

\begin{lemma}\label{lem:almost_smooth_birational_subordinate_ggroup}
	Let $X \to B, Y \to B$ be dominant proper morphisms of varieties over $K$ with geometrically integral generic fibres. Assume that $X$ and $Y$ are almost smooth over $B$. Let $X^{\mathrm{sm}} \subset X, Y^{\mathrm{sm}} \subset Y$ be the smooth loci of the varieties.
	
	If $h: Y \dashrightarrow X$ is a birational map over $B$ then pullback by $h$ defines an isomorphism $\Br_{\mathrm{sub}}(X^{\mathrm{sm}}/B) \cong \Br_{\mathrm{sub}}(Y^{\mathrm{sm}}/B)$.
\end{lemma}

\begin{proof}
	This is immediate from the definition of the subordinate Brauer group and \cite[Lemma 2.2]{Loughran2018Varieties} applied to the fibre over the local ring of each codimension $1$ point of $B$. The conditions in the lemma hold since a dominant map of varieties $X \to B$ with $B$ normal is flat in codimension $1$ \cite[\href{https://stacks.math.columbia.edu/tag/0AUW}{Tag 0AUW}]{Stacks}.
\end{proof}

Let $f:X \to V_k$ be a resolution of singularities, which we assume to be an isomorphism over $U^{\mathrm{top}}$. Consider the map $\pi: U^{\mathrm{top}} \to U^{\mathrm{bot}}$. 

\begin{lemma}
\label{lem:BrSub=BrVert}
	The map $\pi^*:\Br_{\mathrm{sub}}(X/\PPP^{n - 1}) = \Br_{\mathrm{sub}}(U^{\mathrm{top}}/U^{\mathrm{bot}}) \to \Br V_k^{\mathrm{sm}}$ is surjective. 
\end{lemma}
\begin{proof}
	We have $\Br_{\mathrm{sub}}(X/\PPP^{n - 1}) \cong \Br_{\mathrm{sub}}(V_k^{\mathrm{sm}}/\PPP^{n - 1})$ by Lemma \ref{lem:almost_smooth_birational_subordinate_ggroup}. This is isomorphic to $\Br_{\mathrm{sub}}(U^{\mathrm{top}}/U^{\mathrm{bot}})$ since the complements of $U^{\mathrm{top}}$ and $U^{\mathrm{bot}}$ have codimension $2$ by Lemma \ref{lem:properties_U_k}.
	
	Let $\eta$ be the generic point of $\PPP^{n-1}$. The generic fibre $\pi^{-1}(\eta)$ is a product of $k$ smooth conics so the map $\pi^*: \Br \Q(\PPP^{n-1}) \to \Br \pi^{-1}(\eta)$ is surjective by combining \cite[Proposition~7.2.4(a), Proposition~11.2.4]{ColliotTheleneSkorobogatov2021}.
	It follows from purity that $\Br V_k^{\mathrm{sm}} \cong \Br U^{\mathrm{top}} \subseteq \Br \pi^{-1}(\eta) \subseteq \pi^*\Br\Q(\PPP^{n-1})$. Thus $\Br U^{\mathrm{top}}$ is equal to the \emph{vertical Brauer group} $\Br_{\mathrm{vert}}(U^{\mathrm{top}}/U^{\mathrm{bot}})$, see \cite[Definition 11.1.1]{ColliotTheleneSkorobogatov2021} for the definition. 
	
	It remains to show the equality $\pi^*\Br_{\mathrm{sub}}(U^{\mathrm{top}}/U^{\mathrm{bot}}) = \Br_{\mathrm{vert}}(U^{\mathrm{top}}/U^{\mathrm{bot}})$. This is true by \cite[Lemma 2.3]{LRS} as $U^{\mathrm{top}} \to U^{\mathrm{bot}}$ is a flat map of smooth varieties by Lemma \ref{lem:properties_U_k} and every fibre above a codimension one point is a product of reduced conics and thus reduced. We note that the cited lemma assumes implicitly that $U^{\mathrm{top}}$ and $U^{\mathrm{bot}}$ are proper, but that the proof never actually uses the properness.
\end{proof}

\begin{remark}
	With more work one can show that the kernel has size $2^k$ and is generated by the classes $[Q_i] \in \Br\Q(\PPP^{n-1})$, where $Q_i$ is the conic defined by the $i$-th equation in \eqref{eQuadForms}. 
\end{remark}

We can now use the computation of the subordinate Brauer group in \cite[Lemma 5.4]{CKR}. Recall the following notation from the reference. 

Recall the group $\PBM$ from Definition \ref{def:Ralt_PBM}. A bijective mapping between $\PBM$ and $\Br_{\mathrm{sub}}(X/\PPP^{n - 1})/\Br \Q$ is constructed in \cite[Lemma 5.4]{CKR} (the equality between what is denoted $\Br_{\mathrm{sub}}(\PPP^{n-1}, \mathcal{B}) $ in the reference and what we denote $\Br_{\mathrm{sub}}(X/\PPP^{n - 1})$ is \cite[Proposition~2.9]{LRS}). We recall the definition of this map. For $g \in \PBM$ the corresponding Brauer element is given by 
\[
	\alpha_g \coloneqq \sum_{\{i,j\} \subseteq [n-1]} \langle \{i\}, g(j) \rangle (t_i/t_n, t_j/t_n) + \sum_{\ell \in [n-1]} \langle \{-\}, g(\ell) \rangle (-1, t_\ell/t_n) \in \Br \Q(\PPP^{n-1}),
\]
where $\langle S, T \rangle$ is $|S \cap T|$ taken modulo $2$. Combining this with Lemma~\ref{lem:BrSub=BrVert}, we deduce the following proposition.

\begin{prop}
\label{prop:Brauer_group=RALt}
	The map $\PBM \to \Br V_{k}^{\mathrm{sm}}/\Br \Q: g \mapsto \pi^* \alpha_g \bmod \Br \Q$ is a surjective group homomorphism.
\end{prop}

\subsubsection{Volume of Brauer--Manin subset}
Recall \cite[\S13.3]{ColliotTheleneSkorobogatov2021} that for any variety $X$ over $\Q$ there is a pairing $ X(\Adele_{\Q}) \times \Br X \to \Q/\Z$, the Brauer--Manin pairing, whose left kernel is the \emph{Brauer--Manin set} $X(\Adele_{\Q})^{\Br}$. More generally, for any subset $\Omega \subseteq X(\Adele_{\Q})$ we write $\Omega^{\Br} \coloneqq \Omega \cap X(\Adele_{\Q})^{\Br}$. 

We will compute the volume of the Brauer--Manin set in terms of a finite sum of Euler products, we will show that these Euler products are the same Euler products appearing in the character sum. There will be one issue at $\infty$. This is because the points counted in the character sum are only those points contained in a small subset of $V_k(\R)$. At $\infty$ we will thus have to use an upper bound. The following decomposition of the adelic space will be useful.
\begin{definition}
	For $\boldsymbol{\varepsilon} = (\varepsilon_i)_{i \in[n]} \in \{\pm 1\}^n$ define the subspaces
	\begin{align*}
		\PPP^{n - 1}(\R)_{\boldsymbol{\varepsilon}} &\coloneqq \{\mathbf{t} \in \PPP^{n - 1}(\R): \varepsilon_i t_i > 0 \textup{ for all }i \in[n]\}, \\
		V_k^{\mathrm{sm}}(\R)_{\boldsymbol{\varepsilon}} &\coloneqq \{x \in V_k^{\mathrm{sm}}(\R): \pi(x) \in \PPP^{n - 1}(\R)_{\boldsymbol{\varepsilon}}\}, \\
		V_k^{\mathrm{sm}}(\Adele_{\Q})_{\boldsymbol{\varepsilon}} &\coloneqq \{(P_v)_{v} \in V_k^{\mathrm{sm}}(\Adele_{\Q}): P_{\infty} \in V_k^{\mathrm{sm}}(\R)_{\boldsymbol{\varepsilon}}\}.
	\end{align*}
\end{definition}
The condition $\varepsilon_i t_i > 0$ might seem a priori not well-defined. We mean that all of the inequalities hold for some choice of representative $\mathbf{t}$, in particular $\PPP^{n - 1}(\R)_{\boldsymbol{\varepsilon}} = \PPP^{n - 1}(\R)_{-\boldsymbol{\varepsilon}}$.

For each prime $p$ and each $\boldsymbol{\varepsilon}\in\{\pm 1\}^n$, define local integrals
\begin{align*}
\hat{\tau}_p(g)&\coloneqq \int_{V_k^{\mathrm{sm}}(\Q_p)} e(\inv_p(\pi^*\alpha_g(x_p))) \mathrm{d}\tau_{V_k, p}(x_p),\\
\hat{\tau}_{\infty}(g; \boldsymbol{\varepsilon})&\coloneqq \int_{V_k^{\mathrm{sm}}(\R)_{\boldsymbol{\varepsilon}}} e(\inv_\infty(\pi^*\alpha_g(x_\infty))) \mathrm{d}\tau_{V_k, \infty}(x_\infty),
\end{align*} 
where $e(z) = e^{2 \pi i z}$. These integrals converge by Corollary \ref{cor:Tamagawa_measure_converges}. Henceforth we write $\tau_{V_k} \coloneqq \tau_{V_k^{\mathrm{sm}}}$ for the global Tamagawa measure on $V_k^{\mathrm{sm}}(\Adele_\Q)$ induced by the adelic metric of Section~\ref{sec:adelic_metrics} (well-defined by Lemmas~\ref{lem:Hp(X,OX)} and~\ref{lem:Vk_terminal} and Proposition~\ref{prop:Tamagawa_measure_terminal}).
\begin{lemma}\label{lem:BM_set_character_orthogonality}
	For all $\boldsymbol{\varepsilon} \in \{\pm 1\}^n$, we have 
	\[
	\tau_{V_k}\left(V_k^{\mathrm{sm}}(\Adele_\Q)_{\boldsymbol{\varepsilon}}^{\Br} \right) = \frac{1}{\lvert\PBM\rvert}\sum_{g \in \PBM} \hat{\tau}_{\infty}(g; \boldsymbol{\varepsilon})\prod_p \lambda_p \hat{\tau}_p(g),
	\]
	with convergence factors given by $\lambda_p = (1 - \frac{1}{p})^{k + 1}$. Each Euler product on the right-hand side converges absolutely.
\end{lemma}

\begin{proof}
	This follows immediately from Proposition \ref{prop:Brauer_group=RALt} and orthogonality of characters as in \cite[Lemma 3.13]{LRS}. To show convergence of each Euler product we rewrite it as
	\[
		\hat{\tau}_{\infty}(g; \boldsymbol{\varepsilon})\prod_p \lambda_p \hat{\tau}_p(g) = \int_{V_k^{\mathrm{sm}}(\Adele_{\Q})_{\boldsymbol{\varepsilon}}} e\Big(\sum_v \inv_v(\pi^*\alpha_g(x_v))\Big) \mathrm{d}\tau_{V_k}(x).
	\] 
	This integral converges absolutely since $\tau_{V_k}\left(V_k^{\mathrm{sm}}(\Adele_\Q)_{\boldsymbol{\varepsilon}} \right) < \infty$ by Corollary \ref{cor:Tamagawa_measure_converges} and $x \mapsto \sum_v \inv_v(\pi^*\alpha_g(x_v))$ is a continuous function on $V_k^{\mathrm{sm}}(\Adele_{\Q})$ because it is just the Brauer--Manin pairing, which is locally constant.
\end{proof}

We will now compute the integrals $\hat{\tau}_p(g)$ and $\hat{\tau}_{\infty}(g; \boldsymbol{\varepsilon})$ appearing in the above lemma. Recall \eqref{eq:Tvdef}. 
\begin{lemma}
\label{lem:local_BM_equals_T(g,a)}
	For all $g \in \PBM$ and $\mathbf{a} \in (\Z_p \setminus \{0\})^n$ we have $e(\inv_p(\alpha_g(\mathbf{a}))) = T_p(g, \mathbf{a})$.
\end{lemma}

\begin{proof}
	This is \cite[Lemma 5.7]{CKR}.
\end{proof}

\begin{lemma}
\label{lem:Brauer_transform_at_p}
	For all primes $p$ we have
	\begin{align*}
	\hat{\tau}_p(g)
&=
\left(1-\frac1p\right)^{-1}
\int_{\substack{\mathbf t\in\Z_p^n\\
\mathbf{t}\textup{ prim.}}}
\prod_{i=1}^k\tau_p(Q_{i,\mathbf t})
T_p(g,\mathbf t)\mathrm{d}\mathbf{t}\\
&=
\lim_{m\to\infty}
\frac{(1-\frac1p)^{-k-1}}{p^{mn}}
\sum_{\mathbf b\in\adm_p(m)}
\mathfrak S_p(\mathbf b)T_p(g,\mathbf b).
\end{align*}
	In particular, the limit on the right-hand side exists.
\end{lemma}

\begin{proof}
	We have $\pi^*\alpha_g(x_p) = \alpha_g(\pi(x_p))$. The function $x_p \mapsto e(\inv_p(\pi^*\alpha_g(x_p)))$ is $L^1$ as $\tau_{V_k, p}(V_k^{\mathrm{sm}}(\Q_p)) < \infty$ by Corollary \ref{cor:Tamagawa_measure_converges}. By integrating over fibres and lifting the Tamagawa measure on $\PPP^{n-1}$ to the affine cone $\AAA^n \setminus \{0\}$ as in the proof of Lemma \ref{lem:integrate_fibers_affine} we have
\[
	\hat{\tau}_p(g) = \left(1 - \frac{1}{p}\right)^{-1} \int_{\substack{\mathbf{t} \in \Z_p^n \\ \mathbf{t} \textup{ prim.}}} e(\inv_p(\alpha_g(\mathbf{t}))) \prod_{i = 1}^k \tau_p(Q_{i, \mathbf{t}}) \mathrm{d} \mathbf{t}.
\]
	For each $m$ we now restrict the integral to the set of primitive $\mathbf{t} \in \Z_p^n$ such that $\mathbf{t} \bmod p^m \in \adm_p(m)$. These restricted integrals converge to the total integral by the dominated convergence theorem as $\tau_{V_k, p}(V_k^{\mathrm{sm}}(\Q_p)) < \infty$ by Lemma \ref{lem:Vk_terminal} and Corollary \ref{cor:Tamagawa_measure_converges}.
	
	The lemma follows after splitting the integral depending on the reduction of $\mathbf{t}$ in $(\Z/p^m\Z)^n$ and applying Lemma~\ref{lDensityLift} and Lemma \ref{lem:local_BM_equals_T(g,a)} to evaluate the integrand, cf.~the proof of Lemma \ref{lemma:sigmahat}.
\end{proof}

We do something similar at the infinite place.
\begin{lemma}
	For $\mathbf{t} \in \PPP^{n - 1}(\R)_{\boldsymbol{\varepsilon}}$ we have $e(\inv_{\infty}(\alpha_g(\mathbf{t}))) = T_{\infty}(g, \boldsymbol{\varepsilon})$.
\end{lemma}

\begin{proof}
	Note that the Hilbert symbol at $\infty$ only depends on the signs of the variables. This is then \cite[Lemma 5.7]{CKR}.
\end{proof}

This has the following immediate corollary, noting that $T_{\infty}(g, \boldsymbol{\varepsilon}) \in \{\pm 1\}$.

\begin{corollary}
\label{cor:Brauer_transform_at_infinity}
	We have 
	\[
	T_{\infty}(g, \boldsymbol{\varepsilon}) \hat{\tau}_{\infty}(g; \boldsymbol{\varepsilon}) = \tau_{V_k, \infty}(V_k^{\mathrm{sm}}(\R)_{\boldsymbol{\varepsilon}}) \geq 0.
	\]
\end{corollary}

Define for a prime $p$ (including $p = 2$)
\begin{equation}\label{def:sigmahat}
\hat{\sigma}_p(g) \coloneqq \lim_{m \to \infty}\frac{1}{p^{mn}} \sum_{\mathbf{b}\in\adm_p(m)} 
	 \mathfrak{S}_p(\mathbf{b})T_p(g, \mathbf{b}).
\end{equation}
The limit exists by Lemma~\ref{lem:Brauer_transform_at_p}.

This allows us to prove the main result of this section, which is an upper bound for the Tamagawa measure of the Brauer--Manin subset.

\begin{prop}\label{prop:measure_BM_set}
	We have 
	\[
		\tau_{V_k}\left(V_k^{\mathrm{sm}}(\Adele_\Q)^{\Br} \right) \leq \frac{(2^{16}k)^k}{\lvert\PBM\rvert} \sum_{\boldsymbol{\varepsilon} \in \{\pm 1\}^n}\sum_{g \in \PBM} T_{\infty}(g, \boldsymbol{\varepsilon}) \prod_p \hat{\sigma}_p(g).
	\]
	The Euler product is absolutely convergent and the local factors are given by the limits \eqref{def:sigmahat}.
	\end{prop}

\begin{proof}
	The subspaces $V_k^{\mathrm{sm}}(\Adele_{\Q})_{\boldsymbol{\varepsilon}}$ cover $V_k^{\mathrm{sm}}(\Adele_\Q)$ twice up to a set of measure $0$ as $\boldsymbol{\varepsilon}$ ranges over all of $\{\pm 1\}^n$. It thus suffices to bound each $\tau_{V_k}\left(V_k^{\mathrm{sm}}(\Adele_\Q)_{\boldsymbol{\varepsilon}}^{\Br} \right)$ separately.
	
	By combining Lemma \ref{lem:BM_set_character_orthogonality} and Corollary \ref{cor:Brauer_transform_at_infinity} we have
	\[
	\tau_{V_k}\left(V_k^{\mathrm{sm}}(\Adele_\Q)_{\boldsymbol{\varepsilon}}^{\Br} \right) = \frac{\tau_{V_k, \infty}(V_k^{\mathrm{sm}}(\R)_{\boldsymbol{\varepsilon}})}{\lvert\PBM\rvert} \sum_{g \in \PBM} T_{\infty}(g, \boldsymbol{\varepsilon}) \prod_p (1 - \frac{1}{p})^{k + 1} \hat{\tau}_p(g).
	\]	
	The sum over $g$ is non-negative by Lemma~\ref{lemma:charorth} and Lemma~\ref{lem:Brauer_transform_at_p}, so we can use the bound $\tau_{V_k, \infty}(V_{k}^{\mathrm{sm}}(\R)_{\boldsymbol{\varepsilon}}) \leq \tau_\infty(V_{k}) \leq (2^{16} k)^k$ by Lemma~\ref{lRealBound}.
	The proposition then follows from Lemma \ref{lem:Brauer_transform_at_p}.
\end{proof}

We will also need estimates for the Tamagawa volumes of the following local subsets.

\begin{definition}
\label{def:local_conditions}
	Define the following subsets
	\begin{align*}
		\mathcal{W}_{k, 2} &\coloneqq \{x \in V_k^{\mathrm{sm}}(\Q_2): \pi(x) = \mathbf{t} \in \PPP^{n - 1}(\Q_2) \textup{ with } \mathbf{t} \equiv \mathbf{1} \bmod{8}\}, \\
		\mathcal{W}_{k, p} &\coloneqq \{x \in V_k^{\mathrm{sm}}(\Q_p): \pi(x) = \mathbf{t} \in \PPP^{n - 1}(\Q_p) \textup{ with } \mathbf{t} \bmod{p}\in ((\Z/p\Z)^{\times})^n\}.
	\end{align*}
	The modulo conditions here mean that they hold for some integral representative $\mathbf{t} \in \Z_p^n$.
	
	We will consider the following adelic subset
	\[
	\mathcal{W}_k \coloneqq \{(x_v)_v \in V_k^{\mathrm{sm}}(\Adele_\Q): x_2 \in \mathcal{W}_{k, 2} \textup{ and } x_p \in \mathcal{W}_{k, p} \textup{ for all odd primes } p \leq k^{12k}\}.
	\]
\end{definition}

\begin{lemma}\label{lem:measure_local_conditions}
	We have 
	\[
		\tau_{V_k}(\mathcal{W}_k) \leq 4 \cdot (2^{16}k)^k \frac{\mathfrak{S}_2(\mathbf{1})}{2^{3n}} \prod_{2 < p \leq k^{12 k}} \sum_{\substack{\mathbf{b} \in \adm_p(1) \\ \mathbf{b} \in ((\Z/p\Z)^{\times})^n}} \frac{\mathfrak{S}_p(\mathbf{b})}{p^{n}}
		\prod_{p > k^{12k}} \hat{\sigma}_p(0).
	\]
\end{lemma}

\begin{proof}
	The proof is analogous to that of Proposition \ref{prop:measure_BM_set} with the following modifications. Because we are considering the full volume and not that of the Brauer--Manin set we only have to consider the terms corresponding to $g = 0$. On the other hand, because we are considering local conditions we have to compute the Tamagawa measures of $\mathcal{W}_{k, p}$ for $p \leq k^{12k}$. The same argument as in Lemma \ref{lem:Brauer_transform_at_p} shows that 
	\[
		\lambda_2 \tau_{V_k, 2}(\mathcal{W}_{k, 2}) = 4 \frac{\mathfrak{S}_2(\mathbf{1})}{2^{3n}}, \qquad \lambda_p \tau_{V_k, p}(\mathcal{W}_{k, p}) = \sum_{\substack{\mathbf{b} \in \adm_p(1) \\ \mathbf{b} \in ((\Z/p\Z)^{\times})^n}} \frac{\mathfrak{S}_p(\mathbf{b})}{p^{n}} \textup{ for odd primes } p.
	\]
	The $4 = |(\Z/8\Z)^{\times}|$ appearing in $\tau_{V_k, 2}(\mathcal{W}_{k, 2})$ is there because we are considering projective local conditions.
\end{proof}
\subsection{Mori dream space}\label{ss:Mori}
The variety $V_k$ is not Fano for large $k$. We will show in this section that it is a \emph{Mori dream space} \cite[Definition~1.10]{HuKeel2000}. This is a property which implies that the birational geometry of $V_k$ is very well-behaved. The main application for us will be that it has a small birational modification which is Fano. Since most invariants in Manin's conjecture are invariant under small birational modifications this implies that for the sake of Manin's conjecture $V_k$ is almost as good as Fano.

\begin{prop}
\label{prop:Mori_dream}
	The variety $V_k$ is a Mori dream space.
\end{prop}

Before we can prove this we need to recall the definition of a divisorially log terminal divisor \cite[p.~424]{BCHM2010}. For each $i \in [k]$ write $\Delta_i \coloneqq \{Y_i = 0\} \subset V_k$, and set $\Delta \coloneqq \sum_{i = 1}^k \Delta_i$. Choose a log resolution $f: X \to V_k$ of $(V_k, \Delta)$. We may uniquely write 
\[
f^{*}(K_{V_k} + (1 - \delta)\Delta) = K_X + (1 - \delta)\widetilde{\Delta} + \sum_{j} b_j E_j,
\]
where $\widetilde{\Delta}$ is the proper transform of $\Delta$ and $E_j \subset X_{\overline{\Q}}$ is an exceptional divisor of the resolution (loc.~cit.~writes this as $\tilde{D} = (1 - \delta)\Delta$ and $E = \sum_{j} b_j E_j$). Then $K_{V_k} + (1 - \delta)\Delta$, or equivalently the pair $(V_k, (1 - \delta)\Delta)$, is \emph{divisorially log terminal} if $a_j = 1 - b_j > 0$ for all $j$. In other words, $b_j < 1$.
	
It has been understood since work of Igusa \cite{Igusa1974} that there is a relation between these numbers $a_j$ and the convergence of certain $p$-adic integrals, see \cite{Denef1991} for a nice exposition. In the references it is assumed that the total space is smooth, while it is not difficult to extend the arguments to our more general case we were unable to find a clean form of such a generalisation in the literature so we provide a proof for the sake of completeness.

Let $p$ be a prime such that all divisors $E_j$ are defined over $\Q_p$. Fix a $p$-adic metric $\lVert\cdot\rVert_p$ on every line bundle $\mathcal{O}(E_i)$ and on $L_i$ and let $\lVert\cdot\rVert_{K_X}$ be the induced metric on $K_X$ and $\tau_{X, p}$ the corresponding Tamagawa measure on $X(\Q_p)$.
Let $f_{E_j}$, $f_{\Delta_i}$ and $f_{\widetilde{\Delta}}$ be the global sections defining the divisors $E_j$, $\Delta_i$ and $\widetilde{\Delta}$, respectively. We then have the following criterion.

\begin{lemma}
\label{lem:integrability_implies_dlt}
	Let $p$ be a prime such that all divisors $E_j$ are defined over $\Q_p$ and $E_j(\Q_p) \neq \varnothing$. If $0 < \delta < 1$ and
	\[
	\int_{V_k^{\mathrm{sm}}(\Q_p)} \prod_{i = 1}^k \lVert f_{\Delta_i}\rVert_p(x)^{\delta - 1} \mathrm{d}\tau_{V_k, p}(x) < \infty,
	\]
	then the pair $(V_k, (1 - \delta)\Delta)$ is divisorially log terminal.
\end{lemma}

\begin{proof}
	By our choice of $p$-adic measures we have
	\[
	\int_{V_k^{\mathrm{sm}}(\Q_p)} \prod_i \lVert f_{\Delta_i}\rVert_p(x)^{\delta - 1} \mathrm{d}\tau_{V_k, p}(x) = \int_{X(\Q_p)}\lVert f_{\widetilde{\Delta}}\rVert_{p}(x)^{\delta - 1} \prod_j \lVert f_{E_j}\rVert_p(x)^{-b_j} \mathrm{d}\tau_{X, p}(x).
	\]
	The rightmost integral converges absolutely if and only if $1 - \delta < 1$ and $b_j < 1$ for all $j$ by \cite[Lemma 4.1, Corollary 4.4]{ChambertLoirTschinkel2010}.
\end{proof}

\begin{lemma}\label{lem:dlt}
	The pair $(V_k, (1 - \delta)\Delta)$ is divisorially log terminal for all $0 < \delta < 1$.
\end{lemma}

\begin{proof}
	We may choose a $p$ as in the statement of Lemma \ref{lem:integrability_implies_dlt} by spreading out, the Chebotarev density theorem and the Lang--Weil estimates. By integrating over fibres and lifting the Tamagawa measure on $\PPP^{n-1}$ to the affine cone as in the proofs of Lemmas \ref{lem:integrate_fibers}, \ref{lem:integrate_fibers_affine} we have
	\[
	\int_{V_k^{\mathrm{sm}}(\Q_p)} \prod_i \lVert f_{\Delta_i}\rVert_p(x)^{\delta - 1} \mathrm{d}\tau_{V_k, p}(x) = (1 - \frac{1}{p})^{-1} \int_{\substack{\mathbf{t} \in \Z_p^n \\ \mathbf{t} \textup{ prim.}}} \prod_i \int_{Q_{i, \mathbf{t}}(\Q_p)} \hspace{-0.1cm} \lVert f_{\Delta_i}\rVert_p(x_i)^{\delta - 1} \mathrm{d}\tau_{Q_{i,\mathbf{t}}, p}(x_i) \mathrm{d}\mathbf{t}.
	\]	
	The function $\lVert f_{\Delta_i}\rVert_p$ is equal to $\frac{|Y_i|_p}{\max(|X_i|_p, |Y_i|_p, |Z_i|_p)}$ after potentially enlarging $p$ such that $\lVert\cdot\rVert_p$ is the model metric coming from the integral model of $V_k$ defined by the same equations. The right-hand side converges by Lemma \ref{Weighted_Integral_Total}, so we are done by Lemma \ref{lem:integrability_implies_dlt}.
\end{proof}

\begin{proof}[Proof of Proposition \ref{prop:Mori_dream}]
	By \cite[Corollary 1.3.2]{BCHM2010} it suffices to show that $V_k$ is a normal $\Q$-factorial variety and that there exists a $1 > \delta > 0$ such that the pair $(V_k, (1 - \delta)\Delta)$ is divisorially log terminal and  $-K_{V_k} - (1 - \delta)\Delta$ is ample. The variety is factorial by Corollary \ref{cor:factorial} and the pair $(V_k, (1 - \delta)\Delta)$ is divisorially log terminal by Lemma \ref{lem:dlt}.
	
	We have to show that the divisor class $2k \pi^* H + \delta \sum_{i = 1}^k L_i$ is ample. It is contained in the interior of $\Nef_{V_k}$ for $\delta < \frac{2}{2k - 3}$ by Lemma \ref{lem:cones}. It is then ample by \cite[Theorem 1.4.23]{Lazarsfeld}.
\end{proof}

The main consequence of $V_k$ being a Mori dream space is the following. Recall that a \emph{small} birational modification is a birational map that is an isomorphism in codimension $1$ and that a normal variety $X$ is $\Q$-Gorenstein if a multiple of $K_X$ is Cartier. The following argument is likely well-known to experts.

\begin{prop}\label{prop:V_k_Fano}
	The section ring $R(V_k, -K_{V_k}) \coloneqq \bigoplus_{\ell \in \Z_{\geq 0}}\HH^0(V_k, \omega_{V_k}^{-\ell})$ is a finitely generated graded ring. The corresponding projective variety $V_k^{\mathrm{Fano}} \coloneqq \Proj(R(V_k, -K_{V_k}))$ is a normal $\Q$-Gorenstein Fano variety. 
	
	The natural rational map $f:V_k \dashrightarrow V_k^{\mathrm{Fano}}$ is a small birational modification. If $V_k \xleftarrow{p} X \xrightarrow{q} V_k^{\mathrm{Fano}}$ is a joint resolution of singularities then $-p^*K_{V_k} = -q^*K_{V^{\mathrm{Fano}}_k} + E$ for some effective $\Q$-divisor $E$.
\end{prop}
\begin{proof}
	The ring $R(V_k, -K_{V_k}) \otimes_{\Q} \overline{\Q}$ is finitely generated over $\overline{\Q}$ by \cite[Proposition~1.11(4)]{HuKeel2000} since $V_k$ is a Mori dream space. The ring $R(V_k, -K_{V_k})$ is thus finitely generated by descent \cite[\href{https://stacks.math.columbia.edu/tag/00QP}{Tag 00QP}]{Stacks}.
	
	Lemma \ref{lem:cones} implies that $-K_{V_k}$ is movable. By \cite[Proposition~1.11(5)]{HuKeel2000} there thus exists a small birational contraction $V_k \dashrightarrow Y$ where $Y$ is projective and $\Q$-factorial and $-K_{V_k}$ is semi-ample on $Y$. Consider the map $h: Y \to \PPP^N$ defined by a base-point free power $-mK_{V_k}$ of $-K_{V_k}$. The normalisation of its image is $V_k^{\mathrm{Fano}}$ by \cite[Lemma 1.6]{HuKeel2000}. 
	
	The map $Y \to V_k^{\mathrm{Fano}}$ is a small birational modification. Indeed, it is birational since $-K_{V_k}$ is big by Lemma \ref{lem:cones}. Assume now that it contracted a divisor $D \subset Y$. Then the stable base locus of $-K_{V_k} + \epsilon D$ contains $D$ for all rational $\epsilon > 0$ by the projection formula, which contradicts the fact that $-K_{V_k}$ is in the interior of the movable cone by Lemma \ref{lem:cones}. The natural map $f:V_k \dashrightarrow V_k^{\mathrm{Fano}}$ is then given by composing $V_k \dashrightarrow Y \to V_k^{\mathrm{Fano}}$.
	
	The variety $V_k^{\mathrm{Fano}}$ is $\Q$-Gorenstein and Fano since $\mathcal{O}_{\PPP^{N}}(1)|_{V_k^{\mathrm{Fano}}} = -f_*mK_{V_k}$ by construction, and thus ample. But $-f_*K_{V_k} = -K_{V^{\mathrm{Fano}}_k}$ as Weil divisors since $V_k \dashrightarrow V_k^{\mathrm{Fano}}$ is a small birational modification.
	
	 Let $V_k \xleftarrow{p} X \xrightarrow{q} V_k^{\mathrm{Fano}}$ be a joint resolution of singularities and write 
	 \[
	 -p^*K_{V_k} = -q^*K_{V^{\mathrm{Fano}}_k} + E
	 \] 
	 for some $\Q$-Cartier divisor $E$. The map $f:V_k \dashrightarrow V_{k}^{\mathrm{Fano}}$ is small so the $p$- and $q$-exceptional divisors agree, which implies that $p_* E = 0$. The divisor $-K_{V^{\mathrm{Fano}}_k} = -f_*K_{V_k}$ is ample and thus nef so it is $-K_{V_k}$-non-positive by \cite[Lemma 3.6.3]{BCHM2010}. This means exactly that $E$ is effective, see \cite[Definition 3.6.1]{BCHM2010}.
\end{proof}

\begin{lemma}\label{lem:Fano_model_terminal}
	The variety $V_{k}^{\mathrm{Fano}}$ has terminal singularities.
\end{lemma}

\begin{proof}
	Let $V_k \xleftarrow{p} X \xrightarrow{q} V_k^{\mathrm{Fano}}$ be a joint resolution of singularities. By enlarging $X$ we may assume it to be a log resolution of the pair $(V_k, \Delta)$. Let $E_1, \dots, E_m$ be the exceptional prime divisors of these resolutions. For each $i\in[k]$ let $\widetilde{\Delta}_i$ be the proper transform of $\Delta_i$. Let $b_j(V_k), b_j(V_k^{\mathrm{Fano}}), b_j(\Delta_i) \in \Q$ be the numbers such that
	 \begin{align*}
	 	p^*K_{V_k} &= K_X + \sum_{j = 1}^m b_j(V_k) E_j, \\
	 	q^*K_{V_k^{\mathrm{Fano}}} &= K_X + \sum_{j = 1}^m b_j(V_k^{\mathrm{Fano}}) E_j, \\
	 	p^* \Delta_i &= \widetilde{\Delta}_i+\sum_{j=1}^m b_j(\Delta_i) E_j.
	 \end{align*}
	 Note that $b_j(\Delta_i) \in \Z_{\geq 0}$ as the $\Delta_i$ are effective and Cartier. Our goal is to show that $b_j(V_k^{\mathrm{Fano}}) < 0$ for all $1 \leq j \leq m$ by the definition of terminality.
	
	Recall $Z = \{\textup{at least two } Y_i = 0\} \subset V_k$, which contains the stable base locus $\mathbf{B}(-K_{V_k})$ by Lemma \ref{lem:cones}.
	
	Note that the map $f:V_k \dashrightarrow V_k^{\mathrm{Fano}}$ is the \emph{ample model} of $-K_{V_k}$ in the sense of \cite[Definition 3.6.5]{BCHM2010} by \cite[Lemma 3.6.6(4)]{BCHM2010} and Proposition \ref{prop:V_k_Fano}. This implies by definition that if $D$ is a member of the linear system $|-\ell K_{V_k}|$ for some $\ell > 0$ then we have the inequality of $\Q$-divisors
	\begin{equation}\label{eq:inequality_divisors}
		\frac{p^*D}{\ell} \geq q^*K_{V_k^{\mathrm{Fano}}} - p^*K_{V_k} = \sum_{j = 1}^m (b_j(V_k^{\mathrm{Fano}}) - b_j(V_k)) E_j \geq 0.
	\end{equation}
	 
	Let $1 \leq j \leq m$. We will consider two cases depending on the centre $p(E_j) \subseteq V_k$. If $p(E_j) \not \subseteq \mathbf{B}(-K_{V_k})$ then there exists $\ell \in \Z_{> 0}$ and a divisor $D \in |- \ell K_{V_k}|$ whose support does not contain $p(E_j)$. The support of $p^*D$ thus does not contain $E_j$ so \eqref{eq:inequality_divisors} implies that $b_j(V_k^{\mathrm{Fano}}) = b_j(V_k)$ which is negative by Lemma \ref{lem:Vk_terminal}.
	
	The harder case is when $p(E_j) \subseteq \mathbf{B}(-K_{V_k}) \subseteq Z$. So two of the $Y_i$ vanish on $p(E_j)$, by symmetry we may assume $Y_1$ and $Y_2$. It follows that $b_j(\Delta_1), b_j(\Delta_2) \geq 1$.
	
	Let $\ell = 2k - 3$. The divisor class $2 k \ell \pi^* H + \ell L_1 + L_2$ is globally generated by Lemma \ref{lem:cones} so there exists $D \in |2 k \ell \pi^* H + \ell L_1 + L_2|$ whose support does not contain $p(E_j)$. We will apply \eqref{eq:inequality_divisors} to the following divisor, which is effective since $k \geq 2$:
	\[
	D + (\ell - 1)\Delta_2 + \ell \sum_{i = 3}^k \Delta_i \in |- \ell K_{V_k}| = |2k \ell\pi^* H + \ell \sum_{i = 1}^k L_i|.
	\]
	The resulting inequality is as follows, as the support of $p^*D$ does not contain $E_j$:
	\[
		\ell b_j(V_k^{\mathrm{Fano}}) \leq \ell b_j(V_k) + (\ell - 1) b_j(\Delta_2) + \ell \sum_{i = 3}^k b_j(\Delta_i) \leq \ell b_j(V_k) + \ell \sum_{i = 1}^kb_j(\Delta_i) - \ell - 1.
	\]
	 It follows from Lemma \ref{lem:dlt} and unfolding the definition of divisorially log terminal that for all $0 < \delta < 1$ we have $b_j(V_k) + (1 - \delta)\sum_{i = 1}^k b_j(\Delta_i) < 1$. Taking the limit as $\delta \to 0$ and combining this with the above inequality we deduce that $b_j(V_k^{\mathrm{Fano}}) \leq -1/\ell$. 
\end{proof}

\subsection{\texorpdfstring{Heights on $V_k^{\mathrm{Fano}}$}{Heights on Fano model}}
	Our counting result uses the height $H$ from the metric of Section~\ref{sec:adelic_metrics}, but the Batyrev--Tschinkel prediction will be about points of height with respect to a height on $V_k^{\mathrm{Fano}}$. These are not the same, which is why we will construct a sequence of heights $H^{\mathrm{Fano}}_m$ on $V_k^{\mathrm{Fano}}$ such that $H^{\mathrm{Fano}}_m \leq H$ and which converges to $H$ as $m \to \infty$. 
	
	Let $U \subseteq V_k$ be a smooth open subvariety such that $f: V_k \dashrightarrow V_k^{\mathrm{Fano}}$ restricts to an isomorphism $U \cong f(U)$ and such that the complement of $U \subset V_k$ has codimension at least $2$. This exists as $f$ is a small birational modification of normal varieties. We will identify $U$ with $f(U)$.
	
	Let $r \in \Z_{> 0}$ be such that $rK_{V_k^{\mathrm{Fano}}}$ is Cartier. A local metric $\lVert \cdot\rVert^{\mathrm{Fano}}_v$ on $rK_{V_k^{\mathrm{Fano}}}$ can be compared to the metric $\lVert \cdot\rVert_v^r$ of Section~\ref{sec:adelic_metrics} by comparing their restrictions on $U(\Q_v)$. This comparison is independent of the chosen $U$ as $U(\Q_v)$ is dense in $V_k^{\mathrm{sm}}(\Q_v)$.
	
	\begin{lemma}
	\label{lem:existence_approximating_metrics}
		There exists a sequence of adelic metrics $(\lVert \cdot\rVert^{\mathrm{Fano}}_{m, v})_v$ on $r K_{V_k^{\mathrm{Fano}}}$ such that for all $v$ and $m$ we have $\lVert \cdot\rVert^{\mathrm{Fano}}_{m, v} \leq \lVert \cdot\rVert^{\mathrm{Fano}}_{m + 1, v} \leq \lVert \cdot\rVert_v^r$ on $U(\Q_v)$ and as $m \to \infty$ we have $\lVert \cdot\rVert^{\mathrm{Fano}}_{m, v} \to \lVert \cdot\rVert_v^r$ on $U(\Q_v)$ uniformly on compact subsets of $U(\Q_v)$.
	\end{lemma}
	\begin{proof}
		Let $V_k \xleftarrow{p} X \xrightarrow{q} V_k^{\mathrm{Fano}}$ be as in the statement of Proposition \ref{prop:V_k_Fano} and assume moreover that $p$ and $q$ are isomorphisms over $U$. Let $E \coloneqq q^* rK_{V_k^{\mathrm{Fano}}} - p^* r K_{V_k}$ which is an effective Cartier divisor by this proposition.
		
		Let $(\lVert \cdot\rVert'_v)_{v}$ be some adelic metric on $rK_{V_k^{\mathrm{Fano}}}$ and let $\lVert \cdot\rVert_{E, v}$ be the local metric on $E$ such that $p^*\lVert \cdot\rVert_v^r \otimes \lVert \cdot\rVert_{E, v} = q^*\lVert \cdot\rVert'_v $.
		
		Let $f_{E}$ be the global section defining the divisor $E$. This is invertible on $p^{-1}(U) \cong U$ and thus defines an isomorphism $\omega_{V_k}^r|_U \cong \omega_{V_k^{\mathrm{Fano}}}^r|_U$. By multiplying $f_{E}$ by an invertible constant we may assume that this isomorphism is equal to the canonical isomorphism $\omega_{V_k}^r|_U \cong \omega_U^r \cong \omega_{V_k^{\mathrm{Fano}}}^r|_U$ as $\Q[U] = \Q$ by algebraic Hartogs' lemma.
		
		On $U(\Q_v)$ we thus have $ \lVert \cdot\rVert'_v = \lVert f_{E}\rVert_{E, v} \cdot \lVert \cdot\rVert_v^r$. Moreover, $\lVert f_E\rVert_{E, v}$ is a bounded continuous function as it extends to a continuous function on the compact space $X(\Q_v)$. Since these local metrics come from an adelic metric there exists by spreading out a finite set of primes $S$ and a $\Z_S$-model $\mathcal{U}$ of $U$ such that $\lVert f_E\rVert_{E, p}(u) \in p^{\Z_{\leq 0}}$ for all primes $p \not \in S$ and $u \in U(\Q_p)$, with the equality $\lVert f_E\rVert_{E, p}(u) = 1$ if $u \in \mathcal{U}(\Z_p)$.
		
		For each place $v$ fix a sequence $(K_{v,m})_{m}$ of compact subsets $K_{v, m} \subseteq U(\Q_v)$ such that $K_{v,m} \subseteq K_{v, m + 1}$ is a compact neighbourhood and $\bigcup_m K_{v, m} = U(\Q_v)$. Assume moreover that if $p \not \in S$ is prime then $K_{p, m} = \mathcal{U}(\Z_p)$ for all $m$ with $p > 2^m$.
		
		For each place $v$ and $m$ fix a continuous function $h_{v, m}: V_k^{\mathrm{Fano}}(\Q_v) \to \R_{> 0}$ with the following properties. We have $h_{v, m}(u) \leq \lVert f_E\rVert_{E, v}(u)^{-1}$ for all $u \in U(\Q_v)$ with equality if $u \in K_{v, m}$, we have $h_{v, m} \leq h_{v, m + 1}$ and for all primes $p \not \in S$ and $p > 2^m$ we have $h_{p, m}(y) = 1$ for all $y \in V_k^{\mathrm{Fano}}(\Q_p)$. This exists as $V_k^{\mathrm{Fano}}(\Q_v)$ is a compact Hausdorff space and $K_{v, m} \subseteq V_k^{\mathrm{Fano}}(\Q_v)$ is compact.
		
		We then define $\lVert \cdot\rVert^{\mathrm{Fano}}_{m, v} = h_{v, m} \lVert \cdot\rVert'_v$. The inequalities satisfied by $h_{v, m}$ imply the inequalities in the statement. For any compact $K \subseteq U(\Q_v)$ we have $K \subseteq K_{v, m_K}$ for some $m_K$ by construction so $\lVert \cdot\rVert^{\mathrm{Fano}}_{m, v} = \lVert \cdot\rVert_v^{r}$ holds on $K$ for all $m \geq m_K$.
	\end{proof}
	
	Let $\tau^{\mathrm{Fano}}_m$ be the Tamagawa measure on $U(\Adele_{\Q})$ coming from the metric $(\lVert \cdot\rVert^{\mathrm{Fano}}_{m, v})^{1/r}_v$. The inequalities of Lemma \ref{lem:existence_approximating_metrics} imply that $\tau_{V_k} \leq \tau^{\mathrm{Fano}}_{m + 1} \leq \tau^{\mathrm{Fano}}_m$ and uniform convergence on compact subsets implies $\tau^{\mathrm{Fano}}_m \xrightarrow{m \to \infty} \tau_{V_k}$.
	
	\begin{lemma}
	\label{lem:convergence_of_measures}
		For any Borel subset $\Omega \subseteq U(\Adele_{\Q})$ we have $\lim_{m \to \infty}\tau^{\mathrm{Fano}}_m(\Omega) = \tau_{V_k}(\Omega)$.
	\end{lemma}
	
	\begin{proof}
		We may exchange the limit $\tau^{\mathrm{Fano}}_m \xrightarrow{m \to \infty} \tau_{V_k}$ and the integral by dominated convergence as $ \tau^{\mathrm{Fano}}_1(\Omega) < \infty$ by Lemma \ref{lem:Fano_model_terminal} and Corollary \ref{cor:Tamagawa_measure_converges}.
	\end{proof}
	
\subsection{The Batyrev--Tschinkel prediction}
\label{ss:BTconj}
We will now discuss the prediction of Batyrev--Tschinkel in \cite[p.~323: Expected asymptotic formula]{BT98}. We will consider the smooth variety $U$ and the line bundle $\mathcal{L} \coloneqq -K_{V_k}|_U$ to be consistent with their notation. We equip $-K_{V_k^{\mathrm{Fano}}}$ with the adelic metric $(\lVert \cdot\rVert^{\mathrm{Fano}}_{m, v})_v^{-1/r}$ for some $m$, denote the induced height on $V_k^{\mathrm{Fano}}$ by $H_m^{\mathrm{Fano}}$, and restrict it to $\mathcal{L}$ via $f|_U: U \xrightarrow{\sim} f(U)$.

The prediction in \cite{BT98} is stated in terms of language that is no longer in vogue. For example, Batyrev--Tschinkel require the variety to be \emph{strongly saturated} due to certain counterexamples, such as the quadric and cubic surface bundles discussed in \cite[Examples 3.5.2 and 3.5.3]{BT98}. It is now understood that the better way to deal with these and similar counterexamples is to remove a thin subset, see \cite[\S 8]{Peyre2003}. Moreover, they require the variety $U$ to be $\mathcal{L}$-primitive with respect to the chosen line bundle $\mathcal{L}$. The most important examples of such varieties are the smooth loci of $\Q$-Gorenstein canonical Fano varieties, see \cite[Example 2.3.7]{BT98}. The variety $V_k^{\mathrm{Fano}}$ is terminal and Fano by Lemma \ref{lem:Fano_model_terminal} and hence $U$ is $\mathcal{L}$-\emph{primitive} and $a_{\mathcal{L}}(V_k^{\mathrm{Fano}}) = 1$.

There are three rational numbers appearing in the prediction of \cite{BT98}: $\beta_{\mathcal{L}}(V_k^{\mathrm{Fano}})$, $\delta_{\mathcal{L}}(V_k^{\mathrm{Fano}})$ and $\gamma_{\mathcal{L}}(V_k^{\mathrm{Fano}})$, defined in \cite[Definitions 3.4.3, 2.3.16]{BT98}. These are defined in terms of the complement of the non-crepant divisors for a resolution of singularities $X \to V_k^{\mathrm{Fano}}$. Henceforth we assume that this complement is $U$ as $V_k^{\mathrm{Fano}}$ is terminal. Note that all of these invariants only depend on $U$ up to codimension $2$ so we can and will compute these invariants on $V_k^{\mathrm{sm}}$.
\begin{itemize}
\item The number $\beta_{\mathcal{L}}(V_k^{\mathrm{Fano}})$ is the Picard rank of $U$, which is equal to $k+ 1 $ by Lemma \ref{lem:Picard_group_computation}.
\item The cohomological constant $\delta_{\mathcal{L}}(V_k^{\mathrm{Fano}})$ is the size of $\HH^1(\Q, \Pic U_{\overline{\Q}}) = \HH^1(\Q, \Pic \, (V_k^{\mathrm{sm}})_{\overline{\Q}})$. We thus have $\delta_{\mathcal{L}}(V_k^{\mathrm{Fano}}) = 1$ by Lemma \ref{lem:algebraic_Brauer}.
\item The number $\gamma_{\mathcal{L}}(V_k^{\mathrm{Fano}})$ is \emph{Peyre's effective cone constant}. It is at most $ \frac{1}{2k}$ by Lemma \ref{lem:effective_cone_constant}.
\end{itemize}

The last invariant appearing in the expected asymptotic formula is $\tau_{\mathcal{L}}(V_k^{\mathrm{Fano}})$. They define it as the Tamagawa measure of the topological closure $\overline{X(\Q)} \subseteq X(\Adele_{\Q})$ of $X(\Q)$. By Proposition \ref{prop:Tamagawa_measure_terminal} this Tamagawa measure is equal to $\tau_m^{\mathrm{Fano}}(\overline{X(\Q)} \cap U(\Adele_{\Q}))$.
Let $\overline{U(\Q)} \subseteq U(\Adele_{\Q})$ be the topological closure of $U(\Q)$. We have $\overline{U(\Q)} \subseteq \overline{X(\Q)} \cap U(\Adele_{\Q})$ but this inclusion can be strict in general. If the inclusion is strict, then Batyrev--Tschinkel's definition is unreasonable and a better definition is $\tau_{\mathcal{L}}(V_k^{\mathrm{Fano}}) = \tau_m^{\mathrm{Fano}}(\overline{U(\Q)})$.

\begin{remark}\label{rem:tau}
	One can construct a terminal toric variety $V$ such that $\overline{V^{\mathrm{sm}}(\Q)} \neq \overline{X(\Q)} \cap V^{\mathrm{sm}}(\Adele_{\Q})$. A careful reading of Batyrev--Tschinkel's work on Manin's conjecture for toric varieties \cite{BatyrevTschinkel1996HZF} together with \cite[Theorem 4.5]{Loughran2018Varieties} to identify automorphic characters on the torus with Brauer elements shows that their expected asymptotic formula holds for this $V$ only if one uses $\tau_{\mathcal{L}}(V) = \tau_{V^{\mathrm{sm}}}(\overline{V^{\mathrm{sm}}(\Q)})$ as the definition of the Tamagawa number.
\end{remark}

We will call \cite[p.~323: Expected asymptotic formula]{BT98}, taking into account a possible thin subset, with the original definition $\tau_{\mathcal{L}}(V_k^{\mathrm{Fano}}) = \tau_m^{\mathrm{Fano}}(\overline{X(\Q)} \cap U(\Adele_{\Q}))$ the \emph{original Batyrev--Tschinkel asymptotic formula}. The same formula where we instead use this new definition $\tau_{\mathcal{L}}(V_k^{\mathrm{Fano}}) = \tau_m^{\mathrm{Fano}}(\overline{U(\Q)})$ will be called the \emph{modified Batyrev--Tschinkel asymptotic formula}.

We note the following.
\begin{lemma}
	The variety $X$ satisfies weak approximation, i.e.~$\overline{X(\Q)} = X(\Adele_{\Q})$.
\end{lemma}
\begin{proof}
	Weak approximation is a birational invariant of smooth proper varieties, so this follows from Lemma \ref{lem:V_k_rational}.
\end{proof}

The Brauer--Manin set contains the closure so $\tau_m^{\mathrm{Fano}}(\overline{U(\Q)}) \leq \tau_m^{\mathrm{Fano}}(U(\Adele_{\Q})^{\Br})$. Recall the local conditions $\mathcal{W}_{k, p}$ and the adelic subset $\mathcal{W}_k \subset V_k^{\mathrm{sm}}(\Adele_{\Q})$ from Definition~\ref{def:local_conditions}.

\begin{lemma}
\label{lemma:positive}
	We have
	\[
	\tau_{V_k}(V_k^{\mathrm{sm}}(\Adele_{\Q})^{\Br}) > 0, \qquad \tau_{V_k}(\mathcal{W}_k) > 0.
	\]
\end{lemma}

\begin{proof}
	Observe that $V_k^{\mathrm{sm}}(\Adele_{\Q})^{\Br}$ and $\mathcal{W}_k$ are open subsets of $V_k^{\mathrm{sm}}(\Adele_{\Q})$ (using Proposition~\ref{prop:Brauer_group=RALt} for the former) containing $\varpi(\mathbf{1}, (3, 5, 4), \dots, (3, 5, 4)) \in V_k^{\mathrm{sm}}(\Q)$.
\end{proof}

Our counting theorem counts rational points on $V_k$ with respect to $H$ rather than on $V_k^{\mathrm{Fano}}$ with respect to $H_m^{\mathrm{Fano}}$, but the two can be compared to deduce the following upper bounds on $V_k(\Q)$.
\begin{theorem}\label{theorem:BTconstant}\hfill
	\begin{enumerate}
		\item The modified Batyrev--Tschinkel asymptotic formula predicts that there exists a thin subset $\Omega \subset V_k(\Q)$ such that as $B \to \infty$
		\[
		\#\{P \in V_k(\Q) \setminus \Omega: H(P) \leq B\} \leq (2c + o(1)) B (\log B)^{k},
		\]
		where
		\[
		c = \frac{1}{2k \cdot k!}\tau_{V_k}(V_k^{\mathrm{sm}}(\Adele_{\Q})^{\Br})>0.
		\]
		\item The original Batyrev--Tschinkel asymptotic formula for all anticanonical heights predicts that there exists a thin subset $\Omega \subset V_k(\Q)$ such that as $B \to \infty$
		\begin{multline*}
		\#\{P \in V_k(\Q) \setminus \Omega: H(P) \leq B, \ P \in \mathcal{W}_{k, p} \textup{ for all primes } p \leq k^{12k}\} \\ \leq (2c' + o(1)) B (\log B)^{k},
		\end{multline*}
		where
		\[
		c' = \frac{1}{2k \cdot k!}\tau_{V_k}(\mathcal{W}_k)>0.
		\]
	\end{enumerate}
\end{theorem}

\begin{proof}
	The complement of $U \subseteq V_k^{\mathrm{sm}}$ has codimension at least $2$ so their Brauer groups agree by purity \cite[Theorem 3.7.6]{ColliotTheleneSkorobogatov2021}. By Lemma \ref{lem:measure_change_codimension_2} we thus have $\tau_{V_k}(V_k^{\mathrm{sm}}(\Adele_{\Q})^{\Br}) = \tau_U(U(\Adele_{\Q})^{\Br})$ and $\tau_{V_k}(\mathcal{W}_k) = \tau_U(\mathcal{W}_k \cap U(\Adele_{\Q}))$. By Lemma~\ref{lemma:positive}, we have $c, c' > 0$, so by Lemma \ref{lem:convergence_of_measures} we may choose $m$ such that
	\begin{equation}\label{eq:tau_m_bound}
	\frac{1}{2k \cdot k!}\tau_m^{\mathrm{Fano}}(U(\Adele_{\Q})^{\Br}) \leq 2c \quad \textup{and} \quad \frac{1}{2k \cdot k!}\tau_m^{\mathrm{Fano}}(\mathcal{W}_k \cap U(\Adele_{\Q})) \leq 2c'.
	\end{equation}
	
		Substituting $a_{\mathcal{L}}(V_k^{\mathrm{Fano}}) = 1$, $\beta_{\mathcal{L}}(V_k^{\mathrm{Fano}}) = k + 1$, $\delta_{\mathcal{L}}(V_k^{\mathrm{Fano}}) = 1$ and $\gamma_{\mathcal{L}}(V_k^{\mathrm{Fano}}) \leq \frac{1}{2k}$ from above, together with $\tau_m^{\mathrm{Fano}}(\overline{U(\Q)}) \leq \tau_m^{\mathrm{Fano}}(U(\Adele_{\Q})^{\Br})$ and \eqref{eq:tau_m_bound}, the modified Batyrev--Tschinkel asymptotic formula for $H_m^{\mathrm{Fano}}$ predicts that there exists a thin subset $\Omega \subset V_k^{\mathrm{Fano}}(\Q)$ such that as $B \to \infty$
	\[
	\#\{P \in V_k^{\mathrm{Fano}}(\Q) \setminus \Omega: H_m^{\mathrm{Fano}}(P) \leq B\} \leq (2c + o(1)) B (\log B)^{k}.
	\]
	The original Batyrev--Tschinkel asymptotic formula for all anticanonical heights also makes predictions for counting with respect to local conditions at finitely many primes by the abstract equidistribution result \cite[Proposition 2.10]{ChambertLoirTschinkel2010} (in the setting of smooth Fano varieties this goes back to \cite[Proposition 3.3]{Peyre}): since $\mathcal{W}_{k, p}$ is clopen in $V_k^{\mathrm{sm}}(\Q_p)$, the boundary of $f(\mathcal{W}_{k, p} \cap U(\Q_p))$ in $V_k^{\mathrm{Fano}}(\Q_p)$ lies in $(V_k^{\mathrm{Fano}} \setminus f(U))(\Q_p)$, which has $\tau_m^{\mathrm{Fano}}$-measure zero by Lemma~\ref{lem:Fano_model_terminal} and Proposition~\ref{prop:Tamagawa_measure_terminal}. Using $\overline{X(\Q)} = X(\Adele_{\Q})$ together with \eqref{eq:tau_m_bound}, it predicts that there exists a thin subset $\Omega \subset V_k^{\mathrm{Fano}}(\Q)$ such that as $B \to \infty$
\[
	\#\left\{P \in V_k^{\mathrm{Fano}}(\Q) \setminus \Omega: \begin{array}{l}
	H_m^{\mathrm{Fano}}(P) \leq B,\\ P \in f(\mathcal{W}_{k, p} \cap U(\Q_p)) \textup{ for all } p \leq k^{12k}\end{array}\right\} 
	\leq (2c' + o(1)) B (\log B)^{k}.
\]
	
	In each case, $\Omega' \coloneqq (f|_U)^{-1}(\Omega \cap f(U)(\Q)) \cup (V_k \setminus U)(\Q)$ is a thin subset of $V_k(\Q)$ with $V_k(\Q) \setminus \Omega' \subseteq U(\Q)$. On $U(\Q_v)$ we have $\lVert \cdot\rVert^{\mathrm{Fano}}_{m, v} \leq \lVert \cdot\rVert_v^r$ by Lemma~\ref{lem:existence_approximating_metrics}, so $H_m^{\mathrm{Fano}}(f(P)) \leq H(P)$ for all $P \in U(\Q)$, since the two inverses from $K_{V_k}$ to $-K_{V_k}$ and in the height definition cancel. Hence $f|_U$ embeds $\{P \in V_k(\Q) \setminus \Omega': H(P) \leq B\}$ into $\{Q \in V_k^{\mathrm{Fano}}(\Q) \setminus \Omega: H_m^{\mathrm{Fano}}(Q) \leq B\}$, and for $P \in U(\Q)$ we have $P \in \mathcal{W}_{k, p}$ if and only if $f(P) \in f(\mathcal{W}_{k, p} \cap U(\Q_p))$, yielding both claims.
\end{proof}

\section{Denouement}
\label{s:thinsets}
Finally, we are ready to conclude the proofs of our main theorems. We first demonstrate that the number of points coming from any given thin set has a negligible effect on the point count. Then, we show that the efforts of Section \ref{sCounting} have produced a lower bound strictly larger than the main term predicted by Conjecture \ref{conj:BT}. Together, and with Theorems \ref{tLowerBound} and \ref{theorem:BTconstant}, these results establish our contradiction of Conjecture \ref{conj:BT}.

\subsection{Thin sets}
Having produced too many rational points on our variety, it just remains to rule out the possibility that this departure from the predicted behaviour can be explained by a thin set. The argument laid out in Section \ref{sCounting} does not count all points since it only aims for a lower bound. Let $\Omega \subset V_k(\Q)$ be a thin set of rational points on the variety defined by the system \eqref{eQuadForms}. For a prime $p$, let $\Omega(\Z/p\Z)\subseteq V_k(\Z/p\Z)$ be the reduction of $\Omega$ modulo $p$.

Recall the quotient morphism $\varpi: V_k^{\mathrm{prim}} \to V_k(\Q)$ from Remark~\ref{rem:Cox}. We will bound the number of points in $\mathcal{U}_{\mathcal{E}, m}(B; V_k) \cap \varpi^{-1}(\Omega)$ and show that it is an order of magnitude smaller than the main term produced in Theorem~\ref{tLowerBound}.

\begin{theorem}
\label{thm:thinsets}
	Let $k\geq 2$ and $\Omega \subset V_k(\Q)$ be a thin subset. Let $\mathcal{E} \geq 2$ be a real number and let $m \geq 1$. Then, as $B \rightarrow \infty$, we have
	\[
	\#\left(\mathcal{U}_{\mathcal{E}, m}(B; V_k) \cap \varpi^{-1}(\Omega)\right) = o(B (\log B)^k).
	\]
\end{theorem}

\begin{proof}
Without loss of generality, we may assume that $\Omega$ is a single thin set of type I or II.
We will proceed to produce an upper bound on the number of points in $\varpi^{-1}(\Omega)$ by enforcing that $\varpi(\mathbf{t}, \mathbf{x}) \bmod p \in \Omega(\Z/p\Z)$. By \cite[Lemma 3.8]{BLsieving}, we know that there exists a positive density set of primes $\mathcal{P}$ and a constant $C \in [0,1)$ such that 
\begin{equation}
\label{eq:Serre}
\#\Omega(\Z/p\Z) \leq C p^{5k - 1} + O\left( p^{5k-\frac{3}{2}}\right)
\end{equation}
for all $p\in\mathcal{P}$.

For every $\mathcal{E}' \geq \mathcal{E}$ and every $m' \geq m$ sufficiently large in terms of $k$, we see that $\mathcal{U}_{\mathcal{E}, m}(B; V_k)\subseteq\mathcal{U}_{\mathcal{E}', m'}(B; V_k)$ so 
\[
\#( \mathcal{U}_{\mathcal{E}, m}(B; V_k)\cap \varpi^{-1}(\Omega))\leq\#( \mathcal{U}_{\mathcal{E}', m'}(B; V_k)\cap\varpi^{-1}(\Omega)).
\]
Apply Lemma~\ref{lemma:stepU} to all $\mathbf{b} \in \adm(\mathcal{E}', m')$ and all $\boldsymbol{\upsilon} = (\boldsymbol{\upsilon}_1, \dots, \boldsymbol{\upsilon}_k)$ with $\boldsymbol{\upsilon}_i \in Q_{i, \mathbf{b}}(\Z/C(\mathcal{E}', 1)\Z)$ such that $\varpi(\mathbf{b} \bmod p, \boldsymbol{\upsilon} \bmod p) \in \Omega(\Z/p\Z)$ for all $p\leq \mathcal{E}'$. Then using $|T_v(g, \mathbf{b})|\leq 1$ and $\lvert \PBM\rvert=2^{\frac{1}{2}k(k+1)}$, we find that
\[
\#( \mathcal{U}_{\mathcal{E}', m'}(B; V_k)\cap\varpi^{-1}(\Omega))\ll_k
B (\log B)^k 
\prod_{p \leq \mathcal{E}'} \frac{1}{p^{m' n}} \sum_{\substack{\mathbf{b},\boldsymbol{\upsilon}\\ \mathbf{b} \in\adm_p(m')\\ \varpi(\mathbf{b} \bmod p, \boldsymbol{\upsilon}) \in\Omega(\Z/p\Z)}} \mathfrak{S}_p(\mathbf{b},\boldsymbol{\upsilon}).
\]
For primes $p \in \mathcal{P}$ with $p > k^{12k}$ (so $p$ is odd), split the sum according to whether $\mathbf b\in((\Z/p^{m'}\Z)^\times)^{n}$. For $\mathbf b\in((\Z/p^{m'}\Z)^\times)^{n}$, Lemma~\ref{lem:densitywithcongs} gives
\[
\mathfrak S_p(\mathbf b,\boldsymbol\upsilon)
=p^{-k}(1-p^{-1})^k.
\]
Each point in $\Omega(\Z/p\Z)$ has at most $p - 1$ preimages $(\mathbf{b} \bmod p, \boldsymbol{\upsilon})$ under $\varpi$, and each residue class modulo $p$ has $p^{(m' - 1)n}$ lifts to $\adm_p(m')$. The contribution from $\mathbf b\notin((\Z/p^{m'}\Z)^\times)^{n}$ is $O_k(p^{-1})$ by Lemma~\ref{lem:taupestimate} and Lemma~\ref{lDensityLift}. Hence
\[
\frac{1}{p^{m' n}} \sum_{\substack{\mathbf{b},\boldsymbol{\upsilon}\\ \mathbf{b} \in\adm_p(m')\\ \varpi(\mathbf{b} \bmod p,\boldsymbol{\upsilon})\in\Omega(\Z/p\Z)}}\mathfrak{S}_p(\mathbf{b},\boldsymbol{\upsilon})
\leq p^{-5k}(p - 1) \#\Omega(\Z/p\Z)+O_k(p^{-1}),
\]
which is at most $C + O_{k, \Omega}(p^{-1/2})$ for all odd $p\in\mathcal P$ by~\eqref{eq:Serre}. For primes $p\notin\mathcal{P}$ with $p > k^{12k}$, we drop the thin set condition and apply Lemma~\ref{lemma:sigmahat} with $g=0$,
\[
\frac{1}{p^{m' n}}\sum_{\substack{\mathbf{b},\boldsymbol{\upsilon}\\ \mathbf{b} \in\adm_p(m')}} \mathfrak{S}_p(\mathbf{b},\boldsymbol{\upsilon})
\leq\frac{1}{p^{m' n}}\sum_{\mathbf b\in\adm_p(m')}
\mathfrak S_p(\mathbf b)
=1+O_k(p^{-2}).
\]
For primes $p \leq k^{12k}$, we bound
\[
\frac{1}{p^{m' n}} \sum_{\substack{\mathbf{b},\boldsymbol{\upsilon}\\ \mathbf{b} \in\adm_p(m')}} \mathfrak{S}_p(\mathbf{b},\boldsymbol{\upsilon}) \leq \hat{\sigma}_p(0),
\]
which does not depend on $\mathcal{E}'$ or $m'$. Indeed, this bound holds since the left-hand side equals $\int_{A_{m'}}\prod_{i}\sigma_p(Q_{i,\mathbf t})\mathrm{d}\mathbf{t}$ with $A_{m'}\coloneqq\{\mathbf t\in\Z_p^n:\mathbf t\bmod p^{m'}\in\adm_p(m')\}$. Since $A_{m'}\subseteq A_{m'+1}$, the integral is non-decreasing in $m'$, so it is at most its limit $\hat\sigma_p(0)$.

Therefore for every fixed $\mathcal E'\geq\mathcal E$,
\[
\#(\mathcal U_{\mathcal E, m}(B;V_k)\cap\varpi^{-1}(\Omega))
\ll_{k,\Omega}B(\log B)^k
\prod_{\substack{k^{12k} < p \leq \mathcal E'\\p\in\mathcal P}}
\left(C+O_{k,\Omega}(p^{-1/2})\right),
\]
which is sufficient since the Euler product is $o(1)$ as $\mathcal{E}'$ grows to infinity.
\end{proof}

\subsection{Comparing main terms}
The aim of this subsection is to establish the following result, which visibly contradicts the Batyrev--Tschinkel predictions from Theorem~\ref{theorem:BTconstant}.

\begin{theorem}
\label{thm:maintermcomparison}
For every $M > 0$, for sufficiently large $k$ in terms of $M$, and for every thin subset $Z(\Q) \subseteq \Omega \subset V_k(\Q)$, we have for $B$ sufficiently large in terms of $M$, $k$ and $\Omega$
\[
\# \{P \in V_k(\Q)\setminus\Omega : H(P) \leq B\} \geq \frac{M}{k \cdot k!} \tau_{V_k}(V_k^{\mathrm{sm}}(\Adele_{\Q})^{\Br}) \cdot B (\log B)^k,
\]
and
\[
\# \{P \in V_k(\Q)\setminus\Omega : H(P) \leq B, \ P \in \mathcal{W}_{k, p} \textup{ for all primes } p \leq k^{12k}\} \geq \frac{M \tau_{V_k}(\mathcal{W}_k) B (\log B)^k}{k \cdot k!},
\]
where $\tau_{V_k}(V_k^{\mathrm{sm}}(\Adele_{\Q})^{\Br}) > 0$ and $\tau_{V_k}(\mathcal{W}_k) > 0$.
\end{theorem}

\begin{proof}
The positivity of $\tau_{V_k}(V_k^{\mathrm{sm}}(\Adele_{\Q})^{\Br}) $ and $\tau_{V_k}(\mathcal{W}_k) $ follows from Lemma~\ref{lemma:positive}. Write $N(B)$ and $N'(B)$ for the two point counts on the left-hand sides above. For real $\mathcal{E} \geq 2$ and integer $m \geq 3$, note that 
\begin{align*}
\mathcal{U}'_{\mathcal{E},m}(B; V_k) &\subseteq \mathcal{U}_{\mathcal{E},m}(B; V_k), \\
\varpi(\mathcal{U}_{\mathcal{E},m}(B; V_k)) &\subseteq \{P \in V_k(\Q) : H(P) \leq B\}. 
\end{align*}
Also recall that $t_i > 0$ for all $i \in [n]$ implies $\varpi(\mathbf{t}, \mathbf{x}) \in V_k^{\mathrm{sm}}(\Q)$. This yields
\[
\varpi(\mathcal{U}'_{\mathcal{E},m}(B; V_k)) \subseteq \{P \in V_k(\Q) : H(P) \leq B, \ P \in \mathcal{W}_{k, p} \textup{ for all primes } p \leq k^{12k}\}.
\]
Since $\varpi: V_k^{\mathrm{prim}} \to V_k(\Q)$ is $2^{k+1}$-to-$1$, Theorem~\ref{tLowerBound} and Theorem~\ref{thm:thinsets} give, for fixed sufficiently large $\mathcal{E}, m$ and sufficiently large $B$,
\[
N(B) \geq \frac{B (\log B)^k S_{\mathcal{E}, m}(V_k)}{2^{k + 2} k^{c_3k}} \quad \textup{and} \quad N'(B) \geq \frac{B (\log B)^k S'_{\mathcal{E}, m}(V_k)}{2^{k + 2} k^{c_3k}},
\]
with $S_{\mathcal{E}, m}(V_k)$ and $S'_{\mathcal{E}, m}(V_k)$ as in \eqref{eq:defSem} and \eqref{eq:defSemprime}. Set
\[
P_k \coloneqq \lim_{\mathcal{E}\rightarrow\infty}\lim_{m\rightarrow\infty} S_{\mathcal E,m}(V_k), \qquad P'_k \coloneqq \lim_{\mathcal{E}\rightarrow\infty}\lim_{m\rightarrow\infty} S'_{\mathcal E,m}(V_k).
\]
Defining $c_k \coloneqq \frac{\tau_{V_k}(V_k^{\mathrm{sm}}(\Adele_{\Q})^{\Br})}{k \cdot k!}$ and $c_k' \coloneqq \frac{\tau_{V_k}(\mathcal{W}_k)}{k \cdot k!}$, Proposition~\ref{prop:measure_BM_set} yields
	\[
	c_k\leq \frac{(2^{16}k)^k}{k \cdot k! \lvert\PBM\rvert} P_k
	\] 
	by the limit definition of $\hat{\sigma}_p(g)$ from \eqref{def:sigmahat}, while Lemma~\ref{lem:measure_local_conditions} and the $m \to \infty$ limit of the second part of Lemma~\ref{lemma:strictlypositive} (via Lemma~\ref{lDensityLift} and Lemma~\ref{lemma:sigmahat}) yield
	\[
	c_k' \ll \frac{(2^{16}k)^k}{k \cdot k! \lvert\PBM\rvert} P'_k.
	\]
By Lemma~\ref{lemma:strictlypositive}, for fixed $k$ there is $\delta_k>0$, independent of sufficiently large $\mathcal E,m$, such that $S_{\mathcal E,m}(V_k) \geq S'_{\mathcal E,m}(V_k) \geq \delta_k$. The Euler product $\prod_p \hat{\sigma}_p(g)$ is absolutely convergent by Proposition~\ref{prop:measure_BM_set}. For large primes, we can apply Lemma~\ref{lemma:sigmahat}. This allows us to choose sufficiently large $\mathcal E,m$ such that
\[
\left|S_{\mathcal E,m}(V_k)-P_k\right|\leq\frac{\delta_k}{2} \quad \textup{and} \quad \left|S'_{\mathcal E,m}(V_k)-P'_k\right|\leq\frac{\delta_k}{2}.
\]
Therefore $P_k\leq \frac{3}{2}S_{\mathcal E,m}(V_k)$ and $P'_k \leq \frac{3}{2}S'_{\mathcal E,m}(V_k)$, and
\[
\frac{c_k B (\log B)^k}{N(B)} 
\ll \frac{2^{k + 2}(2^{16}k)^k}{k \cdot k! \lvert\PBM\rvert} k^{c_3k}
\ll\frac{k^{O(k)}}{2^{\frac{1}{2}k(k+1)}}=o(1).
\]
Since we can estimate $\frac{c_k' B (\log B)^k}{N'(B)}$ in a similar fashion, this ends the proof.
\end{proof}


\begin{thebibliography}{30}
\bibitem{BatyrevTschinkel1996HZF}
V.~V.~Batyrev and Y.~Tschinkel,
Height zeta functions of toric varieties,
\textit{J. Math. Sci.} \textbf{82} (1996), no.~1, 3220--3239.

\bibitem{BT98}
V.~V.~Batyrev and Y.~Tschinkel, 
Tamagawa numbers of polarized algebraic varieties, 
\textit{Ast\'erisque} no.~251 (1998), 299--340.

\bibitem{BCHM2010}
C.~Birkar, P.~Cascini, C.~D.~Hacon and J.~McKernan,
Existence of minimal models for varieties of log general type, 
\textit{J. Amer. Math. Soc.} \textbf{23} (2010), no.~2, 405--468.

\bibitem{Bosch}
S.~Bosch, W.~L\"utkebohmert and M.~Raynaud,
\textit{N\'eron Models},
Ergebnisse der Mathematik und ihrer Grenzgebiete, 3.~Folge, vol.~21,
Springer-Verlag, Berlin, 1990.

\bibitem{BLsieving}
T.~D.~Browning and D.~Loughran, 
Sieving rational points on varieties, 
\textit{Trans. Amer. Math. Soc.} \textbf{371} (2019), no.~8, 5757--5785.

\bibitem{ChambertLoirTschinkel2010}
A.~Chambert-Loir and Y.~Tschinkel, 
Igusa integrals and volume asymptotics in analytic and adelic geometry, 
\textit{Confluentes Math.} \textbf{2} (2010), no.~3, 351--429.

\bibitem{CKR}
S.~Chan, P.~Koymans and N.~Rome,
Serre's problem for multiple conics,
\textit{Preprint}, arXiv:2504.21792.

\bibitem{ColliotTheleneSkorobogatov2021}
J.-L.~Colliot-Th\'el\`ene and A.~N.~Skorobogatov,
\textit{The Brauer--Grothendieck Group},
Ergebnisse der Mathematik und ihrer Grenzgebiete, 3.~Folge,
vol.~71, Springer, Cham, 2021.

\bibitem{CoxLittleSchenck}
D.~Cox, J.~Little and H.~Schenck,
\textit{Toric Varieties},
Graduate Studies in Mathematics, vol.~124,
American Mathematical Society, Providence, RI, 2011.

\bibitem{Denef1991}
J.~Denef, 
Report on Igusa's local zeta function, 
\textit{Ast\'erisque} no.~201--203 (1991), Exp.~no.~741, 359--386 (1992).

\bibitem{FK4rank}
\'E.~Fouvry and J.~Kl\"uners, 
On the $4$-rank of class groups of quadratic number fields, 
\textit{Invent. Math.} \textbf{167} (2007), no.~3, 455--513.

\bibitem{FMT}
J.~Franke, Y.~I.~Manin and Y.~Tschinkel, 
Rational points of bounded height on Fano varieties, 
\textit{Invent. Math.} \textbf{95} (1989), no.~2, 421--435.

\bibitem{Heath-Brown}
D.~R.~Heath-Brown, 
A new form of the circle method, and its application to quadratic forms,
\textit{J. Reine Angew. Math.} \textbf{481} (1996), 149--206.

\bibitem{HuKeel2000}
Y.~Hu and S.~Keel, 
Mori dream spaces and GIT, 
\textit{Michigan Math. J.} \textbf{48} (2000), 331--348.

\bibitem{Huang}
Z.~Huang,
Quantitative strong approximation for ternary quadratic forms III,
\textit{Preprint}, arXiv:2512.13091.

\bibitem{HSS2026}
Z.~Huang, D.~Schindler and A.~Shute, 
Optimal quantitative weak approximation for projective quadrics, 
\textit{Math. Ann.} \textbf{394} (2026), Paper No.~45.

\bibitem{Igusa1974}
J.~Igusa, 
Complex powers and asymptotic expansions. I. Functions of certain types, 
\textit{J. Reine Angew. Math.} \textbf{268/269} (1974), 110--130.

\bibitem{Lazarsfeld}
R.~Lazarsfeld,
\textit{Positivity in algebraic geometry I. Classical setting: line bundles and linear series},
Ergebnisse der Mathematik und ihrer Grenzgebiete, 3.~Folge, vol.~48,
Springer-Verlag, Berlin, 2004.

\bibitem{LehmannRiedlTanimoto}
B.~Lehmann, E.~Riedl and S.~Tanimoto,
Non-free sections of Fano fibrations,
\textit{Mem. Amer. Math. Soc.} \textbf{319} (2026), no.~1626, 109~pp.

\bibitem{almost}
B.~Lehmann and S.~Tanimoto,
A counterexample to Manin's conjecture for almost Fano varieties,
\textit{Preprint}, arXiv:2611.08890.

\bibitem{Loughran2018Varieties}
D.~Loughran,
The number of varieties in a family which contain a rational point,
\textit{J. Eur. Math. Soc. (JEMS)} \textbf{20} (2018), no.~10, 2539--2588.

\bibitem{LRS}
D.~Loughran, N.~Rome and E.~Sofos,
The leading constant for rational points in families,
\textit{Preprint}, arXiv:2210.13559.

\bibitem{LS}
D.~Loughran and T.~Santens,
Malle's conjecture and Brauer groups of stacks,
\textit{Preprint}, arXiv:2412.04196.

\bibitem{doublesolid}
F.~Okamura, 
Moduli spaces of rational curves on Artin--Mumford double solids,
\textit{Math. Z.} \textbf{313} (2026), no.~4, Paper No.~75, 18~pp.

\bibitem{Peyre}
E.~Peyre, 
Hauteurs et mesures de Tamagawa sur les vari\'et\'es de Fano, 
\textit{Duke Math. J.} \textbf{79} (1995), no.~1, 101--218.

\bibitem{Peyre2003}
E.~Peyre,
Points de hauteur born\'ee, topologie ad\'elique et mesures de Tamagawa,
\textit{J. Th\'eor. Nombres Bordeaux} \textbf{15} (2003), no.~1,
319--349.

\bibitem{PeyreTschinkel}
E.~Peyre and Y.~Tschinkel, 
Tamagawa numbers of diagonal cubic surfaces, numerical evidence, 
\textit{Math. Comp.} \textbf{70} (2001), no.~233, 367--387.

\bibitem{Salberger}
P.~Salberger, 
Tamagawa measures on universal torsors and points of bounded height on Fano varieties, 
\textit{Ast\'erisque} no.~251 (1998), 91--258.

\bibitem{Sofos}
E.~Sofos, 
Uniformly counting rational points on conics, 
\textit{Acta Arith.} \textbf{166} (2014), no.~1, 1--14.

\bibitem{Stacks}
The Stacks Project Authors, 
\textit{Stacks Project}, 
\url{https://stacks.math.columbia.edu}, 2026.

\end{thebibliography}
\end{document}